\documentclass[12pt,leqno]{amsart}
\usepackage{amssymb,amsmath,amsthm}
\usepackage{graphicx}
\usepackage{tikz}
\usetikzlibrary{arrows.meta,positioning}
\usepackage{color}
\usepackage{placeins}
\usepackage[bookmarksopen,bookmarksdepth=3,colorlinks,citecolor=red,pagebackref,hypertexnames=true]{hyperref}
\usepackage[msc-links,nobysame,sorted-cites,initials]{amsrefs}
\usepackage[inline]{enumitem}
\usepackage{mathtools}
\mathtoolsset{showonlyrefs}
\newcommand{\biburl}[1]{{\def\UrlFont{\normalfont\footnotesize}\url{#1}}}

\usepackage[lining]{libertine}
\usepackage{cabin}
\usepackage[libertine]{newtxmath}

\def\GL{{\mathrm{GL}}}
 \def\Pd{{\mathbb{R}}\mathrm{P}^d}
 \def\dpr{\mathsf d}
\def\ord{\mathrm{ord}}
\def\e{\mathrm{e}}
\def\eps{\varepsilon}
\def\d{{\mathrm d}}
\def\dim{\operatorname{dim}}
\def\R {\mathbb{R}}
\def\C {{\mathbb C}}
\def\A {{ \mathrm A}}
\def\Q {{\mathcal Q}}
\def\M {{\mathrm M}}
\def\Z {{\mathbb Z}}
\def\ind{\mathbf{1}}
\def\Span{\mathrm{span}}
\def\Gr{\mathrm{Gr}}
\newcommand{\diam}{\mathrm{diam}}
\newcommand{\supp}{\operatorname{supp}}

 \def\Xint#1{\mathchoice
     {\XXint\displaystyle\textstyle{#1}}%
     {\XXint\textstyle\scriptstyle{#1}}%
     {\XXint\scriptstyle\scriptscriptstyle{#1}}%
     {\XXint\scriptscriptstyle\scriptscriptstyle{#1}}%
     \!\int}
   \def\XXint#1#2#3{{\setbox0=\hbox{$#1{#2#3}{\int}$}
       \vcenter{\hbox{$#2#3$}}\kern-.5\wd0}}
   \def\avgint{\Xint-}

\makeatletter
\DeclareRobustCommand\widecheck[1]{{\mathpalette\@widecheck{#1}}}
\def\@widecheck#1#2{%
    \setbox\z@\hbox{\m@th$#1#2$}%
    \setbox\tw@\hbox{\m@th$#1%
       \widehat{%
          \vrule\@width\z@\@height\ht\z@
          \vrule\@height\z@\@width\wd\z@}$}%
    \dp\tw@-\ht\z@
    \@tempdima\ht\z@ \advance\@tempdima2\ht\tw@ \divide\@tempdima\thr@@
    \setbox\tw@\hbox{%
       \raise\@tempdima\hbox{\scalebox{1}[-1]{\lower\@tempdima\box
\tw@}}}%
    {\ooalign{\box\tw@ \cr \box\z@}}}
\makeatother

\newcounter{counter}

\numberwithin{equation}{section}

\newtheorem{theorem}[counter]{Theorem}
\newtheorem{corollary}[subsection]{Corollary}
\newtheorem{lemma}[subsection]{Lemma}
\newtheorem{proposition}[subsection]{Proposition}
\theoremstyle{definition}
\newtheorem{definition}[subsection]{Definition}
\newtheorem{example}[subsection]{Example}
\newtheorem{remark}[subsection]{Remark}

\theoremstyle{plain}

\numberwithin{figure}{section}

\begin{document}

\title[Directional maximal functions]{Directional maximal functions along algebraic curves}

\author[F.\ Di Plinio]{Francesco Di Plinio}
\address[F.\ Di Plinio]{Dipartimento di Matematica e Applicazioni ``R.\ Caccioppoli'' \\ \newline\indent Universit\`a di Napoli ``Federico II'', Via Cintia, Monte S.\ Angelo, 80126 Napoli, Italy}
\email{\href{mailto:francesco.diplinio@unina.it}{\textnormal{francesco.diplinio@unina.it}}}  
\thanks{F.\ Di Plinio is partially supported by the FRA 2022--2026 Program of University of Napoli Federico II, project ReSinAPAS --
Regularity and Singularity in Analysis, PDEs, and Applied Sciences.}

\author[I. Parissis]{Ioannis Parissis}
\address[I.\ Parissis]{Departamento de Matem\'aticas, Universidad del Pa\'is Vasco, Aptdo. 644, 48080 Bilbao, Spain and Ikerbasque, Basque Foundation for Science, Bilbao, Spain}
\email{\href{mailto:ioannis.parissis@ehu.eus}{\textnormal{ioannis.parissis@ehu.eus}}}
\thanks{I. Parissis is partially supported by grant PID2024-156267NB-I00 funded by MICIU/AEI/10.13039/501100011033 and cofunded by the European Union,
grant IT1875-26 of the Basque Government and IKERBASQUE}

\subjclass[2020]{Primary: 42B25. Secondary: 14P10, 41A05, 41A50, 42B20}
\keywords{Directional maximal operators, algebraic curves, Tchebycheff systems, projective invariance, complex cellular preparation}

\begin{abstract} We prove that the maximal operator over line segments in $N$ arbitrary directions on a real algebraic curve of degree $D$ in $\mathbb R^n$ has $L^2$ operator norm at most $C_{n,\varepsilon}D^{C_n}N^\varepsilon$ for every $\varepsilon>0$, with $C_n$ independent of $\varepsilon$. Combined with polynomial partitioning, this gives the sharp bound $N^{\frac14}$, up to a polylogarithmic factor, for $N$ arbitrary directions in $\mathbb R^3$.  The first main ingredient is a new quantitative strengthening of the property that the curve crosses each affine hyperplane at most $D$ times. This condition controls the overlap of variable-width neighborhoods of hyperplanes. The width attached to each direction is a small multiple of its distance to the osculating hyperplane at an endpoint of its arc. The problem is projectively invariant, and a projective almost orthogonality principle turns the overlap control into the $L^2$ bound. The second main ingredient, building on the complex cellular preparation of Binyamini and Novikov, is a new effective decomposition of algebraic curves. Up to $D^{C_n}$ points, the curve is covered by $D^{C_n}$ arcs whose coordinates, in a suitable parametrization and projective frame, are monomials with distinct exponents times nearly constant units. In these charts we verify the quantitative crossing condition with polynomial dependence on $D$. The framework also applies to nonalgebraic curves with doubling Frenet curvature measures. Consequences include local Sobolev maximal estimates along measurable vector fields with values in an algebraic curve, and the same bounds, up to a logarithmic factor, for directional Calder\'on--Zygmund operators.
\end{abstract}
\maketitle

\section{Introduction} Our topic is maximal averages over straight line segments whose directions vary along a curve. More precisely, for a finite set $V\subset\R^n\setminus\{0\}$, let
\[
\M_V f(x)\coloneqq\sup_{\substack{v\in V\\r>0}} \frac1{2r}\int_{-r}^r|f(x-tv)|\,\d t,\qquad x\in \R^n.
\]
For finite sets $V$ with $\#V\leq N$, lying on a real algebraic curve in $\R^n$ of degree at most $D$, we prove an $L^2$ bound with arbitrarily small power growth in $N$. The estimate simultaneously allows all averaging lengths, arbitrary finite subsets of the curve, and singular or degenerate algebraic loci. The curve enters only through a multiplicative constant, polynomial in $D$ with exponent depending only on the dimension.

A basic feature of the problem is projective invariance. The supremum over lengths absorbs the size of each direction vector, and an invertible linear change of the spatial variables preserves the operator norm. Thus, if $T\in\GL_n(\R)$, then
\[
\left\|\M_{ T(V)}\right\|_{2\to2}=\|\M_V\|_{2\to2}.
\]
Consequently, the natural setting for studying the operator norm of $\M_V$ is real projective space, in the sense that $V$ is regarded as a subset of a curve in $\mathbb R\mathrm P^{n-1}$. This symmetry, absent from the corresponding problem for the single-scale directional average, means that an operator-norm bound for $\M_V$ should be formulated for a class of curves which is invariant under projective transformations. Algebraic curves form such a class, and are particularly suited for the statement of our main result. Moreover, maximal operators along algebraic curves arise naturally in the corresponding problem for \emph{arbitrary directions} in $\R^n$, through polynomial partitioning, as explained below. We also give a projectively invariant geometric condition under which a bound of the same form holds for nonalgebraic curves; see Subsection~\ref{subsec:intro-doubling}.

\subsection{Principal result} Hereafter, an equivalent smooth Fourier model of $\M_V$ is used. Fix a real even $\varphi\in \mathcal C_c^\infty(\R)$ with $\supp\varphi\subset[-1,1]$, whose inverse Fourier transform $\psi\coloneqq\widecheck\varphi$ satisfies
\begin{equation}\label{eq:phi-kernel}
 \psi\geq0\quad\hbox{on }\R,\qquad \psi\geq c_\varphi>0\quad\hbox{on }[-1,1].
\end{equation}
Define for $v\in\R^n$
\[
\widehat{\A_{v,r}f}(\xi)=\varphi(r\xi\cdot v)\widehat f(\xi),\qquad \xi\in \R^n, \qquad \A_Vf=\sup_{\substack{v\in V\\r>0}}|\A_{v,r}f|.
\]
Positivity, the lower bound in \eqref{eq:phi-kernel}, and the rapid decay of $\psi$ give $\M_Vf\lesssim_\varphi\A_V|f|$ and $\A_V f\lesssim_\varphi\M_Vf$. The notation $\A_V$ is also meaningful when $0\in V$, since $\A_{0,r}f=\varphi(0)f$. Constants may depend on the fixed function $\varphi$; in this introduction we suppress this dependence.

\begin{theorem}\label{thm:intro-main}
Let $n\geq2$. There is an exponent $C_n<\infty$, depending only on $n$,
such that for every $\eps>0$ there is $C_{n,\eps}<\infty$ with the
following property. If $Z\subset\R^n$ is a real algebraic set of dimension
at most one and degree at most $D\geq1$, then
\begin{equation}\label{eq:intro-main}
 \sup_{\substack{V\subset Z \\\#V\leq N}} \|\A_V\|_{L^2(\R^n)\to L^2(\R^n)} \leq C_{n,\eps}D^{C_n}N^\eps,\qquad N\geq1.
\end{equation}
Degree is understood as in Section~\ref{sec:algebraic}, counting all irreducible components, including isolated points.
\end{theorem}

Theorem \ref{thm:intro-main}, to our knowledge, is the first subpower $L^2$ bound for positive multiscale averages uniform over arbitrary finite subsets of general algebraic direction curves in dimensions $n\geq3$. In fact, the proof gives the stronger estimate
\begin{equation}\label{eq:intro-subpower}
 \sup_{\substack{V\subset Z\\\#V\leq N}}\|\A_V\|_{L^2(\R^n)\to L^2(\R^n)} \leq C_nD^{C_n} \exp \left(c_n\sqrt{\log(\e+N)}\right).
\end{equation}
Here one may take $c_n=1+\sum_{j=2}^{n-1}\sqrt{\log j}$. This follows from the quantitative form of Theorem~\ref{thm:projective-SECT-maximal} below. We retain \eqref{eq:intro-main} as the main formulation because it expresses the cardinality conclusion simply. The degree exponent in both estimates is independent of $\eps$ and of the coefficients defining the curve. This uniformity includes curves approaching degenerate configurations.

Maximal operators along sets of directions have been studied since the 1970s, and we refer to the survey \cite{DPPsurvey} for their history and open problems. In the plane, by the combined results of \cites{CordobaFefferman1977,Stromberg1977,NSW1978,SjogrenSjolin1981}, maximal operators along lacunary sets of directions of finite order are bounded on $L^p$ for $1<p<\infty$. The necessity of lacunarity was studied by Bateman \cite{Bateman2009}; see also \cite{HRMS}, and almost orthogonality principles for planar directional maximal operators were introduced by Alfonseca, Soria and Vargas \cites{ASV2003,ASV2003b}. In higher dimensions, special lacunary sets of directions were treated by Nagel, Stein and Wainger \cite{NSW1978} and by Carbery \cite{Carbery1988}, and general ones by Parcet and Rogers \cite{PR}.

Theorem \ref{thm:intro-main} is motivated by two classical results. Both concern curves that meet each hyperplane in at most $D$ points, for some $D$; we say that such a curve has the \emph{$D$-crossing property}, or the \emph{finite-crossing property} when $D$ is not specified. This is a generalized form of convexity; see \cite{BaMaPo}. C\'ordoba \cite{Cordoba1982}*{Theorem 2} proved the $L^2$ bound $C_\gamma(\log N)^2$ for the maximal function along tubes of fixed eccentricity $N$, of all lengths, pointing in the directions of a smooth curve $\gamma$ with the finite-crossing property. The single-scale step of his argument uses only this property, and is written for the dyadic samples $\gamma(j2^{-m})$ \cite{Cordoba1982}*{pp. 222--223}. It was made uniform in \cite{DPPAJM}*{Proposition 5.4 and Theorem D}, which give single-scale bounds over arbitrary $N$-point subsets of the curve, with constants depending only on the number of crossings, and only on the degree for smooth algebraic curves on the sphere. On the other hand, Christ, Duoandikoetxea and Rubio de Francia \cite{CDR1986}*{Theorem 2} established a logarithmic bound for the directional Hilbert transform over arbitrary finite subsets of a rectifiable curve with the finite-crossing property. They also expected the analogue of this theorem for the directional maximal operator, and observed that their method gives it for single-scale averages \cite{CDR1986}*{p. 202}. Theorem~\ref{thm:intro-main} provides this analogue for algebraic curves, with the degree playing the role of $D$ in the $D$-crossing property and with the growth \eqref{eq:intro-subpower}, slower than any power of $N$. As explained below, arguments based on crossing counts alone do not extend to averages over segments of arbitrary lengths.

The common geometric input of the contributions \cites{CDR1986,Cordoba1982} described above is the finite-crossing property. For the Hilbert transform, the multiplier is constant on each side of $\xi\cdot v=0$, so, by the finite-crossing property, it has a bounded number of jumps along an ordered curve; since the symbol of the Hilbert transform takes only two values, its total variation is bounded. At a fixed averaging scale, the analogous argument uses band indicators, or superpositions of such indicators, whose boundaries are again hyperplanes. In both cases the relevant functions of $v$ are piecewise constant with jumps on hyperplanes, and this turns crossing counts into square-function estimates. The same mechanism is behind the single-scale curve bounds in \cite{DPPAJM} and the Hilbert-transform estimates of Kim and Pramanik \cite{KimPramanik2022}, which hold in every dimension along such curves, in particular along algebraic curves of bounded degree. Combined with polynomial partitioning, the latter give the bound $N^{1/4}h(N)$, with $h$ an iterated logarithm, for arbitrary directions in $\R^3$, and $N^{\frac{n-2}{2(n-1)}+\eps}$ in $\R^n$, $n\geq3$ \cite{KimPramanik2022}*{Theorems 1.1, 1.2, 2.1 and 3.5}.

For positive averages over all lengths, the geometric question changes. One must control approximate hyperplane incidences at many scales, with the relevant thickness varying along the curve. A count of exact crossings does not supply this control. The new input is a quantitative form of the generalized convexity property described above. A strong form of convexity of a curve in $\R^d$ is that its coordinates, together with the constant function $1$, form a \emph{complete Tchebycheff system}, a classical notion from approximation theory: for each $k\leq d$, every nontrivial linear combination of the first $k+1$ of these functions has at most $k$ zeros. Consequently, no affine hyperplane meets such a curve in more than $d$ points. We show that this crossing property is stable: no hyperplane can pass close to the curve at $d+1$ points taken from separate consecutive arcs, where closeness is measured against the distance from the curve to its osculating hyperplanes at the endpoints of these arcs. The \emph{switch moduli} introduced below measure precisely this stability. Combined with an almost orthogonality principle in the projective setting, a bound on the switch moduli yields the $L^2$ estimate for $\M_V$, uniformly over all subsets $V$ of the curve with $\#V\leq N$, however unevenly they are distributed.

Proving uniform operator norm bounds over arbitrary finite subsets of the curve is essential, and estimates for a fixed sampling of the curve do not provide them. For a fixed smooth finite-type parametrization $\gamma$, Barrionuevo \cite{Barrionuevo1993} states a logarithmic multiscale bound for the dyadic samples $\gamma(j2^{-m})$, with a constant depending on $\gamma$. Finite type means that the derivatives of $\gamma$, including order zero, span the ambient space at every parameter value. This condition is invariant under smooth changes of parameter, though the dyadic sampling condition is not. Reparametrizing to include a prescribed finite set in a dyadic sample changes the parametrization-dependent constant, with no stated control in terms of the cardinality. Thus that result does not furnish an estimate uniform over arbitrary $N$-point subsets of a fixed image curve.\footnote{There is also a gap in \cite{Barrionuevo1993}*{Lemma~1}, which asserts that $|\xi\cdot(\gamma^{(j)}(s)-\gamma^{(j)}(t))|>C_j|s-t|^{L-j}$ whenever $\xi$ is a unit vector, $|s-t|<\delta_j$ and $|\xi\cdot\gamma^{(\ell)}(t)|<c2^{-\ell-1}$ for all $\ell<j$, where $L$, $c$, $\delta_j$ and $C_j$ are positive constants depending only on $\gamma$ and $j$. For $j=1$ this fails for the spherical curve $\gamma(t)=(1,u,u^2)/(1+u^2+u^4)^{1/2}$, $u=t-\tfrac12$, which meets each plane at most four times and is of finite type, as $\det(\gamma,\gamma',\gamma'')=2(1+u^2+u^4)^{-3/2}$: take $\xi=(0,1,0)$, $t=\tfrac12-h$, $s=\tfrac12+h$ and $h$ small. Then $\xi\cdot\gamma$ is odd in $u$, hence small at $t$, while $\xi\cdot\gamma'$ is even, so the left-hand side vanishes.} Theorem~\ref{thm:intro-main} provides exactly this uniformity, together with polynomial control by algebraic degree.

Polynomial partitioning was introduced by Guth and Katz \cite{GuthKatz2015} for the Erd\H{o}s distinct distances problem, and extended to families of varieties by Guth \cite{GuthVarieties2015} with incidence problems in view. It has since become a central tool in Fourier restriction \cites{GuthJAMS2016,GuthActa2018} and in Kakeya-type problems \cite{GuthZahl2018}. In these works the partition takes place in the ambient space of the problem in question. Its use for directional problems through a partition of the set of directions, organized by algebraic dimension and degree, was introduced in \cite{DPPAJM}, where it led to single-scale maximal estimates; it was subsequently used for the almost orthogonality principle for subspace averages in \cite{DP-PAdv}. Here algebraic geometry has a further role. Decompositions of polynomial curves adapted to their torsion \cite{DendrinosWright2010} and preparation theorems for subanalytic functions \cite{CCMRS2018} have been used in Fourier restriction and in the study of oscillatory integrals, respectively. The new ingredient here, which to our knowledge has not appeared before in harmonic analysis, is \emph{effective} preparation, with complexity polynomial in the degree. Namely, effective cellular preparation writes the coordinate functions of the curve, after splitting into polynomially many pieces, as products of monomials and essentially constant real analytic units. The complex cellular theory of Binyamini and Novikov \cite{BN21} supplies this preparation. It originates in the smooth parametrization problem, where it gives the Yomdin--Gromov algebraic lemma \cites{Yom87a,Yom87,Gr87,Burguet} with polynomially many charts and sharpens the parametrization theorem of Pila and Wilkie \cite{PilaWilkie2006} for subanalytic sets. A projective transformation then brings each prepared piece to a representative position in which the switch moduli are bounded polynomially in the degree. The pieces stop at the singular points of the curve and at the branch points of a coordinate projection; since the weights of the canonical charts constructed in Section~\ref{sec:algebraic} are positive on each piece, the pieces also avoid the inflection points of the curve, where the osculating flag degenerates. The preparation accounts as well for such points of the complexified curve lying near, but not on, the real curve. The number of prepared pieces is polynomial in the degree, so the resulting estimate for the curve is uniform, even for sets $V$ with points arbitrarily close to such degenerate points of the complexified curve.

The curve estimate feeds back into polynomial partitioning: for arbitrary directions in three dimensions, the partition splits the directions into those in the cells, which are handled by induction on the cardinality, and those lying on the zero set of the partitioning polynomial, an algebraic curve in the projective plane $\mathbb R\mathrm P^2$ of directions; Theorem~\ref{thm:intro-main} bounds the contribution of the latter at all scales, which leads to Corollary~\ref{cor:intro-arbitrary-3d} below.

\subsection{Applications} As described above, the curve estimate is the central input for cardinality bounds for arbitrary directions in three dimensions. In the plane, the multiscale problem for $N$ arbitrary directions was settled by Katz \cite{KatzDuke}, who proved the sharp logarithmic bound in the number of directions, extending the case of uniformly distributed directions \cite{Stromberg1978}; see \cite{KatzBLMS} for single-scale averages. In three dimensions, sets of $N$ directions can have operator norm at least of order $N^{1/4}$, while upper bounds of this order for maximal averages were previously known only at a single scale or for special sets of directions, such as lacunary sets \cite{PR}. Demeter \cite{Demeter2012} proved the bound $N^{1/4}$, up to logarithmic factors, for single-scale averages along $N$ arbitrary directions in $\R^3$; see also \cite{DPPAJM}*{Theorem B} and \cite{DP-PAdv}*{Theorem F} for quantitative improvements on the lower order term of the operator norm. Jongchon Kim \cite{Kim2021}*{Theorem 1.1} proved multiscale bounds without logarithmic factors for equidistributed directions in all dimensions $n\geq3$, that is, for maximal $\delta$-separated subsets of the sphere, both for maximal averages and for maximal directional singular integrals; in three dimensions his bound is $CN^{1/4}$. Combined with polynomial partitioning, Theorem~\ref{thm:intro-main} removes these restrictions, at the cost of a polylogarithmic factor.

\begin{corollary}\label{cor:intro-arbitrary-3d} There are absolute constants $B,C<\infty$ such that
\begin{equation}\label{eq:intro-arbitrary-3d}
 \sup_{\substack{V\subset \R^3\setminus\{0\} \\\#V\leq N}} \|\A_V\|_{L^2(\R^3)\to L^2(\R^3)} \leq C N^{\frac14}\left(\log(2+N)\right)^B.
\end{equation}
\end{corollary}

As mentioned above, the cardinality exponent $1/4$ is optimal, as shown by equidistributed direction sets \cite{Demeter2012}. The argument combines the case $n=3$ of Theorem~\ref{thm:intro-main} with the algebraic almost orthogonality principle of \cite{DP-PAdv}*{Theorem E}. Keeping track of the polynomial degree dependence yields the polylogarithmic factor in  \eqref{eq:intro-arbitrary-3d}. The proof, a recursion in the number of directions, is given in Section~\ref{sec:applications}.
For $n\geq4$, polynomial partitioning reduces the problem to directions on algebraic varieties of dimension $n-2\geq2$; extending the analysis of this paper from curves to such varieties is a natural question, which we discuss in Remark~\ref{rmrk:higherd}.

An independent motivation comes from differentiation of Sobolev functions along variable directions. Zygmund's problem asks how regularity of a vector field can ensure differentiation by averages over line segments pointing in the directions of the field. A complementary question imposes regularity on the function while allowing a merely measurable field. Bakas,   Roncal and the authors \cite{BDPR2024} obtained differentiation results for Besov functions along measurable choices of codimension-one subspaces. Here the multiscale curve estimate gives local maximal control with arbitrarily small positive Sobolev regularity, for every measurable vector field with range contained in an algebraic curve. No continuity of the vector field is assumed, and the bound is uniform over all such fields.

To make this precise, let $\Gamma\subset {{\mathbb S}^{n-1}}$ be a real algebraic set of dimension at most one and degree at most $D$; the normalization to unit vectors is needed here, since the averaging lengths are now bounded and no longer absorb the length of $v$. Define the local maximal operator
\begin{equation}\label{eq:intro-local-maximal}
 \M_\Gamma^{\mathrm{loc}}f(x) \coloneqq\sup_{\substack{v\in\Gamma\\0<r\leq1}} \frac1{2r}\int_{-r}^r|f(x-tv)|\,\d t,\qquad x\in\R^n.
\end{equation}

\begin{corollary}\label{cor:intro-sobolev} For every $s>0$,
\begin{equation}\label{eq:intro-sobolev}
\|\M_\Gamma^{\mathrm{loc}}f\|_{L^2(\R^n)} \leq C_{n,s}D^{C_n}\|f\|_{H^s(\R^n)}, \qquad f\in\mathcal S(\R^n).
\end{equation}
Consequently, for every measurable vector field $v:\R^n\to\Gamma$, the operator
\[
 \M_{v(\cdot)} ^{\mathrm{loc}}f(x) \coloneqq \sup_{0<r\leq1}\frac1{2r}\int_{-r}^r|f(x-tv(x))|\,\d t,\qquad x\in\R^n,
\]
extends boundedly from $H^s(\R^n)$ to $L^2(\R^n)$, with the same constant, independent of the vector field.
\end{corollary}

The maximal estimate for $\M_{v(\cdot)}^{\mathrm{loc}}$, combined with approximation of $H^s$-functions by smooth functions, yields almost-everywhere differentiation for the averages of $H^s$-functions, for any $s>0$, along directions given by $v(\cdot)$. Although irrelevant for a differentiation result, we note that the restriction to lengths $r\leq1$ cannot be removed. Indeed, by dilation to low frequencies, an $H^s\to L^2$ estimate over all lengths would imply the corresponding $L^2\to L^2$ estimate, which fails as soon as $\Gamma$ contains an arc.  Both corollaries above are proved in Section~\ref{sec:applications}.

The positive maximal estimates also apply to directional singular integrals through an exponential good-$\lambda$ comparison. Let $T$ be a fixed linear Calder\'on--Zygmund operator on $\R$, bounded on $L^2(\R)$, whose kernel satisfies the standard size and H\"older continuity estimates. For $v\in\R^n\setminus\{0\}$, define its action on each parallel line by
\begin{equation}\label{eq:intro-directional-CZ}
 T_vf(y+tv)=T\left[f(y+\cdot v)\right](t), \qquad y\in v^\perp,\quad t\in\R, \qquad T_Vf=\sup_{v\in V}|T_vf|.
\end{equation}
This convention fixes both the origin and the scale of the coordinate on each line. If $T$ is a convolution operator with multiplier $m$, then $T_v$ has multiplier $m(\xi\cdot v)$; in particular, the results below apply to multipliers satisfying the Mikhlin condition $|m^{(k)}(\xi)|\lesssim|\xi|^{-k}$ for $k\leq2$, whose kernels satisfy the standard estimates. The constants may depend on the $L^2$ norm of $T$, its kernel bounds, and its H\"older exponent.

\begin{corollary}  \label{cor:intro-CZ-algebraic} Let $T$ be a Calder\'on--Zygmund operator as above, $n\geq2$, and let $Z\subset\R^n$ be a real algebraic set of dimension at most one and degree at most $D\geq1$. For every $\eps>0$,
\begin{equation}\label{eq:intro-CZ-algebraic}
  \sup_{\substack{V\subset Z\setminus\{0\}\\\#V\leq N}}\|T_V\|_{L^2(\R^n)\to L^2(\R^n)} \leq C_{n,T,\eps}D^{C_n}N^\eps.
\end{equation}
The exponent $C_n$ is independent of $T$ and $\eps$. More precisely, the right-hand side may be replaced by $C_{n,T}D^{C_n}\exp\left(C_n\sqrt{\log(\e+N)}\right)$.
\end{corollary}

\begin{corollary} \label{cor:intro-CZ-3d} For every $\eps>0$,
\begin{equation}\label{eq:intro-CZ-3d}
 \sup_{\substack{V\subset \R^3\setminus\{0\}\\\#V\leq N}}\|T_V\|_{L^2(\R^3)\to L^2(\R^3)} \leq C_{T,\eps}N^{1/4+\eps}.
\end{equation}
In fact, with $B$ as in Corollary~\ref{cor:intro-arbitrary-3d},
\[
\sup_{\substack{V\subset \R^3\setminus\{0\}\\\#V\leq N}} \|T_V\|_{L^2(\R^3)\to L^2(\R^3)}\leq C_T N^{1/4}\left(\log(2+N)\right)^{B+1}.
\]
\end{corollary}

The comparison of the operator norm of $T_V$ with that of $\M_V$ loses only a factor of $\log(2+N)$. This follows from Hunt's exponential good-$\lambda$ method \cite{Hunt1972}, applied on each line as in Demeter and Di Plinio \cite{DemeterDiPlinio2014}*{Proposition 2.2}. For multipliers $m_1,\ldots,m_N$ with $\sup_\xi\sum_{j=1}^N|m_j(\xi)-m_{j-1}(\xi)|<\infty$, where $m_0=0$, Christ, Duoandikoetxea and Rubio de Francia \cite{CDR1986}*{Lemma 5.1} proved an $L^2$ bound of order $\log N$ for the corresponding maximal operator. Along a curve with the finite-crossing property, in any dimension, the directional Hilbert transforms, ordered along the curve, form such a family, which gives their curve estimate. As discussed above, the combination of finite crossing and bounded variation is very special: for a general translation invariant Calder\'on--Zygmund operator, the symbols $m(\xi\cdot v)$ need not have bounded variation along the curve, and this argument is not available. The positive multiscale theorem, through the comparison above, makes the curve estimates, and the three-dimensional estimate for arbitrary directions, available for general Calder\'on--Zygmund operators. The comparison and the proofs are given in Subsection~\ref{subsec:directional-CZ}.

\subsection{A geometric condition beyond algebraicity} \label{subsec:intro-doubling}  Recall that the switch moduli quantify the stability of the crossing property, and that a bound on them, combined with almost orthogonality, yields the $L^2$ estimate. Beyond algebraic curves, the switch moduli can be bounded by geometric conditions. For a planar convex arc, a doubling condition on the measure recording the turning of its tangent suffices; see Remark~\ref{rem:planar-doubling-switch}, which uses the comparable measure $F''(t)\,\d t$ of a normalized graph chart; by Lemma~\ref{lem:iterated-doubling}, the doubling constants of the two measures are polynomially related. Section~\ref{sec:doubling-curvatures} extends this to higher dimensions, with the Frenet curvatures in place of the turning of the tangent.

Working in $\R^n$, $n\geq3$, we write $d=n-1$. Let $\gamma:[0,\ell]\to\mathbb S^d$ be parametrized by arclength, continuous on the closed interval and ${{\mathcal C}^d}$ in its interior, with $\gamma,\gamma',\ldots,\gamma^{(d)}$ linearly independent there. Write $\kappa_1,\ldots,\kappa_{d-1}$ for its spherical Frenet curvatures: if $E_0,\ldots,E_d$ is the orthonormal frame obtained from $\gamma,\gamma',\ldots,\gamma^{(d)}$ by the Gram--Schmidt process, so that $E_0=\gamma$ and $E_1=\gamma'$, then $E_j'=-\kappa_{j-1}E_{j-1}+\kappa_jE_{j+1}$ for $1\leq j\leq d-1$, with $\kappa_0\coloneqq1$; by the linear independence, $\kappa_1,\ldots,\kappa_{d-1}$ are positive in the interior. Suppose that the turning measures $\kappa_i(s)\,\d s$, which we assume to be finite, satisfy
\begin{equation}\label{eq:intro-doubling}
\int_J\kappa_i(s) \,\d s \leq K\min\left\{\int_{J^-} \kappa_i(s) \,\d s,\int_{J^+}\kappa_i(s)\,\d s\right\}
\end{equation}
for every interval $J\subset[0,\ell]$, where $J^-$ and $J^+$ are its two halves. Theorem~\ref{thm:doubling-curvatures} gives
\[
\sup_{\substack{V\subset\gamma([0,\ell])\\\#V\leq N}}\|\A_V\|_{2\to2} \leq C_{d,\eps}K^{C_d} \left(1+\ell+\sum_{i=1}^{d-1}\int_0^\ell \kappa_i(s) \,\d s\right)N^\eps.
\]

The left-hand side is projectively invariant, while the Frenet curvatures on the right-hand side are not. For $L\in\GL_{d+1}(\R)$, let $\gamma_L:[0,\ell_L]\to\mathbb S^d$ be the arclength parametrization of $L\gamma/|L\gamma|$, let $\kappa_{1,L},\ldots,\kappa_{d-1,L}$ be its Frenet curvatures, and let $K_L\in [2,\infty]$ be the least constant in \eqref{eq:intro-doubling} for its turning measures.

\begin{corollary}\label{cor:intro-doubling-projective}
For every $\eps>0$,
\[
\sup_{\substack{[V]\subset[\gamma([0,\ell])]\\\#V\leq N}}\|\A_V\|_{2\to2} \leq C_{d,\eps}\inf_{L\in\GL_{d+1}(\R)}K_L^{C_d} \left(1+\ell_L+\sum_{i=1}^{d-1}\int_0^{\ell_L} \kappa_{i,L}(s) \,\d s\right)N^\eps,
\]
where $[V]$ and $[\gamma([0,\ell])]\subset\mathbb R\mathrm P^d$ denote projectivizations.
\end{corollary}

This follows by applying Theorem~\ref{thm:doubling-curvatures} to each $\gamma_L$. By Proposition~\ref{prop:projective-curvature-doubling}, an invertible linear map followed by radial normalization multiplies the turning measures by bounded positive densities and changes arclength bi-Lipschitzly. Hence the finiteness of the turning measures and of the doubling constant does not depend on $L$, although the value of $K_L$ does.

For arcs which are ${{\mathcal C}^d}$ and nondegenerate on the closed interval $[0,\ell]$, the curvatures are bounded above and away from zero, and \eqref{eq:intro-doubling} holds trivially. However, doubling of the turning measures allows much more: the curvatures may vanish or become unbounded at the endpoints, for instance at power rates, and such behavior enters the estimate only through the doubling constant and the total masses of the turning measures. By subadditivity, finite unions of such arcs are covered as well. This includes arcs contained in lower-dimensional linear subspaces, to which the criterion applies within their span, with the Frenet curvatures computed there, while arcs of great circles are covered by the case $d=1$ of Theorem~\ref{thm:projective-SECT-maximal}; see Lemma~\ref{lem:hyperplane-reduction}. Section~\ref{sec:doubling-curvatures} contains explicit nonalgebraic examples.

For algebraic curves, where the dependence on the degree should be polynomial, we instead bound the switch moduli directly, in suitable projective frames; Section~\ref{sec:doubling-curvatures} explains why a direct use of the doubling constant would lose that polynomial control.

\subsection{The mechanism of the proof} \label{subsec:intro-mechanism} The proof seeks lower-dimensional approximations to the curve which are accurate enough for an almost orthogonality argument. Partition the directions into consecutive blocks and, on each block, replace them by directions on the best approximating curve inside the osculating hyperplane at one of the endpoints of the block. There are two quantities to compare: the error of this replacement and the width of the frequency bands for which overlap can be controlled. For a direction $v$ and its replacement $v'$, let $\delta$ be the angle between them, and fix a small $\eps>0$. Split each average $\A_{v,r}$ into its parts at frequencies $|\xi|$ below and above $(r\eps\delta)^{-1}$. Since the multiplier $\varphi(r\xi\cdot v)$ is supported in $|\xi\cdot v|\leq r^{-1}$, the high-frequency part lives in the band $|\xi\cdot v|\leq\eps\delta|\xi|$ around the hyperplane $v^\perp$. The low-frequency part is dominated by the average in the direction $v'$, at the cost of a factor $\eps^{-d}$; this is the switch lemma of Section~\ref{sec:AO}. Thus a small $\eps$ narrows the bands, which helps to control their overlap across blocks, and its cost falls only on the lower-dimensional comparison.

The strategy outlined above grew out of our analysis of the simplest planar case, the circle, whose constant curvature makes the argument particularly transparent; the same computation applies to short convex arcs with approximately constant curvature $\kappa$. In such a case, the osculating hyperplanes at the endpoints of a block are the tangent lines there, and these are the lower-dimensional replacement curves. For a point at arclength distance $\ell$ from the nearer endpoint of its block, replacing it by the corresponding point of the tangent line there incurs an error of order $\kappa\ell^2$. The critical bands are those nearly tangent to the arc: transversal bands meet it in much shorter pieces, which are controlled by counting crossings. A nearly tangent band of width $\eps\kappa\ell^2$ containing this point meets the arc near it only in a piece of length of order $\sqrt{\eps} \ell$, which for small $\eps$ does not reach the endpoint; hence such bands do not see the neighboring blocks. The argument is rather straightforward in this case, exactly because the curvature cancels: the same scale controls both the approximation error and the separation needed for almost orthogonality. This requires the curvature to be roughly constant at the relevant scales. When the curvature varies rapidly, the two scales no longer match, and this is a serious obstruction to generalizing the argument.

Moving beyond the model case described above, the curvature heuristic has to be replaced by a direct comparison: the error of replacing a point by the tangent line at the nearer endpoint of its block, against the distance of the point from secants through points of other blocks. Uniform algebraic estimates require us to retain this comparison while allowing a separate projective frame on each piece. Approximate constancy of Euclidean curvature is unstable under projective changes of frame: a circle and a hyperbola are projectively equivalent, yet the circle has constant curvature, while the curvature of the hyperbola is largest at its vertex and tends to zero along its branches. We therefore work with these approximation errors themselves rather than with curvature. For example, if the curve is a planar convex graph and one point is chosen in each of three blocks, these errors are the deviations of one of these points from an endpoint tangent of its block and from the secant through the other two points. When the curvature is approximately constant, both deviations are proportional to it, so the curvature cancels from the ratio of the errors, which is then bounded independently of the positions of the points.

Approximation theory gives a useful formulation of this comparison.
Choose a canonical curve $\mathbf X=(X_0,\ldots,X_d)$ and its nested
spaces $\mathcal U_j=\Span\{X_0,\ldots,X_j\}$. The switch modulus
$\mathsf S_j$ bounds the sum of the ratios
\[
\frac{\text{error of endpoint osculating approximation}}{\text{error of interpolation at the other nodes}}
\]
over ordered configurations of $j+1$ nodes in separate blocks. With one node in each of $j+1$ blocks, interpolation at the other $j$ nodes replaces the secant, and osculating approximation at an endpoint of the block replaces the tangent, so these are the two deviations of the previous paragraph. Both errors are measured in the same chart and normalization. This is the residual formula of Proposition~\ref{prop:switch-modulus-product-form}, and the bound is required even when the nodes cluster.

To see how the comparison controls overlap, suppose a frequency $\xi$ belongs to bands from $d+1$ blocks. The incidence function $h_\xi=\xi\cdot\mathbf X$ is then small at the corresponding nodes, relative to their osculating errors. Reconstructing $h_\xi$ from these values by nodal interpolation divides by the interpolation errors in the denominators of the ratios above, so the osculating errors in the numerators combine with them into the switch modulus. In the canonical normalization, for a unit coefficient vector $\xi$, this gives
\[
 1\leq\eps(d+1)\mathsf S_d.
\]
For sufficiently small $\eps$ this is impossible. Thus a bound for $\mathsf S_d$ makes the $d$-crossing property stable at the contact scales required by the analytic argument.

The nested spaces are needed to iterate the argument. The replacement curve in an osculating hyperplane is linearly equivalent to the curve formed by the preceding coordinates, so the comparison term passes from $\mathcal U_d$ to $\mathcal U_{d-1}$, one dimension lower. Restriction to a block preserves the bounds for the switch moduli. These two operations account for the full tower $\mathsf S_2,\ldots,\mathsf S_d$.

The maximal theorem of Section~\ref{sec:projective-switch-classes} requires only that, after some projective change of frame, the curve be a canonical curve with controlled switch moduli. For algebraic curves, such a chart is constructed in Section~\ref{sec:algebraic}. A projective transformation of the curve amounts to a change of basis in the space spanned by its coordinate functions, that is, in the space of incidence functions $h_\xi=\xi\cdot\mathbf X$; choosing the chart means choosing a good basis of this space. On each piece produced by cellular preparation, this basis is chosen so that its members are monomials times essentially constant units, with distinct exponents; this combines the preparation with an induction on the dimension of the space, which alternately differentiates and integrates. Ordered by exponent, such a basis makes the curve comparable to a generalized moment curve $t\mapsto(1,t,t^{\lambda_2},\ldots,t^{\lambda_d})$. Section~\ref{sec:power-unit} then bounds the switch moduli of such curves, using positive integral formulas to compare the two approximation errors uniformly over the nodes, with bounds polynomial in the exponents. Together with the polynomial complexity of preparation, this yields polynomial dependence on the degree.

\subsubsection*{Structure of the paper.} Sections~\ref{sec:AO}--\ref{sec:projective-switch-classes} establish the analytic theorem for curves with controlled switch moduli. Section~\ref{sec:power-unit} verifies this condition for charts comparable to generalized moment curves, called power-unit charts, and Section~\ref{sec:algebraic} constructs those charts for algebraic curves and completes the proof of Theorem~\ref{thm:intro-main}. Section~\ref{sec:doubling-curvatures} is independent of that proof and deals with the switch condition beyond algebraic curves. The applications follow in Section~\ref{sec:applications}.

\section{The abstract almost orthogonality principle} \label{sec:AO}
Theorem~\ref{thm:abstract-ao} takes as input a switch system on a
partitioned subset of $\Pd$. This consists of comparison directions and a bound on the
replacement distance; the overlap of the associated thick bands enters through the overlap number. It reduces
the maximal estimate to those for the comparison sets and the original
blocks, without requiring a parametrization or algebraic structure.
We first record the projective invariance and frequency localization
that make this reduction possible.


\subsection{Directions as projective points} Let $\varphi$ be the function fixed in the introduction and $\psi\coloneqq\widecheck\varphi$ its inverse Fourier transform.
For $v\in\R^{d+1}\setminus\{0\}$, $ r>0 $, define the \emph{smooth directional average}
\begin{equation} \label{eq:multiplier}
\A_{v,r}f(x) \coloneqq   \int_\R f(x-tv)\psi\left(\frac {t} {r}\right)\frac{\d t}{r}=\int_{\R^{d+1}} \varphi(r\xi\cdot v)\widehat f(\xi) \e^{2\pi i x \cdot \xi} \, \d \xi.
\end{equation}
Since $\psi\geq0$, each $\A_{v,r}$ is a positive operator,  and
\begin{equation} \label{eq:A-abs}
  \left|\A_{v,r}f\right|\leq\A_{v,r}|f| .
\end{equation}
Furthermore $\varphi(r\xi\cdot cv)=\varphi(|c|r\xi\cdot v)$ for every $c\in\R\setminus\{0\}$, whence
\begin{equation} \label{eq:homogeneity}
  \A_{cv,r}=\A_{v,|c|r}, \qquad c\neq0 .
\end{equation}
Thus the \emph{family} $\{\A_{v,r}: r>0\}$ depends on $v$ only through the line $\R v$. For this reason, we parametrize our maximal operators using projective points.

Hereafter, $\Pd=\Gr(1,d+1)$ indicates the real projective space in $\R^{d+1}$, that is, the set of all lines through the origin in $\R^{d+1}$. Below, write $[v]$ for the projectivization of $v\in \R^{d+1}\setminus \{0\}$, and set
\begin{equation} \label{eq:Aomega}
\A_{[v]} f\coloneqq\sup_{r>0}\left|\A_{v,r}f\right|,\qquad\A_E\coloneqq\sup_{\omega\in E}\A_\omega, \qquad E\subset\Pd,
\end{equation}
which by \eqref{eq:homogeneity} is unambiguous. Finally, for $Z\subset\Pd$ and an integer $N\geq0$, with $\A_\varnothing\coloneqq0$, we set
\begin{equation} \label{eq:QZN}
\Q_d(Z,N) \coloneqq \sup_{\substack{E\subset Z\\ \#E\leq N}} \left\|\A_E\right\|_{L^2(\R^{d+1})\to L^2(\R^{d+1})} ,
\end{equation}
and record the trivial bounds
\begin{equation} \label{eq:Q-trivial}
\Q_d(Z,N)\leq C_\varphi N, \qquad \Q_d(Z'\cup Z'',N)\leq\Q_d(Z',N)+\Q_d(Z'',N).
\end{equation}

The maximal operator constant \eqref{eq:QZN} turns out to be invariant under the natural symmetries of $\Pd$, which we now recall.  An invertible linear map $T\in\operatorname{GL}_{d+1}(\R)$ sends lines through the origin to lines through the origin, and therefore induces a bijection of $\Pd$,
\[
[v]\longmapsto[Tv],\qquad v\in\R^{d+1}\setminus\{0\}.
\]
Two matrices induce the same bijection precisely when they differ by a nonzero scalar factor, so the group of transformations of $\Pd$ so obtained is the quotient
\[
\operatorname{PGL}_{d+1}(\R) \coloneqq \operatorname{GL}_{d+1}(\R)/\R^{\times} ,
\]
the \emph{projective linear group}.
We write $BZ$ for the image of $Z\subset\Pd$ under $B\in\operatorname{PGL}_{d+1}(\R)$, and say that $A\in\operatorname{GL}_{d+1}(\R)$ \emph{is a representative of} $B$ if $B[v]=[Av]$ for all $v\neq0$.

\begin{lemma}  \label{lem:Q-projective-invariance}
Let $B\in\operatorname{PGL}_{d+1}(\R)$, $Z\subset\Pd$ and $N\geq1$.  Then
\[
\Q_d(BZ,N)=\Q_d(Z,N).
\]
\end{lemma}

\begin{proof} Let $A\in\operatorname{GL}_{d+1}(\R)$ be a representative of $B$. Defining the unitary operator on $L^2(\R^{d+1})$
\[
U_Af\coloneqq\left|\det A\right|^{{\frac12}}f(A\cdot)
\]
leads to the covariance identity
\begin{equation} \label{eq:covariance}
  \A_{Av,r}=U_A^{-1}\A_{v,r}U_A .
\end{equation}
The lemma is then an immediate consequence of the definitions, in particular \eqref{eq:QZN}.
\end{proof}

Hereafter, for $\omega, \omega' \in \Pd$, write
\begin{equation}\label{eq:projective-metric}
\dpr(\omega,\omega') \coloneqq \frac{\left|v\wedge v'\right|}{|v|\,|v'|}=\sin\angle(v,v'),\qquad \omega=[v],\quad\omega'=[v'] .
\end{equation}
This is a metric on $\Pd$, of diameter $1$, and is the restriction to
$\Gr(1,d+1)=\Pd$ of the Grassmannian metric used in \cite{DP-PAdv}.
The metric and the frequency split below are taken in a chosen linear
frame; Lemma~\ref{lem:Q-projective-invariance} allows us to choose this
frame before estimating the maximal operator.

At low frequencies we compare averages in nearby directions; the
complementary frequencies lie in thin bands. To decompose
\eqref{eq:Aomega}, let $\Phi\in \mathcal C^\infty_c(\R^{d+1})$ be a radial
function, fixed once and for all depending only on $d$, with $\ind_{B_1}\leq\Phi\leq\ind_{B_2}$, where $B_r$ denotes the ball of radius $r$ centered at the origin, and define the multiplier operators
\begin{equation} \label{e:LPD}
 P_{\leq R} f(x) \coloneqq \int_{\mathbb R^{d+1}} \widehat f(\xi) \Phi(R^{-1} \xi)\e^{2\pi i x \cdot \xi} \, \d\xi, \quad x\in \mathbb R^{d+1}, \qquad   P_{>R}\coloneqq I-P_{\leq R}.
\end{equation}
Given a \emph{transversal scale} $0<\delta\leq1$, define the \emph{low-} and \emph{high-frequency maximal averages}
\begin{equation} \label{eq:low-high-operators} 
\A^{\flat\delta}_{[v]}f \coloneqq \sup_{\lambda>0}\left|\A_{v,\lambda}P_{\leq(\lambda\delta|v|)^{-1}}f\right|,\qquad\A^{\sharp\delta}_{[v]}f\coloneqq \sup_{\lambda>0}\left|\A_{v,\lambda}P_{>(\lambda\delta|v|)^{-1}}f\right|.
\end{equation}
It is immediate to check that either definition is independent of the representative $v$ by
\eqref{eq:homogeneity}. Also directly from the definitions,
\begin{equation} \label{eq:low-high-splitting}
\A_\omega f   \leq \A^{\flat\delta}_{\omega}f+\A^{\sharp\delta}_{\omega}f, \qquad \omega\in\Pd,\quad 0<\delta\leq1 .
\end{equation}
For $\delta=0$ we set $\A^{\flat0}_{\omega}\coloneqq\A_{\omega}$ and $\A^{\sharp0}_{\omega}\coloneqq0$, so that \eqref{eq:low-high-splitting} holds also for $\delta=0$. Finally, associate to a scale function $\delta:E\to[0,1]$, defined on a set $E\subset\Pd$, the two maximal operators
\begin{equation} \label{eq:low-high-set}
\A^{\flat\delta(\cdot)}_{E}f\coloneqq  \sup_{\omega\in E}\A^{\flat\delta(\omega)}_{\omega}f, \qquad \A^{\sharp\delta(\cdot)}_{E}f\coloneqq\sup_{\omega\in E}\A^{\sharp\delta(\omega)}_{\omega}f,
\end{equation}
abbreviated by $\A^{\flat\delta}_{E}$ and $\A^{\sharp\delta}_{E}$, which, thanks to \eqref{eq:low-high-splitting}, satisfy
\begin{equation} \label{eq:low-high-gen}
\A_Ef\leq\A^{\flat\delta}_{E}f+\A^{\sharp\delta}_{E}f.
\end{equation}

We will use two scales: $\delta(\omega)$ bounds the change of direction, while the smaller transversal scale $\eps\delta(\omega)$ determines the frequency split. The switch lemma permits this choice at a cost $\eps^{-d}$ in the comparison estimate. We can therefore narrow the high-frequency bands to obtain bounded overlap while keeping the same comparison directions.

\begin{lemma} \label{lem:switch} Let $0<\eps\leq1$ and let $E\subset\Pd$ be a finite set. Suppose that the maps $\varsigma:E\to\Pd$ and $\delta:E\to[0,1]$ satisfy
\[
\dpr\left(\omega,\varsigma(\omega)\right)\leq\delta(\omega), \qquad \omega\in E.
\]
Then
\begin{equation} \label{eq:switch}
\left\|\A^{\flat(\eps\delta)}_{E}\right\|_{L^2(\R^{d+1})\to L^2(\R^{d+1})} \lesssim_{d,\varphi} \eps^{-d} \Q_d\left(\varsigma(E),\#E\right) .
\end{equation}
\end{lemma}

\begin{proof} For $\omega\in\Pd$ and $0<\delta\leq1$ we write
\begin{equation}\label{eq:plate}
T_\delta(\omega) \coloneqq \left\{x\in\R^{d+1}: \left|\Pi_\omega x\right|\leq1,\quad \left|\Pi^{\perp}_\omega x\right|\leq\delta\right\},\qquad \left|T_\delta(\omega)\right|\simeq_d\delta^{d},
\end{equation}
for the unit-scale tube along $\omega$ and thickness $\delta$ in the directions of $\omega^\perp$, where $\Pi_\omega$ and $\Pi^\perp_\omega$ denote the orthogonal projections onto the line $\omega$ and onto $\omega^\perp$, and denote by $\M_{\mathrm{HL}}$ the Hardy--Littlewood maximal operator on $\R^{d+1}$. The positivity of $\psi$ and its lower bound on $[-1,1]$ yield the standard estimate
\begin{equation} \label{eq:plate-domination}
\avgint_{x+\rho T_\vartheta(\omega)}|f| \lesssim \M_{\mathrm{HL}}\left(\A_{\omega}|f|\right)(x), \qquad x\in\R^{d+1}, \quad\omega\in\Pd,\quad\rho>0,\quad 0<\vartheta\leq 1.
\end{equation}
If $\omega\in E$ and $\delta(\omega)=0$, then $\varsigma(\omega)=\omega$ and $\A^{\flat0}_\omega f=\A_\omega f\leq\M_{\mathrm{HL}}(\A_{\varsigma(\omega)}|f|)$ almost everywhere, which is the bound obtained below. Fix now $\omega\in E$ with $\delta(\omega)>0$, let $v$ be a unit representative of $\omega$, and write $\omega'\coloneqq\varsigma(\omega)$, $\delta\coloneqq\delta(\omega)$. Proceeding as in \cite{DP-PAdv}*{Lemma 2.4},
\begin{equation}\label{eq:switch-tail}
\left|\A_{v,\lambda}P_{\leq(\lambda\eps\delta)^{-1}}f(x)\right| \lesssim_{d,\varphi} \sum_{k\geq0}2^{-k} \avgint_{x+2^{k}\lambda T_{\eps\delta}(\omega)}|f| \lesssim \eps^{-d}  \M_{\mathrm{HL}}\left(\A_{\omega'}|f|\right)(x), \qquad \lambda>0 ,
\end{equation}
where the second step is obtained by combining the inclusion $T_{\eps\delta}(\omega)\subset 2 T_{\delta}(\omega')$ with \eqref{eq:plate-domination} and summing the series. Taking the supremum over $\omega\in E$ in \eqref{eq:low-high-set} and $\lambda>0$,
\[
\A^{\flat(\eps\delta)}_{E}f \lesssim_{d,\varphi}\eps^{-d} \M_{\mathrm{HL}}\left(\A_{\varsigma(E)}|f|\right) .
\]
Estimate \eqref{eq:switch} now follows from the $L^2(\R^{d+1})$ boundedness of $\M_{\mathrm{HL}}$, the definition \eqref{eq:QZN} and the bound $\#\varsigma(E)\leq\#E$.
\end{proof}

We turn to the complementary component. The rough Fourier restriction of $f$ to a frequency cone $\Gamma$  will also be denoted by $\Gamma f$, abusing notation.  For $\omega=[v]\in\Pd$ and $0<\delta\leq1$, consider the frequency cone
\begin{equation}\label{eq:cone}
\Gamma^{\delta}_{\omega} \coloneqq \left\{\xi\in\R^{d+1}\setminus\{0\}: \left|\xi \cdot v \right|\leq\delta|\xi||v|\right\},\qquad \Gamma_\omega ^0 \coloneqq v^\perp\setminus\{0\}.
\end{equation}

\begin{lemma}  \label{lem:high-frequency-support}
If $E\subset\Pd$, $\delta:E\to[0,1]$, and {$\Gamma\subset\R^{d+1}$} is a measurable cone with $\Gamma^{\delta(\omega)}_{\omega}\subset\Gamma$ for every $\omega\in E$, then
\begin{equation} \label{eq:high-frequency-projection}
  \A^{\sharp\delta}_{E}f=\A^{\sharp\delta}_{E}\Gamma f .
\end{equation}
\end{lemma}

\begin{proof} Fix $\omega\in E$; if $\delta(\omega)=0$, then $\A^{\sharp0}_\omega=0$ and $\omega$ contributes nothing to either side of \eqref{eq:high-frequency-projection}, so we assume $\delta(\omega)>0$. Fix a representative $v$ of $\omega$, and $\lambda>0$, and
abbreviate $\delta\coloneqq\delta(\omega)$.  By \eqref{eq:multiplier} and \eqref{e:LPD}, the operator $\A_{v,\lambda}P_{>(\lambda\delta|v|)^{-1}}$ is the Fourier multiplier with symbol
\[
m_{v,\lambda}(\xi)\coloneqq\varphi(\lambda\xi\cdot v) \left[1-\Phi\left(\lambda\delta|v|\xi\right)\right],\qquad \xi \in \R^{d+1}.
\]
Since $\supp\varphi\subset[-1,1]$ and $\Phi\equiv1$ on $B_1$, the symbol $m_{v,\lambda}$ vanishes unless
\[
\left|\xi\cdot v\right|\leq\lambda^{-1} \qquad\text{and}\qquad |\xi|\geq\left(\lambda\delta|v|\right)^{-1} ,
\]
and on the intersection of these two sets
\[
\left|\xi\cdot v\right|\leq \lambda^{-1} = \delta|v|\left(\lambda\delta|v|\right)^{-1} \leq \delta|v||\xi| .
\]
Therefore $\supp m_{v,\lambda}\subset\Gamma^{\delta(\omega)}_{\omega}\subset\Gamma$, the second inclusion by hypothesis. In particular $m_{v,\lambda}=m_{v,\lambda}\ind_{\Gamma}$ almost everywhere, whence
\[
\A_{v,\lambda}P_{>(\lambda\delta|v|)^{-1}}f=\A_{v,\lambda}P_{>(\lambda\delta|v|)^{-1}}\Gamma f,\qquad f\in L^2(\R^{d+1}) .
\]
As this identity holds for every $\lambda>0$ and every $\omega\in E$, it is preserved by the two suprema in \eqref{eq:low-high-operators} and \eqref{eq:low-high-set}, which is
\eqref{eq:high-frequency-projection}.
\end{proof}

A switch system records the comparison directions and their contact scales on each block.

\begin{definition}\label{def:switch-system} A \emph{switch system} for $Z\subset\Pd$ is a triple $\mathcal S=\left(\{Z_m\}_{m=1}^{M},\varsigma,\delta\right)$
consisting of a partition $Z=Z_1\sqcup\cdots\sqcup Z_{M}$ together with maps
\[
\varsigma:Z\to\Pd, \qquad \delta:Z\to[0,1], \qquad \dpr\left(\omega,\varsigma(\omega)\right)\leq\delta(\omega) \quad(\omega\in Z),
\]
so that $\delta(\omega)=0$ forces $\varsigma(\omega)=\omega$. Given $0<\eps\leq1$, the associated \emph{block cones} and \emph{overlap number} are
\begin{equation} \label{eq:overlap-number}
\Gamma^{\eps}_m(\mathcal S) \coloneqq \bigcup_{\omega\in Z_m}\Gamma^{\eps\delta(\omega)}_{\omega}, \qquad {\mathsf{O}_\eps(\mathcal S)}  \coloneqq \sup_{\xi\in\R^{d+1}\setminus\{0\}} \sum_{m=1}^{M}\ind_{\Gamma^{\eps}_m(\mathcal S)}(\xi) \leq M.
\end{equation}
\end{definition}

The low-frequency comparison and high-frequency localization give the following estimate.

\begin{theorem} \label{thm:abstract-ao}
Let $\mathcal S=\left(\{Z_m\}_{m=1}^{M},\varsigma,\delta\right)$ be a switch system for $Z\subset \Pd$ and  $E\subset Z$ be a finite set. Then, for each $0<\eps \leq 1$,
\begin{equation} \label{eq:claimed-ao}
\begin{split}
\left\|\A_E\right\|_{L^2(\R^{d+1})\to L^2(\R^{d+1})} & \leq C_{d,\varphi}\eps^{-d}  M^{\frac12} \max_{1 \leq m\leq M} \Q_d\left(\varsigma(Z_m),\#(E\cap Z_m)\right)
\\
& \qquad+{\mathsf{O}_\eps(\mathcal S)} ^{\frac12} \max_{1\leq m\leq M}
\Q_d\left(Z_m,\#(E\cap Z_m)\right).
\end{split}
\end{equation}
\end{theorem}

\begin{remark} The two terms on the right-hand side of \eqref{eq:claimed-ao} play different roles. Only the first carries the factor $\eps^{-d}$ and the constant $C_{d,\varphi}$; the coefficient of the second, which refers to the blocks themselves, is exactly $\mathsf O_\eps(\mathcal S)^{\frac12}$. This is what allows \eqref{eq:claimed-ao} to be iterated without accumulating a factor $\eps^{-d}$ at every step. The iteration still pays $\mathsf O_\eps(\mathcal S)^{\frac12}$ per step, and this is the source of the subpower growth in Theorem~\ref{thm:projective-SECT-maximal}.
\end{remark}

\begin{proof}Write
\[
\begin{split} & E_m\coloneqq E\cap Z_m, \quad \Gamma_m\coloneqq\bigcup_{\omega\in E_m}\Gamma^{\eps\delta(\omega)}_{\omega}\subset\Gamma^{\eps}_m(\mathcal S), \quad N_m\coloneqq\#E_m , \quad  \A^\flat\coloneqq\A^{\flat(\eps\delta)}, \quad  \A^\sharp\coloneqq\A^{\sharp(\eps\delta)},
\\
& \mathcal L\coloneqq\max_{1\leq m\leq M}\left\|\A^{\flat}_{E_m}\right\|_{2\to2}, \quad \mathcal H\coloneqq\max_{1\leq m\leq M}\Q_d\left(Z_m,N_m\right).
\end{split}
\]
The low-frequency operator satisfies
\begin{equation}\label{eq:low-bound}
  \left\|\A^{\flat}_{E}\right\|_{2\to2} \leq M^{\frac12}\mathcal L
\end{equation}
by Cauchy--Schwarz. Next, as $P_{>R}=I-P_{\leq R}$, the triangle inequality inside the supremum defining \eqref{eq:low-high-operators} yields the pointwise bound $\A^{\sharp}_{E_m}g\leq\A_{E_m}g+\A^{\flat}_{E_m}g$, whence, by \eqref{eq:QZN} and $E_m\subset Z_m$,
\begin{equation} \label{eq:block-bound}
\left\|\A^{\sharp}_{E_m}\right\|_{2\to2} \leq \mathcal H+\mathcal L.
\end{equation}
Furthermore Lemma~\ref{lem:high-frequency-support} applies on the $m$-th block with $\Gamma=\Gamma_m$, a finite union of measurable cones, and gives $\A^{\sharp}_{E_m}f=\A^{\sharp}_{E_m}\Gamma_mf$.  Therefore, using \eqref{eq:block-bound}, Plancherel, and the definition \eqref{eq:overlap-number} of the overlap number leads to the bound
\[
\left\|\A^{\sharp}_{E}f\right\|_2^2 \leq \sum_{m=1}^{M}\left\|\A^{\sharp}_{E_m}\Gamma_mf\right\|_2^2 \leq \left(\mathcal H+\mathcal L\right)^2 \int_{\R^{d+1}} |\widehat f(\xi) |^2\sum_{m=1}^{M}\ind_{\Gamma_m}(\xi)\,\d\xi \leq {\mathsf{O}_\eps(\mathcal S)} \left(\mathcal H+\mathcal L\right)^2\left\|f\right\|_2^2 .
\]
Combining this with \eqref{eq:low-bound}, invoking \eqref{eq:low-high-gen}, and using ${\mathsf{O}_\eps(\mathcal S)} \leq M$ on the term carrying $\mathcal L$,
\begin{equation} \label{e:switchingpf}
\left\|\A_{E}\right\|_{2\to2} \leq M^{\frac12}\mathcal L +{\mathsf{O}_\eps(\mathcal S)} ^{\frac12}\left(\mathcal H+\mathcal L\right) \leq 2M^{\frac12}\mathcal L+{\mathsf{O}_\eps(\mathcal S)} ^{\frac12}\mathcal H.
\end{equation}
Now, Lemma~\ref{lem:switch} applied to each $E_m$ with the restrictions of $\varsigma$ and $\delta$, together with $\varsigma(E_m)\subset\varsigma(Z_m)$ and $\#\varsigma(E_m)\leq N_m$, gives
\begin{equation} \label{eq:L-bound}
\mathcal L \lesssim_{d,\varphi} \eps^{-d} \max_{1\leq m\leq M}\Q_d\left(\varsigma(Z_m),N_m\right).
\end{equation}
Combining \eqref{eq:L-bound} and \eqref{e:switchingpf} proves \eqref{eq:claimed-ao}. The additional low-frequency term in \eqref{eq:block-bound} is absorbed using $\mathsf O_\eps(\mathcal S)\leq M$; the coefficient of $\mathcal H$ remains exactly $\mathsf O_\eps(\mathcal S)^{\frac12}$.
\end{proof}

\begin{remark}\label{rem:overlap-problem} For curves, the comparison sets will lie in lower-dimensional hyperplanes, while the original blocks will contain fewer directions. To use these two reductions in Theorem~\ref{thm:abstract-ao}, we need an overlap bound independent of the number of blocks. Identify a nonzero $\xi\in\R^{d+1}$ with the projective hyperplane $\xi^\perp\subset\Pd$ that it annihilates. By \eqref{eq:cone},
\begin{equation}
\label{eq:cone-incidence}
\xi\in\Gamma^\delta_\omega \iff\dpr(\omega,\xi^\perp)\leq\delta ,
\end{equation}
where $\dpr(\omega,\xi^\perp)\coloneqq\min_{\omega'\in\xi^\perp}\dpr(\omega,\omega')=|\xi\cdot v|/(|\xi||v|)$ for $\omega=[v]$. Thus $\mathsf O_\eps(\mathcal S)$ counts the largest number of blocks containing a direction $\omega$ within $\eps\delta(\omega)$ of a single hyperplane. In Section~\ref{sec:canonical-ect-overlap}, $\delta$ will be an osculating approximation error, and nodal interpolation will control this thickened incidence count.
\end{remark}

\section{Canonical charts, switch moduli, and thick band overlap} \label{sec:canonical-ect-overlap} We construct the comparison directions and contact scales required by Section~\ref{sec:AO}. The construction follows the heuristics of Subsection~\ref{subsec:intro-mechanism}. We carry it out in a fixed canonical chart; Section~\ref{sec:projective-switch-classes} then allows a separate choice of projective frame for each arc.

\subsection{Canonical charts} \label{subsec:canonical-charts} Throughout this section $I=[0,1]$.  For $w\in L^1(I)$, $w>0$ a.e., we make use of the Volterra operator
\[
\mathcal I_w: L^\infty(I) \to \mathcal{AC}(I), \qquad \mathcal I_w f(t)\coloneqq\int_0^t w(s)f(s)\,\d s,\qquad t\in I.
\]
A \emph{weight vector} is a tuple $\mathbf w=(w_0,\ldots,w_d)$ satisfying
\begin{equation} \label{e:normweights}
w_0\equiv w_1\equiv1, \qquad w_q\in L^1(I),\quad w_q>0\ \text{a.e.},\quad \int_I w_q(t)\,\d t=1, \qquad 2\leq q\leq d .
\end{equation}
The unit-mass normalization is a choice of gauge and will be discussed in
Remark~\ref{rem:normalization-canonical-weights}.
The \emph{canonical coordinates} associated with $\mathbf w$ are
\begin{align}\label{eq:canonical-ECT-coordinates}
X_0&\coloneqq1,\qquad X_1\coloneqq\mathcal I_{w_1}\mathbf1_I=t, \qquad X_j
\coloneqq \mathcal I_{w_1}\mathcal I_{w_2}\cdots\mathcal I_{w_j} \mathbf{1}_I, \qquad 2\leq j\leq d .
\end{align}
Composition in \eqref{eq:canonical-ECT-coordinates} is understood to be read from right to left. Below, when the tuple $\mathbf w$ is understood to be fixed or clear from context, we abuse notation and write $\mathcal I_q $ in place of $\mathcal I_{w_q}$. For $0\leq j\leq d$ set
\begin{equation} \label{eq:canonical-prefix-space}
\mathbf X^{[j]}\coloneqq(X_0,\ldots,X_j):I\to\R^{j+1}, \qquad \mathcal U_j\coloneqq\Span\{X_0,\ldots,X_j\}, \qquad \mathbf X\coloneqq\mathbf X^{[d]} .
\end{equation}
We call $\mathbf X$ the \emph{canonical curve} associated with $\mathbf w$; thus $\mathbf X^{[j]}$ is the canonical curve associated with $(w_0,\ldots,w_j)$. For $J\subset I$, we write
\[
[\mathbf X(J)]\coloneqq\left\{[\mathbf X(t)]: t\in J\right\}\subset\Pd
\]
for the corresponding set of projective directions.

\begin{remark} \label{r:indep} The coordinates $X_0,\ldots,X_d$  defined in \eqref{eq:canonical-ECT-coordinates} are linearly independent.  Indeed, if
\[
\sum_{p=0}^d a_pX_p=0,
\]
evaluation at $t=0$ gives $a_0=0$.  Factoring the remaining identity successively through the Volterra operators $\mathcal I_1,\ldots,\mathcal I_d$, which are injective since $w_q>0$ a.e., and evaluating at $0$ at each step gives $a_1=\cdots=a_d=0$.
\end{remark}

We next introduce generalized derivatives associated with $\mathbf w$.  Set
\begin{equation}\label{eq:generalized-derivatives-coordinates}
\begin{alignedat}{2}
\mathcal D_qX_p&\coloneqq 0,& & \quad 0\leq p<q\leq d,
\\
\mathcal D_qX_q&\coloneqq 1,& &\quad 0\leq q\leq d,
\\
\mathcal D_qX_p&\coloneqq \mathcal I_{q+1}\cdots\mathcal I_p 1,& & \quad 0\leq q<p\leq d,
\end{alignedat}
\end{equation}
and extend $\mathcal D_q$ to $\mathcal U_d$ by linearity, thanks to Remark \ref{r:indep}. Thus $\mathcal D_0h=h$, each $\mathcal D_qh$ is absolutely continuous on $I$, and
\begin{equation} \label{eq:generalized-derivative-recursion}
(\mathcal D_{q-1}h)'=w_q\mathcal D_qh \quad\text{a.e. on }I, \qquad 1\leq q\leq d .
\end{equation}
Since $w_q>0$ a.e., \eqref{eq:generalized-derivative-recursion} determines the continuous function $\mathcal D_qh$ from $\mathcal D_{q-1}h$: it agrees almost everywhere with $(\mathcal D_{q-1}h)'/w_q$. By \eqref{eq:generalized-derivatives-coordinates}, the matrix $\mathcal G_j(t)\coloneqq (\mathcal D_rX_q(t))_{0\leq r,q\leq j}$ is upper triangular with unit diagonal, hence
\begin{equation} \label{eq:generalized-Wronskian-one}
\det\mathcal G_j(t)=1, \qquad t\in I,\quad 0\leq j\leq d .
\end{equation}
The \emph{cumulative coordinates} of the canonical flag are
\begin{equation} \label{eq:cumulative-flag-coordinate}
G_1(t)\coloneqq t, \qquad G_q(t)\coloneqq\int_0^t w_q(s)\,\d s =(\mathcal I_q1)(t)=\mathcal D_{q-1}X_q(t), \qquad 2\leq q\leq d .
\end{equation}

We use interpolation in two ways. Values at $j+1$ distinct nodes recover a function in $\mathcal U_j$; matching its first $j$ generalized derivatives at one point selects an approximation in $\mathcal U_{j-1}$. The latter constructs the comparison curve, while the former controls reconstruction from small nodal values. Both are uniquely determined, as stated below and proved in Subsection~\ref{subsec:unisolvence}.

\begin{lemma} \label{lem:unisolvence} Assume $\mathbf w$ satisfies \eqref{e:normweights}.
\begin{enumerate}[itemsep=.4em]
\item[\emph{(i)}] Let $0\leq j\leq d$ and $0\leq t_0<\cdots<t_j\leq1$. Then the evaluation map
\[
\operatorname{ev}_{\mathbf t}:\mathcal U_j\to\R^{j+1},\qquad h\overset{\operatorname{ev}_{\mathbf t}}{\longmapsto}\left(h(t_0),\ldots,h(t_j)\right),\qquad \mathbf t\coloneqq (t_0,\ldots,t_j),
\]
is a linear isomorphism. Equivalently, for every $(c_0,\ldots,c_j)\in\R^{j+1}$ there is a unique $h\in\mathcal U_j$ such that
\[
h(t_i)=c_i,\qquad 0\leq i\leq j.
\]
\item[\emph{(ii)}] Let $1\leq j\leq d$, $\sigma\in I$, and $h\in\mathcal U_j$. Then there is a unique $\mathcal P_{\sigma,j-1}h\in\mathcal U_{j-1}$ satisfying
\begin{equation} \label{eq:generalized-osculating-interpolant}
\mathcal D_r\left(\mathcal P_{\sigma,j-1}h\right)(\sigma)=\mathcal D_rh(\sigma),\qquad 0\leq r\leq j-1.
\end{equation}
Moreover, the resulting map $\mathcal P_{\sigma,j-1}:\mathcal U_j\to\mathcal U_{j-1}$ is a linear projection.
\end{enumerate}
\end{lemma}

\begin{remark} \label{rmrk:nodalbasis} For $0\leq t_0<\cdots<t_j\leq1$, part (i) of Lemma~\ref{lem:unisolvence} gives a unique basis $\Lambda_0^{\mathbf t},\ldots,\Lambda_j^{\mathbf t}$ of $\mathcal U_j$ satisfying
\[
\Lambda_i^{\mathbf t}(t_k)=\delta_{ik},\qquad 0\leq i,k\leq j.
\]
Consequently, every $h\in\mathcal U_j$ satisfies
\[
h=\sum_{i=0}^j h(t_i)\Lambda_i^{\mathbf t}.
\]
\end{remark}

We call $\mathcal P_{\sigma,j-1}h$ the \emph{osculating interpolant} of $h$ at $\sigma$.

\subsection{The chart-level switch class} \label{subsec:two-classes} We next introduce the quantities which enter the overlap argument. For $2\leq j\leq d$, an \emph{ordered block system} of length $j+1$ is a pair $(\mathbf t,\mathbf J)$, where $\mathbf t=(t_0,\ldots,t_j)$ satisfies $0\leq t_0<\cdots<t_j\leq1$ and $\mathbf J=(J_0,\ldots,J_j)$ is a family of closed subintervals of $I$ such that
\begin{equation} \label{eq:block-system-conditions}
t_i\in J_i\quad(0\leq i\leq j), \qquad \max J_i\leq\min J_{i+1}\quad(0\leq i<j) .
\end{equation}
We denote the collection of all such systems by $\mathfrak B_j(I)$; for a compact interval $I'$, the collection $\mathfrak B_j(I')$ is defined in the same way, with $I'$ in place of $I=[0,1]$. Writing $\partial J\coloneqq\{\inf J,\sup J\}$ for the set of endpoints of an interval $J$, we define, for $2\leq j\leq d$,
\begin{equation} \label{eq:canonical-switch-modulus}
\mathsf S_j\left(\mathbf X^{[j]}\right) \coloneqq \sup_{(\mathbf t,\mathbf J)\in\mathfrak B_j(I)} \sum_{i=0}^{j} \min_{\sigma\in\partial J_i}
\left| 1-\left(\mathcal P_{\sigma,j-1}\Lambda_i^{\mathbf t}\right)(t_i) \right|.
\end{equation}
For $K\geq1$, we write $\mathbf X\in\mathrm{SECT}(d,K)$ if
\begin{equation}\label{eq:SECT-class}
\max_{2\leq j\leq d} \mathsf S_j\left(\mathbf X^{[j]}\right)\leq K .
\end{equation}
For $d=1$, we let $\mathrm{SECT}(1,K)$ consist of the canonical curve $t\mapsto(1,t)$; the condition on the switch moduli is void and the class is independent of $K$.

Each summand in \eqref{eq:canonical-switch-modulus} is a product of two factors. Since $\mathcal P_{\sigma,j-1}$ fixes $\mathcal U_{j-1}$, the osculating error of any $h\in\mathcal U_j$ is $\mathcal D_jh$ times that of $X_j$; see \eqref{eq:residual-identity}. For $h=\Lambda_i^{\mathbf t}$, evaluated at $t_i$, the summand is therefore the product of the \emph{contact error} $\min_{\sigma\in\partial J_i}|X_j(t_i)-(\mathcal P_{\sigma,j-1}X_j)(t_i)|$, the distance at $t_i$ between the curve and an endpoint osculating model, and the \emph{amplification factor} $|\mathcal D_j\Lambda_i^{\mathbf t}|$, which bounds all coefficients of the nodal function by Proposition~\ref{prop:top-coefficient-dominance}. If a frequency lies in $j+1$ block bands, its values at the nodes are at most $\eps$ times the contact errors, and reconstructing it from these values gives $1\leq\eps(j+1)\mathsf S_j(\mathbf X^{[j]})$; see Proposition~\ref{thm:switch-implies-overlap}.

\smallskip\noindent\emph{The moment curve.} For $X_q(t)=t^q/q!$, the contact error at $\sigma$ is $|t-\sigma|^j/j!$, while $|\mathcal D_j\Lambda_i^{\mathbf t}|=j!/\prod_{k\ne i}|t_i-t_k|$. Since the blocks are ordered and $t_k\in J_k$ for every $k$, the interval $J_i$ has an endpoint between $t_i$ and the adjacent node $t_{i\pm1}$ nearest to $t_i$, so the block error is at most $|t_i-t_{i\pm1}|^j/j!$. The product is then bounded by
\[
\frac{|t_i-t_{i\pm1}|^j}{\prod_{k\ne i}|t_i-t_k|}\leq1, \qquad \mathsf S_j\left(\mathbf X^{[j]}\right)\leq j+1.
\]
Every factor in the denominator is at least the nearest-node distance. Thus arbitrarily small sampling gaps cause no loss: the contact error compensates for interpolation amplification.

The definition involves only the flag of spaces $\mathcal U_{j-1}\subset\mathcal U_j$ and the order on the interval; see Proposition~\ref{prop:reparametrization-invariance}. In Proposition~\ref{prop:switch-modulus-product-form} we express the factorization through contact and nodal residuals, in the form used to estimate $\mathsf S_j$.

The point of keeping the full tower of switch moduli is that it supplies both ingredients required by the almost orthogonality argument. To state this explicitly, consider the partition
\begin{equation} \label{e:blocksZM}
\begin{split}
&Z\coloneqq[\mathbf X(I)]=\bigcup_{m=1}^M Z_m, \qquad Z_m\coloneqq[\mathbf X(I_m)], \quad 1\leq m \leq M,
\\
& 0=\omega_0<\cdots<\omega_M=1, \qquad I_m\coloneqq[\omega_{m-1},\omega_m),\ 1 \leq m \leq M-1, \quad I_M\coloneqq[\omega_{M-1},1].
\end{split}
\end{equation}
\begin{proposition} \label{prop:SECT-switch-system} Let $d\geq2$ and $\mathbf X\in\mathrm{SECT}(d,K)$.  For each ordered decomposition \eqref{e:blocksZM} there is an associated switch system $\mathcal S\coloneqq \left(\{Z_m\}_{m=1}^M,\varsigma,\delta\right)$ such that
\begin{equation} \label{eq:SECT-overlap}
\mathsf O_{\frac{1}{2K(d+1)}}(\mathcal S)\leq d.
\end{equation}
Moreover, for every $1\leq m\leq M$, the comparison set $\varsigma(Z_m)$ is contained in the union of two projectivized injective linear images of canonical curves belonging to $\mathrm{SECT}(d-1,K)$.
\end{proposition}

The remainder of the section develops the interpolation and contact estimates needed to prove Proposition~\ref{prop:SECT-switch-system}. The formal proof of the proposition is located in Subsection~\ref{ss:pss}.

\begin{remark} \label{rem:normalization-canonical-weights} The unit-mass normalization in \eqref{e:normweights} may be imposed without loss of generality. Indeed, by Lemma~\ref{lem:Q-projective-invariance} the maximal operator constant $\Q_d(Z,N)$ is invariant under projective changes of coordinates. In particular, we are free to apply a positive diagonal change
\[
X_j\longmapsto \widetilde X_j=\mu_jX_j,\qquad \mu_j>0,\quad \mu_0=\mu_1=1,
\]
which is induced by an element of $\operatorname{GL}_{d+1}(\R)$. Under this change the canonical weights transform according to
\[
w_q\longmapsto \widetilde w_q=\frac{\mu_q}{\mu_{q-1}}w_q,\qquad 2\leq q\leq d.
\]
The factors $\mu_q$, $2\leq q\leq d$, may be chosen recursively so that
\[
\int_I\widetilde w_q(t)\,\d t=1,\qquad 2\leq q\leq d.
\]
Moreover, this diagonal rescaling does not change any of the spaces
\[
\mathcal U_j=\operatorname{span}\{X_0,\ldots,X_j\}.
\]
The nodal functions are therefore unchanged. If $\widetilde{\mathcal D}_r$ denotes the generalized derivative associated with the rescaled weights, then
\[
\widetilde{\mathcal D}_rh=\mu_r^{-1}\mathcal D_rh,\qquad h\in\mathcal U_d,
\]
so the conditions \eqref{eq:generalized-osculating-interpolant} defining $\mathcal P_{\sigma,j-1}$ are unchanged as well. Consequently, the osculating interpolants and the switch moduli are unchanged. In particular, the spaces $\mathcal U_j$, the nodal bases, the osculating interpolants and the switch moduli are defined verbatim for weights $w_q\in L^1(I)$ with $w_q>0$ a.e., without the normalization of their mass, and are unchanged by it; we use this in Example~\ref{ex:moment-curve}.
\end{remark}

\subsection{Zero counting and unisolvence} \label{subsec:unisolvence} Since $\mathcal D_qX_p(0)=\delta_{qp}$ by \eqref{eq:generalized-derivatives-coordinates}, every $h\in\mathcal U_j$ satisfies
\begin{equation} \label{eq:canonical-coefficient-extraction}
h=\sum_{q=0}^{j}\mathcal D_qh(0) X_q, \qquad 0\leq j\leq d;
\end{equation}
in particular, since $\mathcal D_j$ annihilates $X_0,\ldots,X_{j-1}$ and maps $X_j$ to $1$, $\mathcal D_jh$ is the constant top canonical coefficient of $h\in\mathcal U_j$. We equip $\mathcal U_j$ with the Euclidean norm of the canonical coefficient vector,
\begin{equation} \label{eq:canonical-coefficient-norm}
\|h\|_{\mathbf X,j}\coloneqq \left(\sum_{q=0}^j|\mathcal D_qh(0)|^2\right)^{\frac12}, \qquad h\in\mathcal U_j.
\end{equation}
For weights of class $\mathcal C^d$, the canonical coordinates form an extended complete Tchebycheff system \cite{KaSt}*{Chapter I, Example 12, and Chapter XI, Theorem 1.2}, and the facts below are classical: Lemma~\ref{lem:unisolvence} restates the nonvanishing of the determinants in \cite{KaSt}*{Chapter I, Definitions 1.1 and 2.4}, and Proposition~\ref{prop:zero-counting} is \cite{KaSt}*{Chapter I, Corollary 4.1}. Since our weights are only integrable, we include the short proofs.

\begin{lemma} \label{lem:generalized-rolle} Let $h\in\mathcal U_d$ and $1\leq q\leq d$. If $\mathcal D_{q-1}h$ vanishes at two points $a<b$ of $I$, then $\mathcal D_qh$ vanishes at some point of $(a,b)$.
\end{lemma}

\begin{proof} Otherwise $\mathcal D_qh$ has constant sign on $(a,b)$, and \eqref{eq:generalized-derivative-recursion} gives $0=\int_a^bw_q\mathcal D_qh\neq0$, as $w_q>0$ almost everywhere.
\end{proof}

\begin{proposition} \label{prop:zero-counting} Let $0\leq j\leq d$. A nonzero $h\in\mathcal U_j$ has at most $j$ zeros in $I$.
\end{proposition}

\begin{proof} If $h$ has $j+1$ zeros, Lemma~\ref{lem:generalized-rolle} gives at least $j+1-q$ zeros of $\mathcal D_qh$ for each $q\leq j$, in particular a zero of the constant $\mathcal D_jh$, so $h\in\mathcal U_{j-1}$ by \eqref{eq:canonical-coefficient-extraction}. Induction on $j$ concludes.
\end{proof}

\begin{proof}[Proof of Lemma~\textup{\ref{lem:unisolvence}}]
By Proposition~\ref{prop:zero-counting}, $\operatorname{ev}_{\mathbf t}$ is injective on $\mathcal U_j$, which has dimension $j+1$ by Remark~\ref{r:indep}; this is (i). In (ii), the conditions \eqref{eq:generalized-osculating-interpolant} are a linear system for the coefficients in $X_0,\ldots,X_{j-1}$ with the unitriangular matrix $\mathcal G_{j-1}(\sigma)$; uniqueness gives $\mathcal P_{\sigma,j-1}h=h$ for $h\in\mathcal U_{j-1}$.
\end{proof}

Osculating interpolation is the confluent limit of nodal interpolation.

\begin{lemma}\label{lem:osculating-limit} Let $1\leq j\leq d$, $h\in\mathcal U_j$ and $\sigma\in I$. If the $j$-tuples $\boldsymbol\tau^{(n)}$ of distinct points of $I$ converge to $(\sigma,\ldots,\sigma)$, then the interpolants of $h$ in $\mathcal U_{j-1}$ at $\boldsymbol\tau^{(n)}$ converge to $\mathcal P_{\sigma,j-1}h$.
\end{lemma}

\begin{proof} The residual $g_n$ of the $n$-th interpolant vanishes at the $j$ nodes, and $\mathcal D_jg_n=\mathcal D_jh$. By Lemma~\ref{lem:generalized-rolle}, $\mathcal D_\ell g_n$ vanishes at some $\zeta_{\ell,n}$ between the extreme nodes, $0\leq\ell<j$. Descending in $\ell$ with \eqref{eq:generalized-derivative-recursion} and $\int_Iw_q=1$ gives $\sup_I|\mathcal D_\ell g_n|\leq|\mathcal D_jh|$, and then
\[
|\mathcal D_\ell g_n(\sigma)|\leq|\mathcal D_jh|\left|\int_{\zeta_{\ell,n}}^\sigma w_{\ell+1}\right|\longrightarrow0, \qquad 0\leq\ell<j,
\]
and the unitriangularity of $\mathcal G_{j-1}(\sigma)$ concludes.
\end{proof}

\begin{proposition} \label{prop:reparametrization-invariance} Let $2\leq j\leq d$. The switch modulus \eqref{eq:canonical-switch-modulus} depends only on the pair $\mathcal U_{j-1}\subset\mathcal U_j$ and the order of $I$; we write $\mathsf S_j(\mathcal U_{j-1}\subset\mathcal U_j)\coloneqq\mathsf S_j(\mathbf X^{[j]})$. If $\psi:I'\to I$ is an increasing homeomorphism from a compact interval $I'$, then
\[
\mathsf S_j(\mathcal U_{j-1}\circ\psi\subset\mathcal U_j\circ\psi)=\mathsf S_j(\mathcal U_{j-1}\subset\mathcal U_j),
\]
where $\mathcal U_p\circ\psi\coloneqq\{f\circ\psi: f\in\mathcal U_p\}$ and, on the left, $\mathcal P_{\sigma,j-1}$ is the confluent limit of Lemma~\ref{lem:osculating-limit}.
\end{proposition}

\begin{proof} Nodal bases are determined by $\mathcal U_j$ and the nodes and, by Lemma~\ref{lem:osculating-limit}, osculating interpolants by $\mathcal U_{j-1}$ and the base point; this is the first assertion. Composition with $\psi$ carries nodes, ordered block systems, nodal bases and confluent limits in $I'$ to those in $I$, so it preserves each sum in \eqref{eq:canonical-switch-modulus}.
\end{proof}

Osculating replacement keeps $X_0,\ldots,X_{j-1}$ and replaces $X_j$ by its interpolant in $\mathcal U_{j-1}$. Its error is therefore a single scalar residual. The next proposition shows that this residual also determines the osculating error of every $h\in\mathcal U_j$ through its top coefficient $\mathcal D_jh$.

For $2\leq j\leq d$ and $\sigma\in I$, define the \emph{order-$j$ contact residual} based at $\sigma$ by
\begin{equation} \label{eq:canonical-contact-residual}
\mathcal R_j(\cdot;\sigma) \coloneqq X_j-\mathcal P_{\sigma,j-1}X_j  \in \mathcal U_j,
\end{equation}
with $\mathcal P_{\sigma,j-1}$ the osculating interpolant from (ii) of Lemma~\ref{lem:unisolvence}. We introduce the related based Volterra operators
\[
\mathcal I_{w,\sigma}f(t)\coloneqq\int_\sigma^t w(s)f(s)\,\d s, \qquad \mathcal I_{q,\sigma}\coloneqq\mathcal I_{w_q,\sigma}, \qquad \sigma \in I,\quad 1\leq q \leq j.
\]

\begin{proposition} \label{prop:canonical-contact-residual} Let $2\leq j\leq d$ and $\sigma\in I$.  Then the following hold.

\begin{enumerate}[itemsep=.4em]
\item[\emph{(i)}] For each $h\in\mathcal U_j$,
\begin{equation} \label{eq:residual-identity}
h-\mathcal P_{\sigma,j-1}h = \left(\mathcal D_jh\right)\mathcal R_j(\cdot;\sigma) ,
\end{equation}
where, as in \eqref{eq:canonical-coefficient-extraction}, $\mathcal D_jh$ is constant and equal to the top canonical coefficient $\mathcal D_jh(0)$ of $h$.

\item[\emph{(ii)}] There holds
\begin{equation} \label{eq:canonical-contact-integral}
\mathcal R_j(t;\sigma) = \left(\mathcal I_{1,\sigma}\cdots \mathcal I_{j,\sigma}1\right)(t), \qquad t\in I .
\end{equation}

\item[\emph{(iii)}] $\mathcal D_j\mathcal R_j(\cdot;\sigma)=1$, while $\mathcal D_r\mathcal R_j(\cdot;\sigma)(\sigma)=0$ for $0\leq r\leq j-1$, and
\begin{equation} \label{eq:contact-residual-sign}
\operatorname{sgn}\mathcal R_j(t;\sigma) = \operatorname{sgn}\left((t-\sigma)^j\right),\qquad t\in I\setminus\{\sigma\}.
\end{equation}
\end{enumerate}
\end{proposition}

\begin{proof} We first prove {(ii)} and {(iii)}.  Set $V_j\equiv1$ and, recursively,
\[
V_{q-1}\coloneqq\mathcal I_{q,\sigma}V_q, \qquad 1\leq q\leq j.
\]
Thus $V_0$ is the right-hand side of \eqref{eq:canonical-contact-integral}.  A descending induction based on \eqref{eq:generalized-derivative-recursion} gives constants $c_p$, with
$c_j=1$, such that
\[
V_r=\sum_{p=r}^jc_p\mathcal D_rX_p, \qquad 0\leq r\leq j.
\]
In particular $V_0\in\mathcal U_j$, $\mathcal D_rV_0=V_r$, and therefore
\[
\mathcal D_jV_0=1, \qquad \mathcal D_rV_0(\sigma)=0,\quad 0\leq r\leq j-1.
\]
It follows that $X_j-V_0\in\mathcal U_{j-1}$ and
\[
\mathcal D_r(X_j-V_0)(\sigma)=\mathcal D_rX_j(\sigma), \qquad 0\leq r\leq j-1.
\]
Part {(ii)} of Lemma~\ref{lem:unisolvence} now gives $X_j-V_0=\mathcal P_{\sigma,j-1}X_j$, proving
\eqref{eq:canonical-contact-integral} and the first assertions of {(iii)}. The sign property follows at once from \eqref{eq:canonical-contact-integral} and the positivity of the weights.

For {(i)}, set $b=\mathcal D_jh$. Then $h-bX_j\in\mathcal U_{j-1}$, and hence
\[
\mathcal P_{\sigma,j-1}(h-bX_j)=h-bX_j.
\]
Using linearity and \eqref{eq:canonical-contact-residual} gives
\[
h-\mathcal P_{\sigma,j-1}h =b\left(X_j-\mathcal P_{\sigma,j-1}X_j\right) =b\mathcal R_j(\cdot;\sigma),
\]
as claimed.
\end{proof}

For $2\leq j\leq d$ and $\sigma\in I$, define the osculating hyperplane
\begin{equation} \label{eq:prefix-osculating-hyperplane}
H_\sigma^{[j]} \coloneqq \Span\left\{ \mathcal D_0\mathbf X^{[j]}(\sigma),\ldots, \mathcal D_{j-1}\mathbf X^{[j]}(\sigma) \right\} \subset\R^{j+1},
\end{equation}
where the generalized derivatives act componentwise.  The spanning vectors are linearly independent by \eqref{eq:generalized-derivatives-coordinates}, so $H_\sigma^{[j]}$ is indeed a hyperplane. The corresponding \emph{osculating model} is
\begin{equation} \label{eq:osculating-model-curve}
\mathbf X^{\mathrm{osc},[j]}_\sigma \coloneqq \left(X_0,\ldots,X_{j-1}, \mathcal P_{\sigma,j-1}X_j\right).
\end{equation}
Its image is contained in $H_\sigma^{[j]}$.  Indeed, if $\xi$ annihilates $H_\sigma^{[j]}$ and $h=\xi\cdot\mathbf X^{[j]}$, then $\mathcal D_rh(\sigma)=0$ for $0\leq r\leq j-1$, and hence $\mathcal P_{\sigma,j-1}h=0$ by Lemma~\ref{lem:unisolvence}{(ii)}. By linearity this is precisely $\xi\cdot\mathbf X^{\mathrm{osc},[j]}_\sigma=0$.
Finally,
\begin{equation} \label{eq:osculating-model-properties}
\mathbf X^{[j]}(t)-\mathbf X^{\mathrm{osc},[j]}_\sigma(t) = \mathcal R_j(t;\sigma)e_j, \qquad t\in I,
\end{equation}
by \eqref{eq:canonical-contact-residual}, where $e_0,\ldots,e_j$ denotes the standard basis of $\R^{j+1}$. Thus $|\mathcal R_j(t;\sigma)|$ is the distance from $\mathbf X^{[j]}(t)$ to the point $\mathbf X^{\mathrm{osc},[j]}_\sigma(t)$ of the osculating hyperplane. Since $X_0=1$, it also bounds the projective displacement, as used in Lemma~\ref{lem:canonical-switch-system}.

\subsection{Nodal interpolation and switch moduli} We now combine the contact residual with the nodal interpolation basis to obtain the announced factorization of the switch summands. Recall that \(\mathcal U_j\) is equipped with the Euclidean norm of its canonical coefficient vector, given in \eqref{eq:canonical-coefficient-norm}.
In particular, $\|\xi\cdot\mathbf X^{[j]}\|_{\mathbf X,j}=|\xi|$. Let $1\leq j\leq d$ and
$\mathbf t=(t_0,\ldots,t_j)$ with $0\leq t_0<\cdots<t_j\leq1$.  For $0\leq i\leq j$ define the \emph{barycentric weight} of $t_i$ by
\begin{equation} \label{eq:canonical-barycentric-weight}
\beta_i^{(j)}(\mathbf t)\coloneqq \left|\mathcal D_j\Lambda_i^{\mathbf t}\right|.
\end{equation}

The barycentric weights $\beta_i^{(j)}(\mathbf t)$ are strictly positive. Indeed, if $\mathcal D_j\Lambda_i^{\mathbf t}=0$, then $\Lambda_i^{\mathbf t}\in\mathcal U_{j-1}$ and vanishes at the $j$ distinct
nodes $t_k$, $k\neq i$, contradicting Proposition~\ref{prop:zero-counting}. For $0\leq k\leq d$ and $\mathbf s=(s_0,\ldots,s_k)\in I^{k+1}$ set
\begin{equation} \label{eq:raw-canonical-determinant}
\Delta_k(\mathbf s) \coloneqq \det\left(X_q(s_r)\right)_{0\leq r,q\leq k}.
\end{equation}
If $\partial_i\mathbf t$ denotes the tuple obtained by deleting $t_i$, then
\begin{equation} \label{eq:barycentric-raw-determinant}
\beta_i^{(j)}(\mathbf t)=\frac{|\Delta_{j-1}(\partial_i\mathbf t)|}{|\Delta_j(\mathbf t)|}.
\end{equation}
Indeed, writing $\Lambda_i^{\mathbf t}=\sum_{q=0}^jc_qX_q$, the nodal conditions and Cramer's rule give
\[
c_j=(-1)^{i+j}\frac{\Delta_{j-1}(\partial_i\mathbf t)}{\Delta_j(\mathbf t)},
\]
while $|c_j|=|\mathcal D_j\Lambda_i^{\mathbf t}|=\beta_i^{(j)}(\mathbf t)$ by \eqref{eq:canonical-barycentric-weight}. The resulting formula \eqref{eq:barycentric-raw-determinant} is the analogue of the classical barycentric-weight formula for Lagrange interpolation, which motivates the terminology in \eqref{eq:canonical-barycentric-weight}.

For an ordered $j$-tuple $\boldsymbol\tau=(\tau_0,\ldots,\tau_{j-1})$ of distinct points, let $p_{\boldsymbol\tau}\in\mathcal U_{j-1}$ interpolate $X_j$ at those nodes. The \emph{nodal residual} $\mathcal R_j[\boldsymbol\tau]\coloneqq X_j-p_{\boldsymbol\tau}$ is the unique element of $\mathcal U_j$ satisfying
\begin{equation} \label{eq:nodal-residual}
\mathcal D_j\mathcal R_j[\boldsymbol\tau]=1, \qquad \mathcal R_j[\boldsymbol\tau](\tau_r)=0, \quad 0\leq r\leq j-1.
\end{equation}
Equivalently,
\begin{equation} \label{eq:nodal-residual-determinant}
\mathcal R_j[\boldsymbol\tau](u) = \frac{\Delta_j(\tau_0,\ldots,\tau_{j-1},u)}{\Delta_{j-1}(\boldsymbol\tau)}.
\end{equation}
In particular, \eqref{eq:barycentric-raw-determinant} gives
\begin{equation} \label{eq:barycentric-residual}
\beta_i^{(j)}(\mathbf t) =\frac{1}{|\mathcal R_j[\partial_i\mathbf t](t_i)|}.
\end{equation}
For $2\leq j\leq d$ and a closed interval $J\subset I$ set
\begin{equation} \label{eq:order-j-contact-radius}
\rho_j(t;J)\coloneqq\min_{\sigma\in\partial J}|\mathcal R_j(t;\sigma)|, \qquad t\in J.
\end{equation}

\begin{proposition} \label{prop:switch-modulus-product-form} Let $2\leq j\leq d$ and $(\mathbf t,\mathbf J)\in\mathfrak B_j(I)$.  Then, for $0\leq i\leq j$,
\begin{equation}\label{eq:summand-product-form}
\min_{\sigma\in\partial J_i}\left|1-\left(\mathcal P_{\sigma,j-1}\Lambda_i^{\mathbf t}\right)(t_i)\right|=\rho_j(t_i;J_i)\beta_i^{(j)}(\mathbf t)=\min_{\sigma\in\partial J_i}
\frac{|\mathcal R_j(t_i;\sigma)|}{|\mathcal R_j[\partial_i\mathbf t](t_i)|}.
\end{equation}
Consequently,
\begin{equation}\label{eq:switch-modulus-product-form}
\mathsf S_j\left(\mathbf X^{[j]}\right)=\sup_{(\mathbf t,\mathbf J)\in\mathfrak B_j(I)}
\sum_{i=0}^j  \rho_j(t_i;J_i) \beta_i^{(j)}(\mathbf t).
\end{equation}
\end{proposition}

\begin{proof} Apply \eqref{eq:residual-identity} to $h=\Lambda_i^{\mathbf t}$ and evaluate at $t_i$.  Since $\Lambda_i^{\mathbf t}(t_i)=1$,
\[
1- \left(\mathcal P_{\sigma,j-1}\Lambda_i^{\mathbf t}\right)(t_i)=\left(\mathcal D_j\Lambda_i^{\mathbf t}\right)\mathcal R_j(t_i;\sigma).
\]
Taking absolute values and minimizing over $\sigma\in\partial J_i$ gives the first equality in \eqref{eq:summand-product-form}.  The second follows from \eqref{eq:barycentric-residual}. Finally, \eqref{eq:switch-modulus-product-form} follows from \eqref{eq:canonical-switch-modulus}.
\end{proof}

The factorization above uses only the top coefficient $\beta_i^{(j)}$, whereas reconstructing a normalized frequency requires control of the whole coefficient vector of $\Lambda_i^{\mathbf t}$. The next estimate provides this control. Its proof uses the $j$ zeros of the nodal function and the unit mass of the canonical weights to bound each lower generalized derivative by the top one.

\begin{proposition} \label{prop:top-coefficient-dominance}
Let $1\leq j\leq d$ and $\mathbf t=(t_0,\ldots,t_j)$ be as above.  Then, for $0\leq i\leq j$,
\begin{equation}
\label{eq:generalized-derivative-dominance}
\left\| \mathcal D_q\Lambda_i^{\mathbf t} \right\|_{L^\infty(I)} \leq \beta_i^{(j)}(\mathbf t), \qquad 0\leq q\leq j,
\end{equation}
and consequently
\begin{equation} \label{eq:top-coefficient-dominance}
\left\| \Lambda_i^{\mathbf t} \right\|_{\mathbf X,j} \leq \sqrt{j+1}\beta_i^{(j)}(\mathbf t).
\end{equation}
\end{proposition}

\begin{proof} The function $\Lambda_i^{\mathbf t}$ vanishes at the $j$ nodes $t_k$, $k\neq i$.  Lemma~\ref{lem:generalized-rolle}, applied successively between consecutive zeros, shows that for each $0\leq q<j$ there is $\zeta_q\in I$ such that $\mathcal D_q\Lambda_i^{\mathbf t}(\zeta_q)=0$. Hence \eqref{eq:generalized-derivative-recursion} and $\int_Iw_{q+1}=1$ give
\[
\left| \mathcal D_q\Lambda_i^{\mathbf t}(t) \right|= \left| \int_{\zeta_q}^t w_{q+1}(s) \mathcal D_{q+1}\Lambda_i^{\mathbf t}(s)\,\d s \right|
\leq \left\| \mathcal D_{q+1}\Lambda_i^{\mathbf t} \right\|_{L^\infty(I)}.
\]
Iteration, starting with $|\mathcal D_j\Lambda_i^{\mathbf t}|=\beta_i^{(j)}(\mathbf t)$, proves \eqref{eq:generalized-derivative-dominance}. The second estimate follows from \eqref{eq:canonical-coefficient-norm}.
\end{proof}

\subsection{Stability under the switch operations} \label{subsec:sect-closure} The two terms of Theorem~\ref{thm:abstract-ao} require two closure properties. The high-frequency term restricts the original curve to a block; the low-frequency term replaces it by an endpoint osculating model. Both retain the bound $K$: normalization of a subarc changes the basis within each restricted space, so nodal and osculating interpolation are transported together; the osculating model is a linear image of the preceding canonical prefix. This is why the class keeps the full flag $\mathcal U_0\subset\cdots\subset\mathcal U_d$.

More precisely, let $J=[a,b]\subset I$ be nondegenerate, that is, $a<b$, set $\ell_J=b-a$ and
\[
\begin{split}
&\varphi_J(u)\coloneqq a+\ell_Ju,\qquad u\in I,
\\
&m_q(J)\coloneqq \int_Jw_q, \qquad 2\leq q\leq d.\end{split}
\]
The \emph{normalized restriction} of $\mathbf X$ to $J$ is the canonical curve $\mathbf X_J$ on $I$ associated with the weights
\begin{equation} \label{eq:normalized-restricted-weights}
w_0^J\equiv w_1^J\equiv1, \qquad w_q^J(u) \coloneqq \frac{\ell_J}{m_q(J)}w_q(\varphi_J(u)), \qquad 2\leq q\leq d.
\end{equation}
Thus the weights of $\mathbf X_J$ again satisfy
\eqref{e:normweights}.

\begin{lemma}  \label{lem:restriction-monotonicity} For each nondegenerate closed interval $J=[a,b]\subset I$,
\begin{equation}
\label{eq:restriction-monotonicity} \mathsf S_j\left(\mathbf X_J^{[j]}\right) \leq \mathsf S_j\left(\mathbf X^{[j]}\right), \qquad 2\leq j\leq d.
\end{equation}
Moreover, with $Y_0\coloneqq1$, $Y_1\coloneqq X_1-a$ and $Y_q\coloneqq\mathcal R_q(\cdot;a)$ for $2\leq q\leq d$, the canonical coordinates of $\mathbf X_J$ are
\begin{equation}  \label{eq:normalized-restricted-coordinates}
X_q^J= \frac{1}{\ell_Jm_2(J)\cdots m_q(J)}  Y_q\circ\varphi_J, \qquad 1\leq q\leq d,
\end{equation}
with the empty product understood as $1$. In particular, for $0\leq p\leq d$, $\mathbf X^{[p]}\circ\varphi_J$ is the image of $\mathbf X_J^{[p]}$ under a lower-triangular matrix in $\operatorname{GL}_{p+1}(\R)$ with positive diagonal entries.
\end{lemma}

\begin{proof} By \eqref{eq:canonical-contact-residual}, $Y_q\in X_q+\mathcal U_{q-1}$, and hence $\Span\{Y_0,\ldots,Y_p\}=\mathcal U_p$ for every $p$. Moreover, $Y_q=\mathcal I_{1,a}\cdots\mathcal I_{q,a}1$ for $q\geq1$, by \eqref{eq:canonical-contact-integral}. The change of variables $s=\varphi_J(u)$ in the iterated integrals \eqref{eq:canonical-ECT-coordinates} defining $X_q^J$, for which $w_1^J(u)\,\d u=\ell_J^{-1}\d s$ and $w_q^J(u)\,\d u=m_q(J)^{-1}w_q(s)\,\d s$ for $q\geq2$, therefore gives \eqref{eq:normalized-restricted-coordinates}. As $Y_q\in X_q+\mathcal U_{q-1}$, it follows that $X_q\circ\varphi_J$ is $\ell_Jm_2(J)\cdots m_q(J)X_q^J$ plus a combination of $X_0^J,\ldots,X_{q-1}^J$, which is the triangular relation, and that composition with $\varphi_J$ identifies, for each $p$, $\mathcal U_p$ restricted to $J$ with $\mathcal U_p^J$.

Now fix $(\mathbf u,\mathbf K)\in\mathfrak B_j(I)$, with $\mathbf u=(u_0,\ldots,u_j)$ and $\mathbf K=(K_0,\ldots,K_j)$, and set
\[
t_i\coloneqq\varphi_J(u_i),\qquad L_i\coloneqq\varphi_J(K_i),\qquad 0\leq i\leq j.
\]
Writing $\mathbf t=(t_0,\ldots,t_j)$ and $\mathbf L=(L_0,\ldots,L_j)$, we have $(\mathbf t,\mathbf L)\in\mathfrak B_j(I)$ and $L_i\subset J$ for every $i$. Let $\Lambda_i^{\mathbf u,J}$ denote the nodal basis function of $\mathcal U_j^J$ associated with $u_i$. Then
\[
\Lambda_i^{\mathbf u,J}=\Lambda_i^{\mathbf t}\circ\varphi_J,\qquad 0\leq i\leq j.
\]
If $\mathcal P_{\eta,j-1}^J$ denotes the osculating interpolant in $\mathcal U_{j-1}^J$ based at $\eta\in I$, then, with $\sigma=\varphi_J(\eta)$,
\[
\mathcal P_{\eta,j-1}^J(h\circ\varphi_J)=\left(\mathcal P_{\sigma,j-1}h\right)\circ\varphi_J,\qquad h\in\mathcal U_j.
\]
To see this we first note that for $1\leq r\leq j$ and every $h\in \mathcal U_j$, there holds
\begin{equation}\label{eq:Dr}
\mathcal D_r^J(h\circ\varphi_J)=\ell_Jm_2(J)\cdots m_r(J)\left(\mathcal D_rh\right)\circ\varphi_J,
\end{equation}
with the empty product understood as $1$. Identity \eqref{eq:Dr} follows by induction on $r$, using the chain rule, \eqref{eq:normalized-restricted-weights} and the characterization of $\mathcal D_r$ after \eqref{eq:generalized-derivative-recursion}. Applying \eqref{eq:Dr} to $h\in\mathcal U_j$ and to $\mathcal P_{\sigma,j-1}h\in\mathcal U_{j-1}\subset\mathcal U_j$, and using \eqref{eq:generalized-osculating-interpolant}, we obtain
\[
\mathcal D_r^J\left((\mathcal P_{\sigma,j-1}h)\circ\varphi_J\right)(\eta)=\mathcal D_r^J(h\circ\varphi_J)(\eta),\qquad 0\leq r\leq j-1,
\]
the case $r=0$ being simply $\mathcal P_{\sigma,j-1}h(\sigma)=h(\sigma)$. Since $(\mathcal P_{\sigma,j-1}h)\circ\varphi_J\in\mathcal U_{j-1}^J$, the uniqueness in (ii) of Lemma~\ref{lem:unisolvence} gives $(\mathcal P_{\sigma,j-1}h)\circ\varphi_J=\mathcal P_{\eta,j-1}^J(h\circ\varphi_J)$. Consequently,
\[
\min_{\eta\in\partial K_i}\left|1-\left(\mathcal P_{\eta,j-1}^J\Lambda_i^{\mathbf u,J}\right)(u_i)\right|
=
\min_{\sigma\in\partial L_i}\left|1-\left(\mathcal P_{\sigma,j-1}\Lambda_i^{\mathbf t}\right)(t_i)\right|,
\qquad 0\leq i\leq j.
\]
Thus every switch sum for $\mathbf X_J^{[j]}$ agrees with a switch sum for $\mathbf X^{[j]}$, and taking the supremum gives \eqref{eq:restriction-monotonicity}.
\end{proof}

For $d\geq2$, recall that the osculating hyperplane $H_\sigma^{[d]}$ is defined in \eqref{eq:prefix-osculating-hyperplane}. We say that a parametrized curve $\Gamma:I\to H_\sigma^{[d]}$ is \emph{linearly equivalent} to a parametrized curve $\widetilde\Gamma:I\to\R^d$ if there exists a linear isomorphism $A:\R^d\to H_\sigma^{[d]}$ such that $\Gamma(t)=A\widetilde\Gamma(t)$ for all $t\in I$.

\begin{proposition} \label{prop:SECT-closure}
Let $d\geq2$ and $\mathbf X\in\mathrm{SECT}(d,K)$. Then:
\begin{enumerate}[itemsep=.4em]
\item[\emph{(i)}] $\mathbf X^{[m]}\in\mathrm{SECT}(m,K)$ for every $1\leq m\leq d$;
\item[\emph{(ii)}] $\mathbf X_J\in\mathrm{SECT}(d,K)$ for every nondegenerate closed interval $J\subset I$;
\item[\emph{(iii)}] if $J\subset I$ is nondegenerate and $\sigma\in I$, then $\mathbf X^{\mathrm{osc},[d]}_\sigma\circ\varphi_J$ is linearly equivalent to $\mathbf X_J^{[d-1]}$.
\end{enumerate}
\end{proposition}

\begin{proof} Item {(i)} follows from \eqref{eq:SECT-class}, and for $m=1$ from the convention for $\mathrm{SECT}(1,K)$; {(ii)} is Lemma~\ref{lem:restriction-monotonicity}. For {(iii)}, since $\mathcal P_{\sigma,d-1}X_d\in\mathcal U_{d-1}$, there exist $a_0,\ldots,a_{d-1}\in\R$ such that
\[
\mathcal P_{\sigma,d-1}X_d=\sum_{q=0}^{d-1}a_qX_q.
\]
Hence
\[
\mathbf X^{\mathrm{osc},[d]}_\sigma =A_\sigma\mathbf X^{[d-1]},
\]
where
\[
A_\sigma(x_0,\ldots,x_{d-1}) \coloneqq \left(x_0,\ldots,x_{d-1},\sum_{q=0}^{d-1}a_qx_q\right).
\]
The map $A_\sigma$ is injective, and its image, the hyperplane $x_d=\sum_{q<d}a_qx_q$, contains the vectors $\mathcal D_r\mathbf X^{[d]}(\sigma)$, $r<d$, by the osculation conditions \eqref{eq:generalized-osculating-interpolant}; hence $A_\sigma:\R^d\to H_\sigma^{[d]}$ is a linear isomorphism. Therefore
\[
\mathbf X^{\mathrm{osc},[d]}_\sigma\circ\varphi_J= A_\sigma\left(\mathbf X^{[d-1]}\circ\varphi_J\right).
\]
By Lemma~\ref{lem:restriction-monotonicity}, there is an invertible linear map $B_J:\R^d\to\R^d$ such that
\[
\mathbf X^{[d-1]}\circ\varphi_J=B_J\mathbf X_J^{[d-1]}.
\]
Consequently,
\[
\mathbf X^{\mathrm{osc},[d]}_\sigma\circ\varphi_J= A_\sigma B_J\mathbf X_J^{[d-1]},
\]
and $A_\sigma B_J:\R^d\to H_\sigma^{[d]}$ is a linear isomorphism. Thus $\mathbf X^{\mathrm{osc},[d]}_\sigma\circ\varphi_J$ is linearly equivalent to $\mathbf X_J^{[d-1]}$.
\end{proof}

\subsection{Proof of  Proposition~\ref{prop:SECT-switch-system}} \label{ss:pss}
Fix an ordered decomposition \eqref{e:blocksZM}.  For $t\in I_m$, choose $\sigma_m(t)\in\partial I_m$ so that
\[
\rho_d(t;\overline{I_m})=|\mathcal R_d(t;\sigma_m(t))|.
\]
We compare $\mathbf X(t)$ with the corresponding point of the endpoint osculating model and take the contact radius itself as switching scale:
\begin{equation} \label{eq:canonical-switch-map}
\varsigma\left([\mathbf X(t)]\right)\coloneqq \left[\mathbf X^{\mathrm{osc},[d]}_{\sigma_m(t)}(t)\right],\qquad \delta\left([\mathbf X(t)]\right) \coloneqq \rho_d(t;\overline{I_m}).
\end{equation}
Since $X_0\equiv1$ and $X_1(t)=t$, the map $t\mapsto[\mathbf X(t)]$ is injective, so the definitions in \eqref{eq:canonical-switch-map} give well-defined functions $\varsigma$ and $\delta$ on $Z$.

\begin{lemma} \label{lem:canonical-switch-system} The triple $\mathcal S\coloneqq\left(\{Z_m\}_{m=1}^M,\varsigma,\delta\right)$ is a switch system in the sense of Definition~\ref{def:switch-system}. Moreover, for every $1\leq m\leq M$,
\begin{equation} \label{eq:comparison-set-osculating}
\varsigma(Z_m) \subset \bigcup_{\sigma\in\partial I_m} \left[ \mathbf X^{\mathrm{osc},[d]}_\sigma(\overline{I_m}) \right].
\end{equation}
\end{lemma}

\begin{proof}[Proof of Lemma \ref{lem:canonical-switch-system}] Fix $t\in I_m$, and set $x\coloneqq \mathbf X(t)$, $y\coloneqq \mathbf X^{\mathrm{osc},[d]}_{\sigma_m(t)}(t)$. By \eqref{eq:osculating-model-properties},
\[
x-y=\mathcal R_d(t;\sigma_m(t))e_d,
\]
and hence
\[
\dpr\left([x],[y]\right)=\frac{|x\wedge y|}{|x||y|}\leq\frac{|x-y|}{|x|}\leq\rho_d(t;\overline{I_m}),
\]
since $X_0\equiv1$ and therefore $|\mathbf X(t)|\geq1$. Moreover, \eqref{eq:canonical-contact-integral} and \eqref{e:normweights} give $\rho_d(t;\overline{I_m})\leq1$. Thus $\dpr(\omega,\varsigma(\omega))\leq\delta(\omega)$ with $\delta:Z\to[0,1]$. If $\delta(\omega)=0$, then $\mathcal R_d(t;\sigma_m(t))=0$ by the choice of $\sigma_m(t)$, so $x=y$ and $\varsigma(\omega)=\omega$, as required. The inclusion \eqref{eq:comparison-set-osculating} is immediate from the definition of $\varsigma$.
\end{proof}

It remains to control the overlap number \eqref{eq:overlap-number}. At zero thickness, incidence of one projective hyperplane with $d+1$ ordered blocks would produce $d+1$ distinct zeros of $h(t)=\xi\cdot\mathbf X(t)$, contradicting Proposition~\ref{prop:zero-counting}.  At positive thickness those zeros are replaced by values of $h$ controlled by the corresponding contact radii. The switch modulus is precisely what makes the same argument quantitative.

\begin{proposition} \label{thm:switch-implies-overlap} Let $0<\eps\leq1$ and suppose that
\begin{equation} \label{eq:epsilon-switch-overlap}
\eps(d+1) \mathsf S_d(\mathbf X) <1.
\end{equation}
Then the switch system $\mathcal S$ associated with \eqref{e:blocksZM} satisfies
$\mathsf O_\eps(\mathcal S)\leq d$.
\end{proposition}

\begin{proof} Suppose instead that some
$\xi\in\R^{d+1}\setminus\{0\}$ belongs to $d+1$ distinct block cones. By homogeneity assume $|\xi|=1$. There are then $1\leq m_0<\cdots<m_d\leq M$ and
$t_0<\cdots<t_d$ with $t_i\in I_{m_i}$ such that
\[
\xi\in \Gamma^{\eps\rho_d(t_i;\overline{I_{m_i}})}_{[\mathbf X(t_i)]}, \qquad 0\leq i\leq d.
\]
The function $h(t)\coloneqq \xi \cdot \mathbf X(t) \in \mathcal U_d$ has $\|h\|_{\mathbf X,d}=1$, $\xi$ being its canonical coefficient vector. By the definition of the cones,
\begin{equation}\label{eq:block-incidence-small-value}
|h(t_i)|\leq \eps\rho_d(t_i;\overline{I_{m_i}})|\mathbf X(t_i)|\leq\eps\sqrt{d+1}\rho_d(t_i;\overline{I_{m_i}})
\end{equation}
as $|\mathbf X(t)|\leq\sqrt{d+1}$ by \eqref{eq:canonical-ECT-coordinates} and \eqref{e:normweights}. We expand $h$ by nodal interpolation at $\mathbf t=(t_0,\ldots,t_d)$:
\[
h=\sum_{i=0}^dh(t_i)\Lambda_i^{\mathbf t}.
\]
Taking canonical coefficient norms \eqref{eq:canonical-coefficient-norm} and using the triangle inequality, then Proposition~\ref{prop:top-coefficient-dominance} together with \eqref{eq:block-incidence-small-value}, we obtain
\begin{equation}\label{eq:overlap-switch-bound}
1=\|h\|_{\mathbf X,d}\leq\sum_{i=0}^d|h(t_i)|\,\|\Lambda_i^{\mathbf t}\|_{\mathbf X,d} \leq \eps(d+1)\sum_{i=0}^d \rho_d(t_i;\overline{I_{m_i}})\beta_i^{(d)}(\mathbf t) \leq \eps(d+1)\mathsf S_d(\mathbf X),
\end{equation}
where the last bound follows from \eqref{eq:switch-modulus-product-form}, since the pairs $(t_i,\overline{I_{m_i}})$, $0\leq i\leq d$, form an admissible ordered block system in $\mathfrak B_d(I)$.  A contradiction with  \eqref{eq:epsilon-switch-overlap}  is thus reached, and the proof is complete.
\end{proof}

\begin{proof}[Proof of Proposition~\textup{\ref{prop:SECT-switch-system}}, conclusion] Lemma~\ref{lem:canonical-switch-system} gives the switch system
$\mathcal S$ associated with \eqref{e:blocksZM}. Since $\mathsf S_d(\mathbf X)\leq K$, Proposition~\ref{thm:switch-implies-overlap} applied with $\eps=\frac{1}{2K(d+1)}$ gives \eqref{eq:SECT-overlap}. By \eqref{eq:comparison-set-osculating}, each comparison set $\varsigma(Z_m)$ is contained in the union of the two endpoint osculating models over $I_m$. Proposition~\ref{prop:SECT-closure}{(iii)} identifies each of these, after transport to the reference interval, with a projectivized linear image of $\mathbf X_{\overline{I_m}}^{[d-1]}$, while {(i)}--{(ii)} of Proposition~\ref{prop:SECT-closure} give $\mathbf X_{\overline{I_m}}^{[d-1]}\in\mathrm{SECT}(d-1,K)$. This proves both assertions.
\end{proof}

\begin{remark}  \label{rem:planar-doubling-switch} Let $\mathbf X(t)=(1,t,F(t))$ be canonical. Assume that $\d\mu=F''(t)\,\d t$ is $K$-doubling, in the sense that for each interval $J\subset[0,1]$,
\[
K \min\left\{\mu(J^-), \mu(J^+)\right\} \geq \mu(J)
\]
where $J^-$ and $J^+$ denote the left and right halves of $J$. We claim that
\[
\mathsf S_2(\mathbf X)\leq1+2K.
\]
Fix $t_0\in (0,1)$, and $h,k>0$ such that $t_0-h,t_0+k\in [0,1]$ as well.  The difference of secant slopes $s$ then satisfies
\[
\begin{split}
&\quad s\coloneqq\frac{F(t_0+k)-F(t_0)}{k}-
\frac{F(t_0)-F(t_0-h)}{h}
\\
&= \int_{t_0-h}^{t_0}\left(u-(t_0-h)\right) \frac{\d\mu(u)}{h} +\int_{t_0}^{t_0+k}(t_0+k-u)\frac{\d\mu(u)}{k} \geq\frac{\mu([t_0-h,t_0]) + \mu([t_0,t_0+k]) }{2K}.
\end{split}
\]
The three absolute nodal residuals at $t_0-h,t_0,t_0+k$ are $hs$, $hks/(h+k)$ and $ks$, respectively. By positivity of the contact integral, the block tangent errors are bounded by those at the adjacent nodes; see Lemma~\ref{lem:switch-node-reduction}. The resulting upper bounds for the two outer-node terms sum to $1$, by the displayed identity for $s$. The middle-node tangent error is at most
\[
\min\left\{h\mu([t_0-h,t_0]),  k\mu([t_0,t_0+k])\right\}.
\]
Thus, for every ordered block system with nodes $\boldsymbol\tau=(t_0-h,t_0,t_0+k)$,
\[
\begin{aligned}
\sum_{i=0}^2\rho_2(\tau_i;J_i)\beta_i^{(2)}(\boldsymbol\tau) &\leq 1+\frac{h+k}{hks} \min\left\{h\mu([t_0-h,t_0]), k\mu([t_0,t_0+k])\right\}
\\
&\leq 1+\frac{\mu([t_0-h,t_0])+\mu([t_0,t_0+k])}{s} \leq 1+2K.
\end{aligned}
\]
Here $\mu$ measures the change of slope in the chosen graph chart. Doubling prevents this turning mass from concentrating where its contribution to the secant error is too small. Section~\ref{sec:doubling-curvatures} develops the analogous comparison for all canonical weights and relates it to Frenet curvature measures.
\end{remark}

\section{The projective \texorpdfstring{$\mathrm{SECT}$}{SECT} class and the maximal theorem} \label{sec:projective-switch-classes} The distances and coefficient norms in Section~\ref{sec:canonical-ect-overlap} were evaluated in a canonical chart. By Lemma~\ref{lem:Q-projective-invariance}, we may choose that chart before estimating $\Q_d$. We therefore require the existence of a projective frame in which the curve is canonical and its switch
moduli are bounded.

\begin{definition} \label{def:projective-SECT} Let $d\geq2$ and $K\geq1$. We write $\mathcal C\in\mathrm{PSECT}(d,K)$ if there are a projective transformation $B\in\operatorname{PGL}_{d+1}(\R)$ and a canonical curve $\mathbf X\in\mathrm{SECT}(d,K)$ such that
\begin{equation} \label{eq:projective-SECT-chart}
B\mathcal C=[\mathbf X(I)].
\end{equation}
For $d=1$, the condition on the switch moduli is void, and $\mathrm{PSECT}(1,K)$ is the projective saturation of the canonical curve $t\mapsto(1,t)$, $t\in I$; in particular, the class is independent of $K$.
\end{definition}

Equivalently,
\begin{equation}\label{eq:projective-SECT-saturation}
\mathrm{PSECT}(d,K)=\left\{ B[\mathbf X(I)]: B\in\operatorname{PGL}_{d+1}(\R),
\mathbf X\in\mathrm{SECT}(d,K)\right\}.
\end{equation}
Thus $\mathrm{PSECT}(d,K)$ is projectively invariant by definition. The condition is the existence of a canonical chart switch moduli; no invariance of the numerical values of the switch moduli under arbitrary changes of chart is asserted. For $d\geq2$, Proposition~\ref{prop:SECT-closure} and Lemma~\ref{lem:restriction-monotonicity} imply that $\mathrm{PSECT}(d,K)$ is closed under passage to nondegenerate closed subarcs; for $d=1$ this follows directly from the definition. Moreover, Proposition~\ref{prop:SECT-switch-system} shows that, in a canonical chart, each comparison set is contained in the union of two projectivized injective linear images of curves in $\mathrm{SECT}(d-1,K)$.

\begin{theorem} \label{thm:projective-SECT-maximal} Let $d\geq1$, $K\geq1$ and $\alpha>0$.  There is $C_{d,\alpha,\varphi}<\infty$ such that
\begin{equation} \label{eq:projective-SECT-maximal}
\sup_{\mathcal C\in\mathrm{PSECT}(d,K)} \Q_d(\mathcal C,N) \leq C_{d,\alpha,\varphi} K^{\frac{d(d+1)}2} N^\alpha, \qquad N\geq1.
\end{equation}
More precisely, put
\[
 p_d=\sum_{j=2}^d j=\frac{d(d+1)}2-1,
 \qquad b_d=1+\sum_{j=2}^d\sqrt{\log j},
\]
where both sums are empty for $d=1$, so $p_1=0$ and $b_1=1$. Then
\begin{equation}\label{eq:projective-SECT-subpower} 
\sup_{\mathcal C\in\mathrm{PSECT}(d,K)}\Q_d(\mathcal C,N) \leq C_{d,\varphi}K^{p_d} \exp \left(b_d\sqrt{\log(2+N)}\right).
\end{equation}
\end{theorem}

Each endpoint osculating model is a linear image of a canonical prefix. The following lemma transfers its maximal norm to the lower-dimensional space without loss, independently of the linear embedding.

\begin{lemma} \label{lem:hyperplane-reduction}
Let $\iota:\R^d\to\R^{d+1}$ be an injective linear map and let $Z\subset{\mathbb R}\mathrm P^{d-1}$. Writing $[\iota][v]\coloneqq[\iota v]$, one has
\begin{equation} \label{eq:hyperplane-reduction}
\Q_d\left([\iota]Z,N\right) =\Q_{d-1}(Z,N), \qquad N\geq1.
\end{equation}
\end{lemma}

\begin{proof}Since $\iota$ is linear and injective, we may choose $A\in\mathrm{GL}_{d+1}(\R)$ such that $A(v,0)=\iota v$ for every $v\in\R^d$, and let $j(v)=(v,0)$. The projective transformation induced by $A$ maps $[j]Z$ onto $[\iota]Z$. Hence, by Lemma~\ref{lem:Q-projective-invariance},
\[
\Q_d([\iota]Z,N)=\Q_d([j]Z,N).
\]
Thus it remains to prove the claim for the standard embedding $j$.

For a finite $E\subset Z$, let $E_0$ be its image under this embedding. Writing $x=(x',x_{d+1})$, the defining multipliers give
\[
\A_{E_0}f(x',x_{d+1})=\A_E^{(d)}\left(f(\cdot,x_{d+1})\right)(x'),
\]
where $\A^{(d)}$ denotes the corresponding maximal operator on $\R^d$. Fubini's theorem gives
\[
\|\A_{E_0}\|_{L^2(\R^{d+1})\to L^2(\R^{d+1})} \leq \|\A_E^{(d)}\|_{L^2(\R^d)\to L^2(\R^d)}.
\]
Testing on tensor products gives the converse. Taking the supremum over $E$ gives \eqref{eq:hyperplane-reduction}.
\end{proof}

For $d\geq1$, $K\geq1$ and $N\geq1$, set
\begin{equation} \label{eq:projective-SECT-counting}
\mathcal N_d(K;N) \coloneqq \sup_{\mathcal C\in\mathrm{PSECT}(d,K)} \Q_d(\mathcal C,N).
\end{equation}
By \eqref{eq:Q-trivial}, $\mathcal N_d(K;N)\leq C_\varphi N$.

\begin{proposition} \label{prop:projective-SECT-recursion}
Let $d\geq2$ and $K\geq1$. For each integer $L\geq2$ and every $N\geq1$,
\begin{equation}\label{eq:projective-SECT-recursion}
\mathcal N_d(K;N)\leq\sqrt d \mathcal N_d\left(K;\left\lceil\frac NL\right\rceil\right)+C_{d,\varphi}K^{d}L^{\frac12}\mathcal N_{d-1}\left(K;\left\lceil\frac NL\right\rceil\right).
\end{equation}
The coefficient of the self-referential term does not depend on $K$, on $L$, or on $\varphi$.
\end{proposition}

\begin{proof}Fix $\mathcal C\in\mathrm{PSECT}(d,K)$ and $E\subset\mathcal C$ with $\#E\leq N$. Choose a chart in \eqref{eq:projective-SECT-chart}. By Lemma~\ref{lem:Q-projective-invariance}, we may assume $\mathcal C=[\mathbf X(I)]$, with $
\mathbf X\in\mathrm{SECT}(d,K)$. Let
\[
T\coloneqq\left\{t\in I: [\mathbf X(t)]\in E\right\}.
\]
Since $X_0=1$ and $X_1(t)=t$, the projective parametrization $t\mapsto[\mathbf X(t)]$ is injective. Dividing the ordered set $T$ into at most $L$ consecutive groups, and choosing the interior cut-points outside $T$, gives an ordered decomposition \eqref{e:blocksZM}, with $M\leq L$, such that
\begin{equation}\label{eq:projective-SECT-block-cardinality}
\#(E\cap Z_m)\leq n,\qquad n\coloneqq\left\lceil\frac NL\right\rceil .
\end{equation}
Let $\mathcal S=(\{Z_m\}_{m=1}^M,\varsigma,\delta)$ be the corresponding  switch system as in Proposition~\ref{prop:SECT-switch-system}. With $\eps_{d,K}\coloneqq\frac{1}{2K(d+1)}$, that proposition gives $\mathsf O_{\eps_{d,K}}(\mathcal S)\leq d$. Theorem~\ref{thm:abstract-ao}, applied with this value of $\eps$, with
$\eps^{-d}=\left(2(d+1)\right)^dK^d$, with $M\leq L$, and with $\mathsf O_{\eps_{d,K}}(\mathcal S)\leq d$ in the coefficient of the second term, therefore yields
\begin{equation}\label{eq:projective-SECT-AO}
\|\A_E\|_{2\to2}\leq C_{d,\varphi}K^{d}L^{\frac12}\max_m\Q_d\left(\varsigma(Z_m),n\right)+\sqrt d\max_m\Q_d(Z_m,n).
\end{equation}
By (ii) of Proposition~\ref{prop:SECT-closure} and Lemma~\ref{lem:restriction-monotonicity}, the closure $\overline{Z_m}$ is projectively equivalent to the projectivization of a member of $\mathrm{SECT}(d,K)$, and hence belongs to $\mathrm{PSECT}(d,K)$. Therefore
\begin{equation}\label{eq:projective-SECT-block-bound}
\Q_d(Z_m,n)\leq\mathcal N_d(K;n).
\end{equation}
On the other hand, Proposition~\ref{prop:SECT-switch-system} places $\varsigma(Z_m)$ in the union of two projectivized images, under injective linear maps $\R^d\to\R^{d+1}$, of curves in $\mathrm{SECT}(d-1,K)$. By \eqref{eq:Q-trivial} and Lemma~\ref{lem:hyperplane-reduction},
\begin{equation} \label{eq:projective-SECT-comparison-bound}
\Q_d\left(\varsigma(Z_m),n\right)\leq2\mathcal N_{d-1}(K;n).
\end{equation}
Substituting \eqref{eq:projective-SECT-block-bound} and \eqref{eq:projective-SECT-comparison-bound} into \eqref{eq:projective-SECT-AO}, and then taking the supremum over $E$ and $\mathcal C$ proves \eqref{eq:projective-SECT-recursion}.
\end{proof}

\begin{proof}[Proof of Theorem~\ref{thm:projective-SECT-maximal}]
We prove \eqref{eq:projective-SECT-subpower} by induction on $d$.
For $d=1$ the bound $\mathcal N_1(K;N)\lesssim_\varphi\log(2+N)$
is the planar theorem of Katz \cite{KatzDuke}. Since
$\log(2+N)\leq\exp(\sqrt{\log(2+N)})$, it implies the stated bound
with $p_1=0$ and $b_1=1$.
Fix $d\geq2$, $N\geq2$, and an integer $L\geq2$. Put
$a_d=\sqrt d$ and $q=\lceil\log N/\log L\rceil$.
Successive applications of Proposition~\ref{prop:projective-SECT-recursion}
give
\begin{equation}\label{eq:projective-SECT-C0}
\begin{split}
\mathcal N_d(K;N) &\leq a_d^q\mathcal N_d(K;1)+C_{d,\varphi}K^d L^{1/2} \sum_{i=0}^{q-1}a_d^i \mathcal N_{d-1} \left(K;\left\lceil N/L^{i+1}\right\rceil\right)
\\
&\leq C_{d,\varphi}a_d^q \left(1+K^dL^{1/2}\mathcal N_{d-1}(K;N)\right).
\end{split}
\end{equation}
Here nested ceilings collapse because $L$ is an integer, and the geometric sum is at most $C_da_d^q$ since $a_d>1$. The single-direction bound is uniform in $K$. Choose
\begin{equation}\label{eq:projective-SECT-L-choice}
 L=\left\lceil\exp \left(\sqrt{(\log d)\log N}\right)\right\rceil.
\end{equation}
Writing $u=\sqrt{(\log d)\log N}$, we have $u\leq\log L\leq u+\log2$ and
\[
 q\log a_d+\tfrac12\log L \leq\frac{(\log N)\log a_d}{u}+\log a_d +\tfrac12u+\tfrac12\log2 \leq\sqrt{(\log d)\log N}+C_d.
\]
Insert the inductive bound in \eqref{eq:projective-SECT-C0}. Since $p_d=d+p_{d-1}$ and $b_d=b_{d-1}+\sqrt{\log d}$, this proves \eqref{eq:projective-SECT-subpower}; $N=1$ follows from the single-direction estimate. Finally,
\[
 \exp \left(b_d\sqrt{\log(2+N)}\right) \leq C_{d,\alpha}N^\alpha
\]
for every $\alpha>0$, and $p_d\leq d(d+1)/2$. This gives \eqref{eq:projective-SECT-maximal} as well.
\end{proof}

\section{Power-unit charts and direct switch control} \label{sec:power-unit} The switch moduli compare osculating and nodal approximation errors in the same canonical flag. We now bound them when the canonical weights are powers times positive units of bounded multiplicative oscillation. The positive integral formulas for the residuals and interpolation determinants let us compare both errors with those of a generalized moment curve, uniformly even when the nodes cluster. The resulting polynomial bound supplies the overlap estimate of Proposition~\ref{thm:switch-implies-overlap} for the charts constructed in Section~\ref{sec:algebraic}.

\begin{definition}\label{def:power-unit-class} Let $d\geq2$ and $\Lambda,R\geq1$. A canonical curve $\mathbf X$ belongs to $\mathrm{PU}(d,\Lambda,R)$ if its normalized canonical weights have the form
\begin{equation} \label{eq:power-unit-exponents}
w_q(t)=c_qt^{a_q-1}u_q(t), \qquad \Lambda^{-1}\leq a_q\leq\Lambda, \qquad 2\leq q\leq d,
\end{equation}
where the units $u_q$ are positive measurable functions satisfying
\begin{equation} \label{eq:power-unit-units}
\sup_Iu_q\leq R\inf_Iu_q,
\end{equation}
and $c_q>0$ is determined by $\int_Iw_q=1$.  We set $a_1\coloneqq1$, and call
\begin{equation} \label{eq:muntz-exponents}
\lambda_0\coloneqq0, \qquad \lambda_q\coloneqq a_1+\cdots+a_q, \qquad 1\leq q\leq d,
\end{equation}
the \emph{exponent vector} $\boldsymbol\lambda\coloneqq(\lambda_0,\ldots,\lambda_d)$ of $\mathbf X$. We write $\mathcal C\in\mathrm{PPU}(d,\Lambda,R)$ if $B\mathcal C=[\mathbf X(I)]$ for some $B\in\operatorname{PGL}_{d+1}(\R)$ and some $\mathbf X\in\mathrm{PU}(d,\Lambda,R)$.
\end{definition}

\begin{remark} \label{rem:PU-dimension-one} For $d=1$ the canonical curve is $(1,t)$ and there are no weights $w_q$ to prescribe. We therefore declare $\mathrm{PU}(1,\Lambda,R)$ to consist of this single curve for all $\Lambda,R\geq1$, and $\mathrm{PPU}(1,\Lambda,R)$ to be its $\operatorname{PGL}_2(\R)$-orbit; thus $\mathrm{PPU}(1,\Lambda,R)=\mathrm{PSECT}(1,K)$ for all $\Lambda,R,K\geq1$. With this convention Theorem~\ref{thm:projective-SECT-maximal}, which is stated for $d\geq1$, covers the degenerate case $e_\nu=1$ of Theorem~\ref{thm:algebraic-PPU-decomposition}.
\end{remark}

The estimate of the switch quantities depends on the largest exponent and on the smallest gap between consecutive exponents. We call
\begin{equation} \label{eq:muntz-conditioning}
\delta_{\boldsymbol\lambda}\coloneqq\min_{1\leq q\leq d}(\lambda_q-\lambda_{q-1}),
\qquad \kappa_{\boldsymbol\lambda}\coloneqq \frac{1+\lambda_d}{\delta_{\boldsymbol\lambda}}
\end{equation}
the \emph{conditioning} of $\boldsymbol\lambda$. Since $a_1=1$, \eqref{eq:power-unit-exponents} gives $\lambda_d\leq1+(d-1)\Lambda$ and $\delta_{\boldsymbol\lambda}\geq\Lambda^{-1}$, whence
\begin{equation}\label{eq:kappa-vs-Lambda}
\kappa_{\boldsymbol\lambda}\leq(d+1)\Lambda^2 .
\end{equation}

\begin{theorem} \label{thm:PU-PSECT}
Let $d\geq2$ and let $\mathbf X\in\mathrm{PU}(d,\Lambda,R)$ have exponent vector $\boldsymbol\lambda$.  Then
\begin{equation}
\label{eq:PU-SECT-sharp}
\max_{2\leq j\leq d} \mathsf S_j\left(\mathbf X^{[j]}\right) \leq C_dR^{d(d-1)} \kappa_{\boldsymbol\lambda}^{d(d+2)}.
\end{equation}
Consequently, for all $\Lambda,R\geq1$,
\begin{equation}
\label{eq:PU-SECT-polynomial}
\mathrm{PU}(d,\Lambda,R) \subset \mathrm{SECT}\left(d,K_d(\Lambda,R)\right), \qquad \mathrm{PPU}(d,\Lambda,R) \subset \mathrm{PSECT}\left(d,K_d(\Lambda,R)\right),
\end{equation}
where $K_d(\Lambda,R)\coloneqq C_dR^{d(d-1)}\Lambda^{2d(d+2)}.$
\end{theorem}

Combining \eqref{eq:PU-SECT-polynomial} with Theorem~\ref{thm:projective-SECT-maximal} gives the corresponding maximal estimate.

\begin{corollary} \label{cor:PPU-maximal} Let $1\leq e\leq d$, let $\Lambda,R\geq1$ and let $\alpha>0$. Suppose $\mathcal C\subset\Pd$ lies in an $e$-dimensional projective subspace $\mathbb R \mathrm{P}^e\simeq \Pi\subset\Pd$ and belongs to $\mathrm{PPU}(e,\Lambda,R)$ as a subset of $\Pi$. Then
\begin{equation}\label{eq:PPU-maximal}
\Q_d(\mathcal C,N) \lesssim_{d,\alpha,\varphi} (\Lambda R)^{C_d} N^{\alpha}, \qquad N\geq1 ,
\end{equation}
with $C_d$ depending on $d$ alone, and in particular not on $\alpha$.
\end{corollary}
Theorem~\ref{thm:PU-PSECT} and Corollary~\ref{cor:PPU-maximal} are proved in Subsection~\ref{subsec:PU-SECT-proof}.

\begin{remark} \label{rem:PU-interface} Only two features of a power-unit chart enter \eqref{eq:PU-SECT-sharp}: the multiplicative oscillation $R$ of the units and the conditioning $\kappa_{\boldsymbol\lambda}$. Thus, for the later algebraic application, it is enough to produce power-unit charts with polynomial control of these two quantities. Section~\ref{sec:algebraic} supplies such control after a projective change of frame in the span of each algebraic arc. More precisely, it gives $D^{C_d}$ pieces, each contained in a projective power-unit curve with
\[
\Lambda+R\leq D^{C_d}.
\]
Each piece is in general only a subarc of such a curve, not such a curve itself: on annular preparation cells the parameter does not reach $0$, and Lemma~\ref{lem:fill-annular-hole} extends the arc to a full power-unit curve. This is harmless: $\Q_d(\cdot,N)$ is monotone under inclusion, and if $e$ is the dimension of the span of a piece, then by \eqref{eq:PU-SECT-polynomial}, or Remark~\ref{rem:PU-dimension-one} if $e=1$, the power-unit curve containing it belongs to $\mathrm{PSECT}(e,K)$ with $K$ polynomial in $D$.
\end{remark}

The proof first replaces block endpoints by adjacent nodes (Subsection~\ref{subsec:residual-reductions}), then compares positive contact and determinant integrals to remove the units with loss $R^{j(j-1)}$ (Subsection~\ref{subsec:unit-reduction}). A power reparametrization normalizes the least exponent gap to $1$, after which generalized Vandermonde estimates give polynomial control (Subsection~\ref{subsec:muntz-switch}).

\subsection{Residual reductions} \label{subsec:residual-reductions} We begin with two elementary consequences of the positivity of the contact integral that allow us to eliminate the block variables from the definition of the switch moduli.

\begin{lemma}  \label{lem:contact-monotonicity} Let $2\leq j\leq d$ and $t\in I$.  Then
\[
\sigma\longmapsto|\mathcal R_j(t;\sigma)|
\]
is nonincreasing on $[0,t]$ and nondecreasing on $[t,1]$.
\end{lemma}

\begin{proof} By \eqref{eq:canonical-contact-integral}, $\mathcal R_j(t;\sigma)$ is, up to its fixed sign, the integral of the positive function
\[
w_2(s_2)\cdots w_j(s_j)
\]
over the simplex between $\sigma$ and $t$. The simplex increases as $\sigma$ moves away from $t$.
\end{proof}

\begin{lemma}  \label{lem:switch-node-reduction} Let $2\leq j\leq d$. Then
\begin{equation} \label{eq:switch-node-form}
\mathsf S_j(\mathbf X^{[j]}) \leq \sup_{0\leq t_0<\cdots<t_j\leq1} \sum_{i=0}^j \min_{\substack{k\in\{0,\ldots,j\}\\ |k-i|=1}}|\mathcal R_j(t_i;t_k)| \beta_i^{(j)}(\mathbf t),
\end{equation}
where at $i=0$ and $i=j$ the minimum contains only one term.
\end{lemma}

\begin{proof} Fix
$(\mathbf t,\mathbf J)\in\mathfrak B_j(I)$ and write $J_i=[\ell_i,r_i]$. For $1\leq i\leq j-1$, the ordering of the blocks and the condition $t_k\in J_k$ give
\[
t_{i-1}\leq \ell_i\leq t_i\leq r_i\leq t_{i+1}.
\]
Hence Lemma~\ref{lem:contact-monotonicity} yields
\[
|\mathcal R_j(t_i;\ell_i)|\leq|\mathcal R_j(t_i;t_{i-1})|,\qquad|\mathcal R_j(t_i;r_i)|\leq|\mathcal R_j(t_i;t_{i+1})|.
\]
Therefore
\[
\rho_j(t_i;J_i)\leq\min\left\{|\mathcal R_j(t_i;t_{i-1})|,|\mathcal R_j(t_i;t_{i+1})|\right\}.
\]
The same argument at the two endpoints gives
\[
\rho_j(t_0;J_0)\leq|\mathcal R_j(t_0;t_1)|,\qquad\rho_j(t_j;J_j)\leq|\mathcal R_j(t_j;t_{j-1})|.
\]
Multiplying by $\beta_i^{(j)}(\mathbf t)$, summing in $i$, and taking the supremum over ordered block systems proves \eqref{eq:switch-node-form}, in view of \eqref{eq:switch-modulus-product-form}.
\end{proof}

\subsection{Units and reparametrization}\label{subsec:unit-reduction} To remove the units, we compare the positive integrals representing contact residuals and canonical determinants. Positivity allows bounds for each weight to pass directly to these integrals, independently of the separation of the nodes. This gives a multiplicative comparison of the errors themselves. The determinant formula we need is the following. For $1\leq k\leq d$ and $0\leq t_0<\cdots<t_k\leq1$, let $\mathcal T_k(\mathbf t)$ be the region in the variables $s_m^{(q)}$ defined by
\begin{equation}\label{eq:staircase-region}
\begin{split}
&t_{m-1}\leq s_m^{(1)} \leq t_m,\qquad 1\leq m\leq k,
\\
& s_m^{(q-1)} \leq s_m^{(q)} \leq s_{m+1}^{(q-1)}, \qquad  2\leq q\leq k,\quad 1\leq m\leq k-q+1 .
\end{split}
\end{equation}

\begin{lemma}\label{lem:staircase-vandermonde} For $1\leq k\leq d$ and $\mathbf t$ as above,
\begin{equation}\label{eq:staircase-vandermonde}
\Delta_k(\mathbf t)=\int_{\mathcal T_k(\mathbf t)}\prod_{q=2}^k \prod_{m=1}^{k-q+1}w_q(s_m^{(q)}) \,\d\mathbf s.
\end{equation}
\end{lemma}

\begin{proof} We show, by downward induction on $0\leq p\leq k-1$, the claim
\[
\det\left(\mathcal D_pX_{p+q}(t_r)\right)_{0\leq r,q\leq k-p}=\int_{\mathcal T_{k-p}(\mathbf t)}\prod_{j=1}^{k-p}\prod_{m=1}^{k-p-j+1}w_{p+j}(s_m^{(j)}) \,\d\mathbf s
\]
for all $\mathbf t=(t_0,\ldots,t_{k-p})$ with $0\leq t_0<\cdots<t_{k-p}\leq1$; for $p=0$ the left-hand side is $\Delta_k(\mathbf t)$, and since $w_1\equiv1$ this is \eqref{eq:staircase-vandermonde}. Subtract from each row of the matrix on the left the preceding one, and expand along the first column, which becomes $(1,0,\ldots,0)$ as $\mathcal D_pX_p\equiv1$ by \eqref{eq:generalized-derivatives-coordinates}. By \eqref{eq:generalized-derivative-recursion}, the resulting matrix has entries
\[
\mathcal D_pX_{p+q}(t_r)-\mathcal D_pX_{p+q}(t_{r-1})=\int_{t_{r-1}}^{t_r}w_{p+1}\mathcal D_{p+1}X_{p+q},\qquad 1\leq r,q\leq k-p.
\]
For $p=k-1$ this is the single entry $\int_{t_0}^{t_1}w_k$, as $\mathcal D_kX_k\equiv1$, which proves the claim. For $p<k-1$, multilinearity in the rows and Fubini's theorem, which applies since the continuous functions $\mathcal D_{p+1}X_{p+q}$ are bounded, turn the left-hand side of the claim into
\[
\int_{t_0}^{t_1}\cdots\int_{t_{k-p-1}}^{t_{k-p}}\prod_{m=1}^{k-p}w_{p+1}(s_m)\det\left(\mathcal D_{p+1}X_{p+1+q}(s_m)\right)_{\substack{1\leq m\leq k-p\\0\leq q\leq k-p-1}}\,\d s_1\cdots\d s_{k-p}.
\]
For almost every $(s_1,\ldots,s_{k-p})$ in the domain of integration, $s_1<\cdots<s_{k-p}$, and the induction hypothesis at these nodes writes the last determinant as an integral over $\mathcal T_{k-p-1}(s_1,\ldots,s_{k-p})$. Setting $s_m^{(1)}\coloneqq s_m$, the levels of this region become the levels $2,\ldots,k-p$ of $\mathcal T_{k-p}(\mathbf t)$ in \eqref{eq:staircase-region}, which proves the claim.
\end{proof}

\begin{example} \label{ex:moment-curve} For an exponent vector $\boldsymbol\lambda$ as in \eqref{eq:muntz-exponents}, the \emph{generalized moment curve} is
\begin{equation}\label{eq:moment-curve}
\mathbf M_{\boldsymbol\lambda}(t)\coloneqq\left(1,t,t^{\lambda_2},\ldots,t^{\lambda_d}\right),\qquad t\in I.
\end{equation}
The weights $\widetilde w_q(t)\coloneqq a_qt^{a_q-1}$, $2\leq q\leq d$, satisfy \eqref{e:normweights}, and each coordinate of the associated canonical curve is a positive multiple of $t^{\lambda_q}$. Thus this canonical curve is obtained from $\mathbf M_{\boldsymbol\lambda}$ by a positive diagonal change of coordinates, and the two curves have the same switch moduli by Remark~\ref{rem:normalization-canonical-weights}.
\end{example}

We now compare a power-unit chart with the generalized moment curve $\mathbf M_{\boldsymbol\lambda}$ and transfer the switch bounds quantitatively.

\begin{proposition}  \label{prop:PU-to-moment}
Let $\mathbf X\in\mathrm{PU}(d,\Lambda,R)$, with weights
\[
w_q(t)=c_qt^{a_q-1}u_q(t), \qquad 2\leq q\leq d,
\]
and let $\boldsymbol\lambda$ be defined by \eqref{eq:muntz-exponents}.  Then
\begin{equation} \label{eq:PU-to-moment}
\mathsf S_j(\mathbf X^{[j]}) \leq R^{j(j-1)} \mathsf S_j(\mathbf M_{\boldsymbol\lambda}^{[j]}), \qquad 2\leq j\leq d.
\end{equation}
\end{proposition}

\begin{proof} Let $\widetilde{\mathbf X}$ be the canonical curve with the weights $\widetilde w_q$ of Example~\ref{ex:moment-curve}; it has the same switch moduli as $\mathbf M_{\boldsymbol\lambda}$. We write $\widetilde\Delta_k$, $\widetilde{\mathcal R}_j$, $\widetilde\rho_j$ and $\widetilde\beta_i^{(j)}$ for the quantities \eqref{eq:raw-canonical-determinant}, \eqref{eq:canonical-contact-residual}, \eqref{eq:order-j-contact-radius} and \eqref{eq:canonical-barycentric-weight}, respectively, associated with the curve $\widetilde{\mathbf X}$. The measure
$\d\widetilde\mu_q(t)\coloneqq\widetilde w_q(t)\,\d t$ is a probability measure, and the normalization $\int_Iw_q=1$ gives
\[
\frac{c_q}{a_q}=\left( \int_Iu_q\,\d\widetilde\mu_q\right)^{-1}\implies  \frac{w_q(t)}{\widetilde w_q(t)}=\frac{u_q(t)}{\int_Iu_q\,\d\widetilde\mu_q}\in[R^{-1},R].
\]
We compare the switch summands of $\mathbf X$ and $\widetilde{\mathbf X}$.  In the staircase representation \eqref{eq:staircase-vandermonde}, the determinant $\Delta_k$ contains
$\sum_{q=2}^k(k-q+1)=\binom{k}{2}$ weight factors. Therefore
\[
R^{-\binom{k}{2}}\widetilde\Delta_k \leq \Delta_k\leq R^{\binom{k}{2}}\widetilde\Delta_k.
\]
Likewise, the contact integral \eqref{eq:canonical-contact-integral} contains $w_2,\ldots,w_j$ once each, and hence
\[
|\mathcal R_j(t;\sigma)| \leq R^{j-1} |\widetilde{\mathcal R}_j(t;\sigma)|.
\]
Using \eqref{eq:barycentric-raw-determinant}, we obtain
\[
\rho_j(t_i;J_i)\beta_i^{(j)}(\mathbf t) \leq R^{(j-1)+\binom{j-1}{2}+\binom j2} \widetilde\rho_j(t_i;J_i) \widetilde\beta_i^{(j)}(\mathbf t).
\]
Since $(j-1)+\binom{j-1}{2}+\binom j2=j(j-1)$, summing over $i$ and taking the supremum over ordered block systems yields \eqref{eq:PU-to-moment}, in view of \eqref{eq:switch-modulus-product-form}.
\end{proof}

We next exploit the reparametrization invariance of Proposition~\ref{prop:reparametrization-invariance} to normalize the exponents. For the remainder of the section we allow an arbitrary strictly increasing exponent vector
\[
0=\lambda_0<\lambda_1<\cdots<\lambda_j
\]
and write
\begin{equation}\label{e:muntzflag}
\mathcal V_p(\boldsymbol\lambda)\coloneqq \Span\{t^{\lambda_0},\ldots,t^{\lambda_p}\}, \qquad p=j-1,j.
\end{equation}
For a general exponent vector, the first nonconstant coordinate is $t^{\lambda_1}$ rather than $t$, so this flag is not in the canonical normalization used above. We therefore work with the reparametrization-invariant order-$j$ switch modulus of the flag $\mathcal V_{j-1}(\boldsymbol\lambda)\subset\mathcal V_j(\boldsymbol\lambda)$, and denote it by $\mathsf S_j(\boldsymbol\lambda)$. By Proposition~\ref{prop:reparametrization-invariance}, composition with $\psi(s)=s^{1/\lambda_1}$ identifies this flag with $\mathcal V_{j-1}(\boldsymbol\lambda/\lambda_1)\subset\mathcal V_j(\boldsymbol\lambda/\lambda_1)$. The latter has first nonconstant exponent $1$, and is therefore spanned by the coordinates of $\mathbf M_{\boldsymbol\lambda/\lambda_1}^{[j]}$, hence, by Example~\ref{ex:moment-curve}, by those of a canonical curve. In particular, when $\lambda_1=1$, $\mathsf S_j(\boldsymbol\lambda)$ agrees with the switch modulus of the prefix $\mathbf M_{\boldsymbol\lambda}^{[j]}=(1,t,t^{\lambda_2},\ldots,t^{\lambda_j})$ of the generalized moment curve \eqref{eq:moment-curve}.

\begin{proposition}  \label{prop:muntz-scaling}
For each $c>0$, $\mathsf S_j(c\boldsymbol\lambda)=\mathsf S_j(\boldsymbol\lambda).$
\end{proposition}

\begin{proof} Put $\psi_c(s)\coloneqq s^c$, so that $t^{x}\circ\psi_c(s)=s^{cx}$. Then, for $p=j-1,j$,
\[
\mathcal V_p(\boldsymbol\lambda)=\mathcal V_p(\boldsymbol\lambda/\lambda_1)\circ\psi_{\lambda_1},\qquad \mathcal V_p(c\boldsymbol\lambda)=\mathcal V_p(\boldsymbol\lambda/\lambda_1)\circ\psi_{c\lambda_1},
\]
and, by Example~\ref{ex:moment-curve}, the flag $\mathcal V_{j-1}(\boldsymbol\lambda/\lambda_1)\subset\mathcal V_j(\boldsymbol\lambda/\lambda_1)$ is spanned by the coordinates of a canonical curve. Since $\psi_{\lambda_1}$ and $\psi_{c\lambda_1}$ are increasing homeomorphisms of $I$, Proposition~\ref{prop:reparametrization-invariance} shows that both switch moduli equal that of this canonical flag.
\end{proof}

\subsection{The switch estimate for pure powers}\label{subsec:muntz-switch} For pure powers, nodal interpolation is expressed by generalized Vandermonde determinants at distinct nodes, while contact interpolation uses their confluent limits, with repeated nodes. We estimate both by the same leading-monomial bound. After normalizing the least exponent gap, its loss is polynomial in the largest exponent. Keeping the resulting residuals in a ratio will expose the cancellation needed for the switch estimate.

Fix $0=\lambda_0<\lambda_1<\cdots<\lambda_j$. For $0\leq p\leq j$ and $\mathbf y=(y_0,\ldots,y_p)$ with $0\leq y_0<\cdots<y_p$, set
\begin{equation} \label{eq:generalized-schur}
\mathfrak s_p(\mathbf y)\coloneqq\frac{\det\left(y_r^{\lambda_q}\right)_{0\leq r,q\leq p}}{\prod_{0\leq r<s\leq p}(y_s-y_r)}.
\end{equation}
These quotients, the generalized Schur functions, are positive and symmetric. When positive nodes coalesce, we use the standard confluent convention, replacing repeated evaluation rows by the corresponding normalized derivative rows. Thus, if the distinct values taken by $\mathbf y$ are $\tau_1<\cdots<\tau_m$ with multiplicities $l_1,\ldots,l_m$, we set
\[
\mathfrak s_p(\mathbf y)\coloneqq
\frac{\det\left(\frac{1}{r!}\frac{\d^r}{\d t^r}t^{\lambda_q}\big|_{t=\tau_\nu}\right)_{\substack{1\leq\nu\leq m,\ 0\leq r<l_\nu\\0\leq q\leq p}}}
{\det\left(\frac{1}{r!}\frac{\d^r}{\d t^r}t^q\big|_{t=\tau_\nu}\right)_{\substack{1\leq\nu\leq m,\ 0\leq r<l_\nu\\0\leq q\leq p}}}.
\]
For the ordinary Vandermonde determinant and its confluent form, see for example~\cite{Shu}*{Example~2.22}. For positive nodes, this agrees with the continuous extension of \eqref{eq:generalized-schur} from distinct configurations.

We first work under the normalized gap condition
\begin{equation}\label{eq:muntz-normalized-gaps}
\lambda_q-\lambda_{q-1}\geq1,\qquad1\leq q\leq j.
\end{equation}
In particular,
\begin{equation}\label{eq:muntz-excess-monotone}
0=\lambda_0\leq\lambda_1-1\leq\cdots\leq\lambda_j-j.
\end{equation}

\begin{lemma} \label{lem:muntz-leading-monomial} Assume \eqref{eq:muntz-normalized-gaps}. Then, for every $0\leq p\leq j$ and $0<y_0\leq\cdots\leq y_p$,
\begin{equation}\label{eq:muntz-leading-monomial}
\prod_{q=0}^{p}y_q^{\lambda_q-q}\leq\mathfrak s_p(y_0,\ldots,y_p)\leq\lambda_j^{\binom{j+1}{2}} \prod_{q=0}^{p}y_q^{\lambda_q-q}.
\end{equation}
\end{lemma}

\begin{proof} For strictly increasing arguments, Proposition~5.10 of \cite{KhareTao}, applied with
\[
\mathbf n=(\lambda_0,\ldots,\lambda_p), \qquad \mathbf n_{\min}=(0,1,\ldots,p),
\]
and divided by the ordinary Vandermonde, gives
\[
\prod_{q=0}^{p}y_q^{\lambda_q-q} \leq \mathfrak s_p(y_0,\ldots,y_p) \leq
\left(\prod_{0\leq r<s\leq p}\frac{\lambda_s-\lambda_r}{s-r}
\right) \prod_{q=0}^{p}y_q^{\lambda_q-q}.
\]
The product in parentheses is at most $\lambda_j^{\binom{j+1}{2}}$.  The confluent cases follow by continuity.
\end{proof}

Recall from \eqref{eq:summand-product-form} and \eqref{eq:barycentric-residual} that the switch summands are ratios of a contact residual and a nodal residual. For the flag $\mathcal V_{j-1}(\boldsymbol\lambda)\subset\mathcal V_j(\boldsymbol\lambda)$ of \eqref{e:muntzflag}, set $X_j\coloneqq t^{\lambda_j}$, let $\mathcal R_j[\boldsymbol\tau]$ be the element of $X_j+\mathcal V_{j-1}(\boldsymbol\lambda)$ vanishing at the nodes $\boldsymbol\tau$, and let $\mathcal R_j(\cdot;\sigma)\coloneqq X_j-\mathcal P_{\sigma,j-1}X_j$, with $\mathcal P_{\sigma,j-1}$ understood as in Proposition~\ref{prop:reparametrization-invariance}. Composition with $\psi(s)=s^{1/\lambda_1}$ carries both residuals to the same positive multiple of the corresponding residuals of the canonical curve of Example~\ref{ex:moment-curve} with exponents $\boldsymbol\lambda/\lambda_1$. Hence their ratios are the switch summands of $\mathsf S_j(\boldsymbol\lambda)$, and \eqref{eq:switch-node-form} holds for $\mathsf S_j(\boldsymbol\lambda)$. We now estimate these ratios at the two adjacent nodes.

\begin{lemma} \label{lem:muntz-adjacent-residuals}
Assume \eqref{eq:muntz-normalized-gaps} and let $0<t_0<\cdots<t_j\leq1$. Then for $1\leq i\leq j-1$,
\begin{align}
\frac{|\mathcal R_j(t_i;t_{i-1})|}{|\mathcal R_j[\partial_i\mathbf t](t_i)|}&\leq\lambda_j^{j(j+1)} \left(\frac{t_{i}-t_{i-1}}{t_{i+1}-t_i}\right)^{j-i}, \label{eq:muntz-left-residual}
\\
\frac{|\mathcal R_j(t_i;t_{i+1})|}{|\mathcal R_j[\partial_i\mathbf t](t_i)|}
&\leq\lambda_j^{j(j+1)} \left(\frac{t_{i+1}-t_i}{t_{i}-t_{i-1}}\right)^i
\left(\frac{t_{i+1}}{t_i}\right)^{\lambda_j}. \label{eq:muntz-right-residual}
\end{align}
At the extreme nodes,
\begin{equation}\label{eq:muntz-edge-residuals}
\frac{|\mathcal R_j(t_0;t_1)|}{|\mathcal R_j[\partial_0\mathbf t](t_0)|},\quad\frac{|\mathcal R_j(t_j;t_{j-1})|}{|\mathcal R_j[\partial_j\mathbf t](t_j)|
}\leq\lambda_j^{j(j+1)}.
\end{equation}
\end{lemma}

\begin{proof} We first write the residuals in the monomial basis. Specializing \eqref{eq:nodal-residual-determinant} to the basis $1,t^{\lambda_1},\ldots,t^{\lambda_j}$ and factoring the generalized Vandermonde by \eqref{eq:generalized-schur} gives, up to sign,
\[
\mathcal R_j[\boldsymbol\tau](u)=\prod_{r=0}^{j-1}(u-\tau_r) \frac{\mathfrak s_j(\boldsymbol\tau,u)}{\mathfrak s_{j-1}(\boldsymbol\tau)}.
\]
To pass from nodal to contact residuals, let $P_{\boldsymbol\tau}X_j\in\mathcal V_{j-1}(\boldsymbol\lambda)$ denote the interpolant of $X_j$ at the nodes $\boldsymbol\tau$. By definition, $\mathcal R_j[\boldsymbol\tau]=X_j-P_{\boldsymbol\tau}X_j$. If $\boldsymbol\tau\to\sigma^{(j)}\coloneqq(\sigma,\ldots,\sigma)$, Lemma~\ref{lem:osculating-limit}, transported as in Proposition~\ref{prop:reparametrization-invariance}, gives
\[
P_{\boldsymbol\tau}X_j\longrightarrow\mathcal P_{\sigma,j-1}X_j.
\]
Hence $\mathcal R_j[\boldsymbol\tau]\longrightarrow\mathcal R_j(\cdot;\sigma)$. Passing to the limit in the determinant formula above therefore yields
\begin{equation}\label{eq:muntz-residual-factorization}
\frac{|\mathcal R_j(t_i;\sigma)|}{|\mathcal R_j[\partial_i\mathbf t](t_i)|
}=\frac{|t_i-\sigma|^j}{\prod_{k\neq i}|t_i-t_k|}\frac{\mathfrak s_j(\sigma^{(j)},t_i)}{\mathfrak s_j(\mathbf t)}\frac{\mathfrak s_{j-1}(\partial_i\mathbf t)
}{\mathfrak s_{j-1}(\sigma^{(j)})}.
\end{equation}
Apply Lemma~\ref{lem:muntz-leading-monomial} to the four Schur factors.
The two numerator estimates together cost at most $\lambda_j^{j(j+1)}$;
the denominator estimates cost no additional constant. Put
$E_r=\lambda_r-r$, so $0=E_0\leq\cdots\leq E_j$.
For $\sigma=t_{i-1}$ the leading-monomial quotient is
\[
A_i^-\coloneqq \frac{t_i^{E_j}\prod_{r=0}^{i-1}t_r^{E_r} \prod_{r=i}^{j-1}t_{r+1}^{E_r}} {\prod_{r=0}^jt_r^{E_r}} =\prod_{r=i}^{j-1}\left(\frac{t_{r+1}}{t_r}\right)^{E_r-E_j} \leq1.
\]
For $\sigma=t_{i+1}$, sorting the confluent tuple puts $t_i$ first, with exponent $E_0=0$, and replaces the factor $t_i^{E_j}$ by $t_{i+1}^{E_j}$. Thus its quotient is
\[
A_i^+=A_i^-\left(\frac{t_{i+1}}{t_i}\right)^{E_j} \leq\left(\frac{t_{i+1}}{t_i}\right)^{\lambda_j}.
\]
It remains to estimate the ordinary Vandermonde factor in \eqref{eq:muntz-residual-factorization}. Since the $i$ nodes to the left of $t_i$ are at distance at least $t_{i}-t_{i-1}$, and the $j-i$ nodes to its right are at distance at least $t_{i+1}-t_i$,
\[
\prod_{k\neq i}|t_i-t_k|\geq(t_{i}-t_{i-1})^{i}(t_{i+1}-t_i)^{j-i}.
\]
With $\sigma=t_{i-1}$ or $t_{i+1}$ this gives, respectively,
\[
\frac{|t_i-\sigma|^j}{\prod_{k\neq i}|t_i-t_k|}\leq\left(\frac{t_{i}-t_{i-1}}{t_{i+1}-t_i}\right)^{j-i}, \qquad \frac{|t_i-\sigma|^j}{\prod_{k\neq i}|t_i-t_k|}\leq\left(\frac{t_{i+1}-t_i}{t_{i}-t_{i-1}}\right)^i.
\]
Combining these bounds with \eqref{eq:muntz-residual-factorization} proves \eqref{eq:muntz-left-residual} and \eqref{eq:muntz-right-residual}. At $i=0$, the factor $(t_1/t_0)^{E_j}$ cancels the $r=0$ factor in $A_0^-$, since $E_0=0$, leaving
\[
 A_0^+=\prod_{r=1}^{j-1} \left(\frac{t_{r+1}}{t_r}\right)^{E_r-E_j}\leq1.
\]
At $i=j$, the product defining $A_j^-$ is empty and equals $1$. The ordinary Vandermonde factor is at most $1$ at either extreme, which proves \eqref{eq:muntz-edge-residuals}.
\end{proof}

\begin{proposition} \label{prop:muntz-switch} Let $0=\lambda_0<\lambda_1<\cdots<\lambda_j$. Then
\begin{equation}\label{eq:muntz-switch}
\mathsf S_j(\boldsymbol\lambda)\leq C_j\left(\frac{\lambda_j}{\delta_{\boldsymbol\lambda}}\right)^{j(j+2)}, \qquad \delta_{\boldsymbol\lambda}\coloneqq\min_{1\leq q\leq j}(\lambda_q-\lambda_{q-1}).
\end{equation}
\end{proposition}

\begin{proof} By Proposition~\ref{prop:muntz-scaling}, we may replace $\boldsymbol\lambda$ by $\delta_{\boldsymbol\lambda}^{-1}\boldsymbol\lambda$ and therefore assume \eqref{eq:muntz-normalized-gaps}. It is enough to prove
\begin{equation} \label{eq:normalized-Muntz-switch}
\mathsf S_j(\boldsymbol\lambda) \leq C_j\lambda_j^{j(j+2)}.
\end{equation}
Notice that the gap normalization implies $\lambda_j\geq j$. The estimates also pass to a distinct configuration with $t_0=0$. Indeed, after the power reparametrization to canonical coordinates, the contact integral \eqref{eq:canonical-contact-integral} is continuous in both endpoints. The nodal residual is continuous as well, by its determinant formula, and is nonzero at the omitted node by unisolvence. Thus the residual ratios extend continuously as $t_0\downarrow0$; no ordinary derivative of $t^{\lambda_q}$ at zero is needed. It therefore suffices to estimate configurations with $t_0>0$. By \eqref{eq:switch-node-form} and \eqref{eq:barycentric-residual}, each switch summand is bounded by the better of the two adjacent residual ratios from Lemma~\ref{lem:muntz-adjacent-residuals}. The point is that these two bounds compensate for each other. Indeed, fix an interior node $t_i$ and write
\[
a=t_i-t_{i-1}, \qquad b=t_{i+1}-t_i.
\]
If $a\leq b$, the left-hand estimate \eqref{eq:muntz-left-residual} is already bounded by $\lambda_j^{j(j+1)}$.  Suppose instead that $b<a$ and put $r=b/a\in(0,1)$. Since $a\leq t_i$,
\[
\frac{t_{i+1}}{t_i}=1+\frac{b}{t_i}\leq1+r.
\]
Thus the two estimates of Lemma~\ref{lem:muntz-adjacent-residuals}, after removing their common factor $\lambda_j^{j(j+1)}$, are bounded respectively by
\[
r^{-(j-i)} \qquad\text{and}\qquad r^i(1+r)^{\lambda_j}.
\]
These quantities cannot both be large. In fact, if $r\leq(1+\lambda_j)^{-1}$, then
\[
r^i(1+r)^{\lambda_j} \leq \e^{\lambda_j r} \leq \e,
\]
whereas if $r>(1+\lambda_j)^{-1}$, then
\[
r^{-(j-i)}\leq(1+\lambda_j)^{j-i}\leq (1+\lambda_j)^j.
\]
Consequently, at every interior node,
\[
\min_{\varepsilon\in\{-1,1\}}\frac{|\mathcal R_j(t_i;t_{i+\varepsilon})|}{|\mathcal R_j[\partial_i\mathbf t](t_i)|}\leq \e \lambda_j^{j(j+1)}(1+\lambda_j)^j.
\]
The extreme nodes satisfy the stronger bound \eqref{eq:muntz-edge-residuals}. Summing the $j+1$ switch summands and using $\lambda_j\geq j$ gives
\[
\mathsf S_j(\boldsymbol\lambda)\leq C_j\lambda_j^{j(j+1)}(1+\lambda_j)^j \leq C_j\lambda_j^{j(j+2)},
\]
which proves \eqref{eq:normalized-Muntz-switch}, and hence \eqref{eq:muntz-switch}.
\end{proof}

\subsection{Proofs of Theorem~\ref{thm:PU-PSECT} and Corollary~\ref{cor:PPU-maximal}} \label{subsec:PU-SECT-proof} We combine the comparison with the generalized moment curve, Proposition~\ref{prop:PU-to-moment}, with the switch estimate for pure powers, Proposition~\ref{prop:muntz-switch}; the corollary then follows from Theorem~\ref{thm:projective-SECT-maximal}.

\begin{proof}[Proof of Theorem~\ref{thm:PU-PSECT}] Suppose $\mathbf X\in\mathrm{PU}(d,\Lambda,R)$ has exponent vector $\boldsymbol\lambda$, and fix $2\leq j\leq d$.  By Proposition~\ref{prop:PU-to-moment},
\[
\mathsf S_j\left(\mathbf X^{[j]}\right) \leq R^{j(j-1)} \mathsf S_j\left(\mathbf M_{\boldsymbol\lambda}^{[j]}\right)=R^{j(j-1)} \mathsf S_j(\lambda_0,\ldots,\lambda_j),
\]
the equality holding since $\lambda_1=1$. Proposition~\ref{prop:muntz-switch}, applied to the truncated vector $(\lambda_0,\ldots,\lambda_j)$, bounds the second factor by $C_j(\lambda_j/\delta)^{j(j+2)}$, where $\delta$ is the least gap of the truncation. Truncation neither decreases the least gap nor increases the top exponent, so $\lambda_j/\delta\leq\lambda_d/\delta_{\boldsymbol\lambda} \leq\kappa_{\boldsymbol\lambda}$.  Taking the maximum over $2\leq j\leq d$ gives \eqref{eq:PU-SECT-sharp}, and \eqref{eq:kappa-vs-Lambda} then gives \eqref{eq:PU-SECT-polynomial}. The projective inclusion follows from Definition~\ref{def:power-unit-class} and Definition~\ref{def:projective-SECT}.
\end{proof}

\begin{proof}[Proof of Corollary~\ref{cor:PPU-maximal}] Apply \eqref{eq:PU-SECT-polynomial} in $\Pi$, or Remark~\ref{rem:PU-dimension-one} if $e=1$, together with Theorem~\ref{thm:projective-SECT-maximal}, and return to $\Pd$ by $d-e$ applications of Lemma~\ref{lem:hyperplane-reduction}. Since $K_e(\Lambda,R)\lesssim_e(\Lambda R)^{C_e}$ for $e\geq2$, while for $e=1$ one may take $K=1$, the resulting constant is $C_{d,\alpha,\varphi}(\Lambda R)^{C_d}$ after enlarging $C_d$.
\end{proof}

\section{Algebraic curves and projective power-unit charts} \label{sec:algebraic}This section supplies the geometric input required in Remark~\ref{rem:PU-interface}. An algebraic curve of degree $D$ is covered, up to polynomially many exceptional points, by polynomially many subarcs of projective power-unit curves with polynomially controlled parameters.

\medskip\noindent\emph{Conventions.} A \emph{real projective algebraic set} is a set
\[
Z=\left\{[x]\in\Pd: p_1(x)=\cdots=p_m(x)=0\right\}
\]
with $p_j\in\R[x_0,\ldots,x_d]$ homogeneous. Its Zariski closure $Z^{\C}\subset\mathbb C\mathrm P^d$ is invariant under complex conjugation. A \emph{complex projective algebraic set} is the set of common zeros in $\mathbb C\mathrm P^d$ of finitely many homogeneous polynomials, and an \emph{algebraic curve} is a complex projective algebraic set of dimension one. For a complex projective algebraic set $Y$, we write $Y(\R)\coloneqq Y\cap\Pd$ for its set of real points. Then $Z=Z^{\C}(\R)$, since $Z^{\C}$ lies in the complex zero set of the $p_j$. We write $\dim Z$ and $\deg Z$ for the dimension and degree of $Z^{\C}$. An algebraic set is \emph{irreducible} if it is not the union of two proper algebraic subsets. As in \cite{BN21}*{\S2.2.3}, the \emph{degree} of an irreducible complex algebraic set is its number of intersection points with a generic linear space of complementary dimension, and the degree of a reducible one is the sum of the degrees of its irreducible components, whatever their dimensions; in particular a point has degree $1$. The dimension and degree of a real algebraic subset of $\R^n$ are those of the projective closure of its Zariski closure in $\mathbb C^n$. A \emph{real algebraic arc} in $Z$ is a semialgebraic subset of $Z$ which is an embedded real-analytic submanifold of $\Pd$ homeomorphic to an open interval. The \emph{span} of a subset of $\Pd$ is the smallest real projective subspace containing it. For an algebraic curve $Y$, $\operatorname{Sm}(Y)$ denotes its set of smooth points. The algebraic complexity of a holomorphic function is that of \cite{BN21}*{\S2.2.3}: the least degree of an algebraic variety whose intersection with the product of the domain and $\C$ contains its graph as an analytic component, and, for a map, the maximum over its components; we use polynomially equivalent notions without further comment.

The dimension $d$ is fixed throughout. The symbols $C_d$ and $C_e$ denote constants depending only on $d$, respectively on $e\leq d$, and further subscripts, as in $C_{e,\mu}$, indicate further dependence; the value of every such constant may change from one occurrence to the next. We write $D^{C_d}$ for a quantity bounded by a power of $D\geq2$ with exponent depending only on $d$, so that a constant depending on $d$ times a fixed power of $D$ is bounded by $D^{C_d}$. Moreover, $\operatorname{poly}_{a,b,\ldots}(t_1,\ldots,t_k)$ denotes a quantity bounded by $C(1+t_1+\cdots+t_k)^C$ with $C\geq1$ depending only on the subscripts $a,b,\ldots$. Bounds $\delta^{C}$ in the degree $\delta$ of a curve are read with $\delta$ replaced by $\max\{\delta,2\}$.

\subsection{Statement of the decomposition} \label{subsec:algebraic-statement} The main result of the section is the following quantitative decomposition theorem.

\begin{theorem} \label{thm:algebraic-PPU-decomposition} For each $d\geq1$ there is $C_d<\infty$ with the following property. Let $D\geq2$ and let $Z\subset\Pd$ be a real projective algebraic set with
\[
\dim Z\leq1,\qquad \deg Z\leq D.
\]
There are a finite set $E\subset Z$ and real algebraic arcs $\mathcal C_1,\ldots,\mathcal C_M\subset Z$ such that
\begin{equation}\label{eq:algebraic-PPU-decomposition}
Z=E\cup\mathcal C_1\cup\cdots\cup\mathcal C_M,\qquad \#E+M\leq D^{C_d}.
\end{equation}
Moreover, there are parameters
\[
1\leq\Lambda_D,R_D\leq D^{C_d},
\]
depending only on $d$ and $D$, such that the following holds. For every $1\leq\nu\leq M$, let
\[
\Pi_\nu\coloneqq\operatorname{span}(\mathcal C_\nu)\subset\Pd,\qquad e_\nu\coloneqq\dim\Pi_\nu,\qquad 1\leq e_\nu\leq d.
\]
There is a projective isomorphism
\[
\phi_\nu:\Pi_\nu\longrightarrow\mathbb R\mathrm P^{e_\nu}
\]
and a curve
\[
\widetilde{\mathcal C}_\nu\in\mathrm{PPU}(e_\nu,\Lambda_D,R_D)
\]
such that
\begin{equation}\label{eq:algebraic-PPU-extension}
\phi_\nu(\mathcal C_\nu)\subset\widetilde{\mathcal C}_\nu.
\end{equation}
Since $\mathrm{PPU}(e,\Lambda,R)$ is invariant under $\operatorname{PGL}_{e+1}(\R)$, the statement holds for every projective isomorphism $\Pi_\nu\to\mathbb R\mathrm P^{e_\nu}$, after replacing $\widetilde{\mathcal C}_\nu$ by its projective image.
\end{theorem}

Each arc $\mathcal C_\nu$ of Theorem~\ref{thm:algebraic-PPU-decomposition} has to lie on a canonical curve, since the analysis of Section~\ref{sec:power-unit} applies to such curves, and we are free to use projective transformations to achieve this. In a given parametrization, the coordinate functions of $\mathcal C_\nu$ span a finite-dimensional space of functions, and projective transformations are changes of basis in this space; see the discussion in Subsection~\ref{subsec:intro-mechanism}. In general this space has no basis of monomials, but we construct, at the same time, a parametrization and a basis consisting of monomials times essentially constant analytic units, with distinct exponents. The weights defining the corresponding canonical coordinates are then powers times units, as Section~\ref{sec:power-unit} requires. This is done in four steps.

\begin{enumerate}[leftmargin=*]
\item \emph{Real arcs.} Proposition~\ref{prop:algebraic-monotone-arcs} cuts the real locus into polynomially many arcs avoiding the singular points of the curve, each a graph over one affine coordinate $x$.
\item \emph{Complex charts.} Lemma~\ref{lem:effective-scalar-lifting} splits each arc, up to finitely many points, into polynomially many pieces. Each piece is the image of the real segment of a punctured disc or annulus under a holomorphic map from the disc or annulus into the complex curve, which takes real values on that segment.
\item \emph{Monomial bases.} Proposition~\ref{prop:monomial-projective-frame} splits these pieces further and, on each, constructs the basis described above, using the monomial preparation of Proposition~\ref{prop:effective-cellular-preparation} and an induction on the dimension.
\item \emph{Power-unit weights.} Proposition~\ref{prop:monomial-frame-PU} proves that the weights are powers times units, and Lemma~\ref{lem:fill-annular-hole} then shows that the arc lies on a projective power-unit curve.
\end{enumerate}
Apart from B\'ezout's theorem and the Milnor--Thom bound, most of the algebraic geometric input of this section is contained in Lemma~\ref{lem:plane-model} below. From \cite{BN21} we use the Cellular Parametrization Theorem \cite{BN21}*{Theorem 8}, in Lemma~\ref{lem:real-cellular-cover}, and the one-dimensional monomialization lemma \cite{BN21}*{Lemma 46}, in Proposition~\ref{prop:effective-cellular-preparation}.

\subsection{Real arcs} \label{subsec:algebraic-differential-calculus} We write $\mathbb C\mathrm P^e$ for the complex projective space of complex lines in $\C^{e+1}$, with homogeneous coordinates $x_0,\ldots,x_e$, and $\operatorname{Sing}(Y)\coloneqq Y\setminus\operatorname{Sm}(Y)$. An irreducible algebraic curve is \emph{nondegenerate} if it is not contained in a projective hyperplane.

\begin{lemma}\label{lem:plane-model} Let $e\geq2$, let $Y\subset\mathbb C\mathrm P^e$ be a conjugation-invariant, nondegenerate, irreducible algebraic curve of degree $\delta$, and let $i\neq\ell$. There is a linear form $y$ on $\C^{e+1}$ with real coefficients such that the map
\[
p:\ Y\longrightarrow\mathbb C\mathrm P^2,\qquad p([x])\coloneqq[x_i:x_\ell:y(x)],
\]
is well defined on all of $Y$ and has the following properties.
\begin{enumerate}[itemsep=.4em]
\item[\emph{(i)}] Its image $C$ is an irreducible, conjugation-invariant curve. In the affine coordinates $(t,w)=(x_\ell/x_i,y/x_i)$,
\[
C\cap\C^2=\left\{(t,w)\in\C^2: f(t,w)=0\right\}
\]
for a real polynomial $f$, irreducible over $\C$, of degree at most $\delta$. The points of $Y$ with $x_i=0$ are mapped to the line at infinity.
\item[\emph{(ii)}] The fibers of $p$ have at most $\delta$ points, and $p$ maps $Y\cap p^{-1}(\operatorname{Sm}(C))$ biholomorphically onto $\operatorname{Sm}(C)$.
\item[\emph{(iii)}] $\partial_wf$ does not vanish identically on $C$, and
\[
\#\left\{(t,w)\in\C^2: f=\partial_wf=0\right\}\leq\delta^2,\qquad \#\operatorname{Sing}(Y)\leq\delta^3.
\]
\end{enumerate}
\end{lemma}

\begin{proof} For $\tau\in\C$, let $Y_\tau$ be the set of points of $Y\cap\{x_i\neq0\}$ with $x_\ell=\tau x_i$. By nondegeneracy, $Y$ is not contained in the hyperplane $\{x_\ell=\tau x_i\}$, so $\#Y_\tau\leq\delta$. Fix $\tau_0$ for which $k\coloneqq\#Y_{\tau_0}$ is maximal. The set $Y\cap\{x_i=x_\ell=0\}$ is finite, since $Y\not\subset\{x_i=0\}$. We choose a real linear form $y$ satisfying the conditions
\[
y\neq ax_i+bx_\ell\ \text{ for all } a,b\in\R,\qquad y\neq0\ \text{ on } Y\cap\{x_i=x_\ell=0\},\qquad y/x_i\ \text{ injective on } Y_{\tau_0}.
\]
Each condition fails only on finitely many proper linear subspaces of the space of real linear forms, the first one since $e\geq2$; hence almost every real linear form satisfies all three, and we fix one.

Since $x_i$, $x_\ell$ and $y$ are linear, $p([x])$ does not depend on the representative $x$, and the vector $(x_i,x_\ell,y(x))$ is nonzero, by the second condition if $x_i=x_\ell=0$. Hence $p$ is well defined on all of $Y$, and it commutes with conjugation since $y$ has real coefficients. It maps the points of $Y$ with $x_i=0$ to the line at infinity, and the others to the points with affine coordinates $(t,w)=(x_\ell/x_i,y/x_i)$.

Since $p$ is induced by the linear map $x\mapsto(x_i,x_\ell,y(x))$ from $\C^{e+1}$ to $\C^3$, the image $C\coloneqq p(Y)$ is an irreducible conjugation-invariant curve, so that $C\cap\C^2=\{f=0\}$ for an irreducible $f\in\R[t,w]$. A line of $\mathbb C\mathrm P^2$, and a point of $C$, pull back under $p$ into hyperplanes which do not contain $Y$, by the first condition and nondegeneracy. Hence lines meet $C$ in at most $\delta$ points, so that $m\coloneqq\deg f\leq\delta$, and the fibers of $p$ have at most $\delta$ points. If $\partial_wf$ vanished on $C$, then $f$ would divide $\partial_wf$, which has smaller degree in $w$, so $\partial_wf\equiv0$ and $f=t-c$; then $x_\ell=cx_i$ at all but finitely many points of $Y$, and $Y$ would lie in the hyperplane $\{x_\ell=cx_i\}$. Hence $\#\{f=\partial_wf=0\}\leq m(m-1)$ by B\'ezout.

Write $f=\sum_{j\leq k'}a_j(t)w^j$, with $a_{k'}\not\equiv0$. For $\tau\in\C$, the set $C_\tau$ of points of $C\cap\C^2$ with $t=\tau$ consists of the points $(\tau,w)$ with $w$ a root of $f(\tau,\cdot)$, which is not identically zero since $f$ is not a multiple of $t-\tau$. Hence $\#C_\tau\leq k'$, with equality unless $a_{k'}(\tau)=0$ or $f(\tau,\cdot)$ has a double root $w$, so that $(\tau,w)\in\{f=\partial_wf=0\}$. Both happen only for $\tau$ in a finite set $T$. The map $p$ sends $Y_\tau$ onto $C_\tau$, and it is injective on $Y_{\tau_0}$ by the choice of $y$, so that $k=\#C_{\tau_0}\leq k'$. For $\tau\notin T$ we therefore have
\[
k\leq k'=\#C_\tau\leq\#Y_\tau\leq k,
\]
so that $\#Y_\tau=\#C_\tau$ and $p$ is injective on $Y_\tau$. Two points of $Y$ with the same image lie in the same set $Y_\tau$ or in the finite set $Y\cap\{x_i=0\}$, so only finitely many points of $Y$ share their image with another point.

It remains to prove the second assertion of (ii). Let $F\subset Y$ consist of the singular points of $Y$ and of the points $q\in Y$ with $p(q)=p(q')$ for some $q'\in Y\setminus\{q\}$; as shown above, $F$ is finite, and $p$ is injective on $Y\setminus F$.

\begin{figure}[ht]
\centering
\begin{tikzpicture}[>=stealth]
\draw[thick] (0,4) .. controls (2,4.8) and (3.5,3.6) .. (5,4.1) .. controls (6.5,4.6) and (8,3.8) .. (10,4.2);
\node at (-0.4,4) {$Y$};
\node[red] at (1.5,4.30) {$\times$};
\node[red] at (8.3,4.10) {$\times$};
\node[red,above] at (1.5,4.45) {\small$F$};
\node[red,above] at (8.3,4.25) {\small$F$};
\fill (5,4.1) circle (1.5pt);
\draw[dashed] (5,4.1) circle (0.75);
\node[above] at (5,4.1) {\small$h(Q)$};
\draw[thick] (0,0) .. controls (2,0.8) and (3.5,-0.4) .. (5.5,0.1) .. controls (7,0.5) and (8.5,-0.2) .. (10,0.2);
\node at (-0.4,0) {$C$};
\draw[dashed] (6,0.22) circle (1.1);
\node at (6.90,1.35) {\small$\mathcal N$};
\fill (6,0.19) circle (1.5pt);
\node[below right] at (6,0.2) {\small$P$};
\fill (5.2,0.04) circle (1.5pt);
\node[below] at (5.2,0) {\small$Q$};
\fill[red] (1.5,0.28) circle (1.5pt);
\fill[red] (8.6,0.08) circle (1.5pt);
\node[red,below] at (1.5,0.2) {\small$p(F)$};
\node[red,below] at (8.6,-0.00) {\small$p(F)$};
\draw[->] (5,3.55) -- node[right] {$p$} (5.2,0.3);
\draw[->,dashed] (5.35,0.35) to[bend right=45] node[pos=0.55,right] {$h$} (5.25,3.6);
\draw (1.2,2.2) circle (0.5);
\fill (1.2,2.2) circle (1.2pt);
\node[below] at (1.2,1.65) {\small$\Delta'$};
\node[above] at (1.2,2.22) {\scriptsize$z_0$};
\node[above] at (9,2.22) {\scriptsize$0$};
\draw[->] (1.55,2.55) .. controls (2.8,3.0) .. node[above left] {$\psi$} (4.55,4.05);
\draw (9,2.2) circle (0.5);
\fill (9,2.2) circle (1.2pt);
\node[right] at (9.55,2.2) {\small$\Delta$};
\draw[->] (8.65,1.85) .. controls (8.1,1.3) .. node[above right, xshift=-8pt] {$\Phi$} (6.95,0.85);
\draw[->,dashed] (1.75,2.2) -- node[pos=0.3,below] {\small$\varphi=\Phi^{-1}\circ p\circ\psi$} (8.45,2.2);
\end{tikzpicture}
\caption{The maps in the proof of Lemma~\ref{lem:plane-model}.}
\label{fig:plane-model-inverse}
\end{figure}
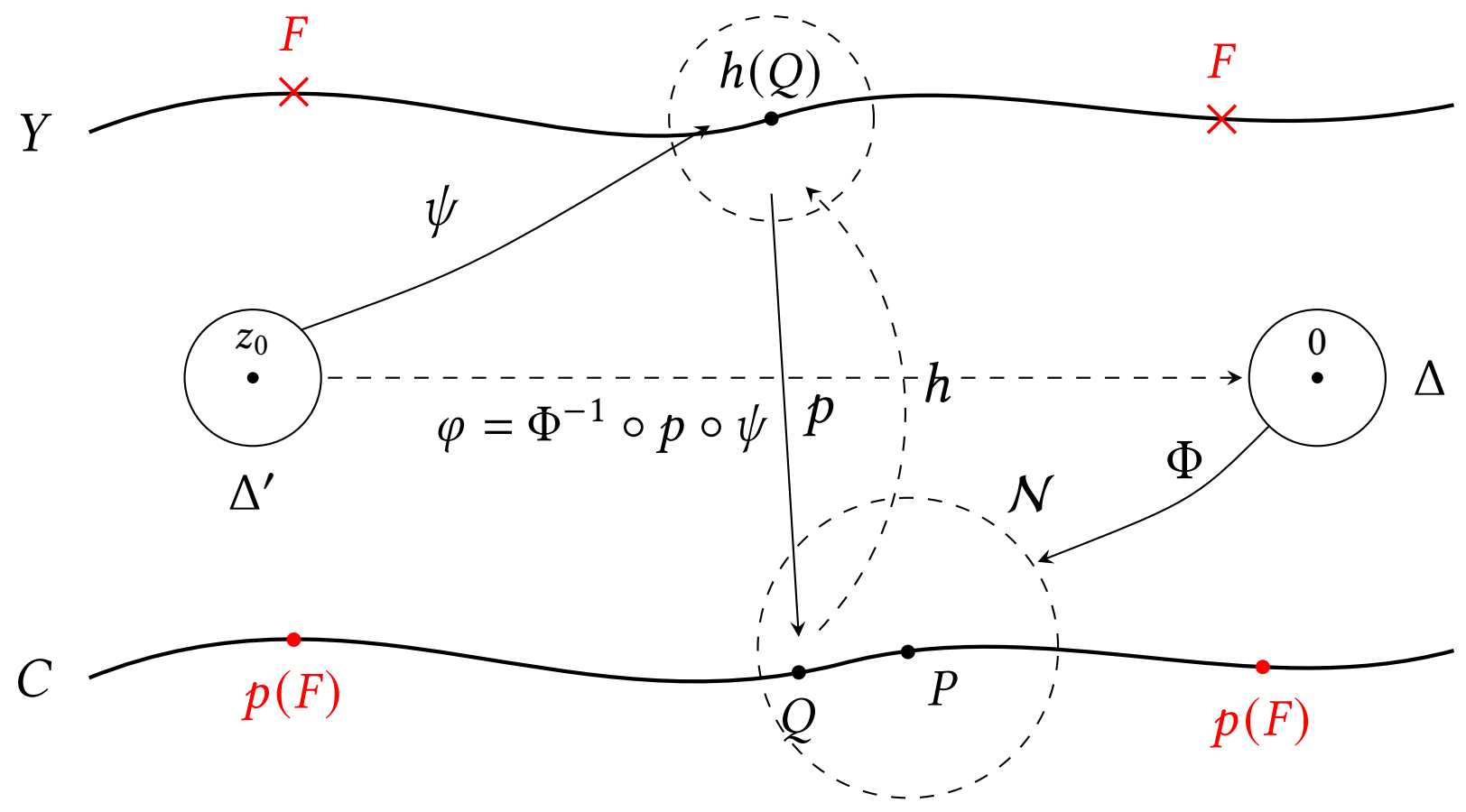

Fix $P\in\operatorname{Sm}(C)$, and a biholomorphic map $\Phi$ from a disc $\Delta$ centered at $0$ onto a neighborhood $\mathcal N$ of $P$ in $C$, with $\Phi(0)=P$, such that $\mathcal N\setminus\{P\}$ does not meet the finite set $p(F)$; see Figure~\ref{fig:plane-model-inverse}. Each $Q\in\mathcal N\setminus\{P\}$ then has a single preimage $h(Q)\notin F$, a smooth point of $Y$. If $\psi$ is a biholomorphic map from a disc $\Delta'$ onto a neighborhood of $h(Q)$ in $Y$, with $\psi(z_0)=h(Q)$, the holomorphic function $\varphi\coloneqq\Phi^{-1}\circ p\circ\psi$ satisfies $\varphi'(z_0)\neq0$, since otherwise the points of $\mathcal N$ near $Q$ would have at least two preimages near $h(Q)$. Hence $h\circ\Phi=\psi\circ\varphi^{-1}$ near $\Phi^{-1}(Q)$, and $h\circ\Phi$ is holomorphic on $\Delta\setminus\{0\}$.

As $\Delta\ni w\to0$, the points $h(\Phi(w))$ approach the finite set $p^{-1}(P)$, by compactness of $Y$ and continuity of $p$. After shrinking $\Delta$, $h\circ\Phi$ is therefore bounded in an affine chart, and by Riemann's removable singularity theorem it extends to a holomorphic map $g:\Delta\to Y$ with $p\circ g=\Phi$. If $q\in Y\setminus F$ and $p(q)\in\mathcal N$, then $q=g(\Phi^{-1}(p(q)))$, since $q$ is the only preimage of $p(q)$. As $F$ is finite and $Y$ has no isolated points, every point of $Y$ near $p^{-1}(P)$ is a point $q\in Y\setminus F$ with $p(q)\in\mathcal N$, or a limit of such points. Hence $p^{-1}(P)=\{g(0)\}$, and $Y$ coincides with $g(\Delta)$ near $g(0)$. Since $p\circ g=\Phi$ is an injective immersion, so is $g$; thus $g(0)\in\operatorname{Sm}(Y)$, and $p$ is a biholomorphism near $g(0)$, with inverse $g\circ\Phi^{-1}$.

Consequently, $p$ maps $\operatorname{Sing}(Y)$ into $\operatorname{Sing}(C)$, and since its fibers have at most $\delta$ points, $\#\operatorname{Sing}(Y)\leq\delta\,\#\operatorname{Sing}(C)$. The singular points of $C$ in the affine chart lie in $\{f=\partial_wf=0\}$, which has at most $m(m-1)$ points, and $C$ meets the line at infinity in at most $m$ points. Hence $\#\operatorname{Sing}(Y)\leq\delta(m(m-1)+m)\leq\delta^3$.
\end{proof}

\begin{remark}The set $\{f=\partial_wf=0\}$ in (iii) may contain smooth points of $C$ as well as real singular points of $C$ without real preimages. For instance, if $Y\subset\mathbb C\mathrm P^3$ is the twisted cubic, consisting of the points $[1:s:s^2:s^3]$, $s\in\C$, together with $[0:0:0:1]$, $i=0$, $\ell=2$ and $y=x_1+x_3$, then $(t,w)=(s^2,s+s^3)$ and $f=w^2-t(1+t)^2$. Here $\{f=\partial_wf=0\}$ consists of $(0,0)$, a smooth point of $C$ where $t$ is not a local coordinate, and of the node $(-1,0)$, a real point of $C$ whose preimages $s=\pm i$ are not real.
\end{remark}

We now use the plane model of Lemma~\ref{lem:plane-model} to cut the real points of $Y$ into polynomially many arcs, each a graph over an affine coordinate.

\begin{proposition} \label{prop:algebraic-monotone-arcs} Let $e\geq2$ and let $Y\subset\mathbb C\mathrm P^e$ be a conjugation-invariant, nondegenerate, irreducible algebraic curve of degree $\delta$. Then $Y(\R)$ is the union of at most $\delta^{C_e}$ points and at most $\delta^{C_e}$ real algebraic arcs $\Gamma\subset\operatorname{Sm}(Y)$. Each arc $\Gamma$ lies in the unit cube of a standard affine chart $\{x_i\neq0\}$, and one of the affine coordinates $x=x_\ell/x_i$ of that chart is injective on $\Gamma$.
\end{proposition}

\begin{proof} For indices $0\leq a<b\leq e$ of the homogeneous coordinates, nondegeneracy shows that $Y$ is contained in neither of the hyperplanes $\{x_a=\pm x_b\}$, so by B\'ezout the set $S_\ast$ of points of $Y(\R)$ at which $|x_a|=|x_b|$ for some $a\neq b$ satisfies $\#S_\ast\leq e(e+1)\delta$. Every point of $Y(\R)\setminus S_\ast$ lies in exactly one of the sets
\[
\mathcal R_i\coloneqq\left\{|x_i|>|x_a|\text{ for every }a\neq i\right\},\qquad 0\leq i\leq e,
\]
and $\mathcal R_i$ is contained in the unit cube of the chart $\{x_i\neq0\}$. Fix $i$ and $\ell\neq i$, put $x\coloneqq x_\ell/x_i$, and let $p$, $C$ and $f$ be given by Lemma~\ref{lem:plane-model}. The set
\[
\Sigma_i\coloneqq\left\{(t,w)\in\R^2: f=\partial_wf=0\right\}
\]
has at most $\delta^2$ points, so the set $S_i$ of points of $Y(\R)\cap\mathcal R_i$ mapped into $\Sigma_i$ has at most $\delta^3$ points.

Let $W_i\coloneqq(Y(\R)\cap\mathcal R_i)\setminus(S_i\cup S_\ast)$ and $q\in W_i$. Then $p(q)$ is a real point of $C$ with $\partial_wf\neq0$, so $p(q)\in\operatorname{Sm}(C)$ and, by the implicit function theorem, near $p(q)$ the curve $C$ is a graph $w=h(t)$, with $h$ holomorphic and real for real $t$. Since $p$ is a biholomorphism over $\operatorname{Sm}(C)$ commuting with conjugation, $q\in\operatorname{Sm}(Y)$ and near $q$ the set $Y(\R)$ is the curve $x\mapsto p^{-1}(x,h(x))$, as $t\circ p=x$. Thus $W_i$ is an embedded one-dimensional real-analytic submanifold of $\mathbb R\mathrm P^e$, contained in $\operatorname{Sm}(Y)$, on which $x$ is a local coordinate. A connected component $\Gamma$ of $W_i$ is not compact, since $x$ would attain its maximum on it; hence $\Gamma$ is homeomorphic to an open interval and $x$ is strictly monotone on it. Components of semialgebraic sets are semialgebraic \cite{BPR}, so $\Gamma$ is a real algebraic arc.

It remains to count the components. Let
\[
K_i\coloneqq\left\{(t,w)\in\R^2: f(t,w)=0\right\}\setminus\left(\Sigma_i\cup p(S_\ast)\right).
\]
Since $\partial_wf\neq0$ on $K_i$, we have $K_i\subset\operatorname{Sm}(C)$, where $p$ is a biholomorphism commuting with conjugation. Hence $p$ maps $W_i$ into $K_i$, as $W_i$ avoids $S_i$ and $S_\ast$, and the preimage of a component of $K_i$ is a connected subset of $Y(\R)\setminus S_\ast$, hence lies in a single region $\mathcal R_j$, and in $W_i$ if $j=i$. As $K_i$ is a one-dimensional manifold, the preimages lying in $W_i$ are open and cover it, so $W_i$ has at most as many components as $K_i$. By the Milnor--Thom theorem \cites{Milnor64,BPR} the real zero set of $f$ has at most $\delta(2\delta-1)$ components, and removing a point $P$ adds at most as many components as the number of points in which a small circle centered at $P$ meets $C$, which is at most $2\delta$ by B\'ezout, since $C$ passes through $P$ and is therefore not the circle. Since $\#(\Sigma_i\cup p(S_\ast))\leq\delta^2+e(e+1)\delta$, the set $W_i$ has at most $\delta^{C_e}$ components. Every point of $Y(\R)$ lies in $S_\ast\cup S_0\cup\cdots\cup S_e$ or in one of the sets $W_0,\ldots,W_e$, and summing over $i$ proves the proposition.
\end{proof}

\subsection{Cells} \label{subsec:algebraic-cellular-preparation} For $0\leq r<R$ and $\mu>1$ write
\[
A(r,R)\coloneqq\{z\in\C: r<|z|<R\},\qquad \mathcal L^{[\mu]}\coloneqq A\left(\frac r\mu,\mu R\right)\quad\text{if }\mathcal L=A(r,R),
\]
so that $A(0,R)$ is a punctured disc. For a disc $\mathcal B=\{|z|<R\}$ we write $\mathcal B^{[\mu]}\coloneqq\{|z|<\mu R\}$; replacing it by the punctured disc $A(0,R)$ does not change its positive real part. We call discs, punctured discs and annuli centered at $0$ \emph{cells}; together with points, they are the one-dimensional real cells of \cite{BN21}*{Definition 10}. On cells we write
\[
\mathfrak D\coloneqq z\frac{\d}{\d z}.
\]
A \emph{real cell} is a translate $\theta+\mathcal B$ of a cell $\mathcal B$, with $\theta\in\R$. Its extensions are $\theta+\mathcal B^{[\mu]}$, and its \emph{segment} is the interval $\theta+(\mathcal B\cap(0,\infty))$. The translation $z\mapsto\theta+z$ carries $\mathcal B$, its extensions and $\mathcal B\cap(0,\infty)$ onto $\theta+\mathcal B$, its extensions and its segment. It preserves complexity and reality on the real axis, and it reduces each statement below about one real cell to the case $\theta=0$. A holomorphic function $f$ on a connected open set invariant under conjugation, real on a nondegenerate interval, is real at all real points of the set, by the identity theorem applied to $f-\overline{f(\bar\cdot)}$. In one variable, the $\{\sigma\}$-extensions of \cite{BN21}*{\S2.1, (29)} are extensions in this sense: the $\{\sigma\}$-extension of a disc, a punctured disc or an annulus $\mathcal B$ is $\mathcal B^{[\lambda]}$, where
\begin{equation} \label{eq:sigma-extension}
\lambda=\e^{\pi^2/(2\sigma)}\quad\text{(punctured discs, annuli)},\qquad \lambda-\lambda^{-1}=\frac{2\pi}{\sigma}\quad\text{(discs)}.
\end{equation}
The parameter $\sigma$ is used only when citing \cite{BN21}.

To parametrize the arcs of Proposition~\ref{prop:algebraic-monotone-arcs} holomorphically, we need polynomially many cells around a real interval which, together with some room around them, avoid a finite set of points of the complex plane. The following lemma, the one-variable form of the Cellular Parametrization Theorem of \cite{BN21}, provides them. It is applied in Lemma~\ref{lem:effective-scalar-lifting} and in Proposition~\ref{prop:effective-cellular-preparation} below.

\begin{lemma} \label{lem:real-cellular-cover} Let $P\in\R[z]$ be a nonzero polynomial of degree at most $\beta$, and $\nu>1$ be fixed. Let $\mathcal L=A(r,R)$ be a punctured disc or an annulus, and let $\lambda\geq \e^2$. Then there are
\begin{enumerate}[itemsep=.4em]
\item[\emph{(a)}] a finite set $E\subset\mathcal L\cap(0,\infty)$,
\item[\emph{(b)}] segments $I_\alpha$ of real cells $\theta_\alpha+\mathcal B_\alpha$,
\end{enumerate}
such that
\begin{enumerate}[itemsep=.4em]
\item[\emph{(i)}] $P$ has no zeros on $\theta_\alpha+\mathcal B_\alpha^{[\lambda]}\subset\mathcal L^{[\nu]}$ if $\mathcal B_\alpha$ is a punctured disc or an annulus, or on $\theta_\alpha+\mathcal B_\alpha^{[\log\lambda]}\subset\mathcal L^{[\nu]}$ if $\mathcal B_\alpha$ is a disc;
\item[\emph{(ii)}] the $I_\alpha$ cover $\mathcal L\cap(0,\infty)$ up to $E$:
\[
\mathcal L\cap(0,\infty)=E\cup\bigcup_\alpha I_\alpha,\qquad I_\alpha\subset\mathcal L\cap(0,\infty);
\]
\item[\emph{(iii)}] $\#E+\#\{\alpha\}\leq\operatorname{poly}_{\nu}(\beta,\log\lambda)$.
\end{enumerate}
\end{lemma}

\begin{proof} Put
\[
Q(z)\coloneqq P(z)(z-r)(z-R),
\]
a nonzero real polynomial of degree at most $\beta+2$, bounded on $\mathcal L^{[\nu]}$. By \eqref{eq:sigma-extension}, $\mathcal L^{[\nu]}$ is the $\{\rho\}$-extension of $\mathcal L$ for $\rho\coloneqq\pi^2/(2\log\nu)$. The Cellular Parametrization Theorem \cite{BN21}*{Theorem 8}, in its real algebraic form, applied to $Q$ with target $\mathcal L^{[\nu]}$ and extension parameter
\[
\sigma\coloneqq\frac{\pi^2}{2\log\lambda},
\]
gives a real cellular cover of $\mathcal L$ by $\operatorname{poly}_\nu(\beta,\sigma^{-1})=\operatorname{poly}_\nu(\beta,\log\lambda)$ maps $f_j$, compatible with $Q$, that is \cite{BN21}*{Definition 15}, such that $Q\circ f_j$ vanishes identically or nowhere. Each $f_j$ is defined on the $\{\sigma\}$-extension of a one-dimensional real cell $\mathcal B_j$ of \cite{BN21}, which is a point, a disc, a punctured disc or an annulus. In the last three cases, by \eqref{eq:sigma-extension}, this extension is $\mathcal B_j^{[\lambda_j]}$, with $\lambda_j=\lambda$, or, if $\mathcal B_j$ is a disc, $\lambda_j-\lambda_j^{-1}=\frac4\pi\log\lambda$, so that $\lambda_j\geq\log\lambda$. Each $f_j$ is prepared \cite{BN21}*{Theorem 8, Definition 11}, that is, $f_j(w)=\theta_j+w^{q_j}$ with an integer $q_j\geq1$, and $\theta_j\in\R$ since $f_j$ is real. As $w\mapsto w^{q_j}$ maps $\mathcal B_j\cap(0,\infty)$ onto the positive real part of the image cell $\mathcal B_j'$, and $\mathcal B_j^{[\lambda_j]}$ onto $(\mathcal B_j')^{[\lambda_j^{q_j}]}\supset(\mathcal B_j')^{[\lambda_j]}$, we may replace $\mathcal B_j$ by $\mathcal B_j'$, and $f_j$ by the translate $w\mapsto\theta_j+w$, whose values on $(\mathcal B_j')^{[\lambda_j]}$ are values of the original $f_j$ on $\mathcal B_j^{[\lambda_j]}$; in particular, compatibility is preserved. Being a real cover means \cite{BN21}*{Definition 13} that $\mathcal L\cap(0,\infty)$ is covered by the points $f_j(\mathcal B_j)$, for $\mathcal B_j$ a point, and by the intervals $f_j(\mathcal B_j\cap(0,\infty))$ for the other $j$. For the latter, $\theta_j+\mathcal B_j$ is a real cell, with $\theta_j+\mathcal B_j^{[\lambda_j]}\subset\mathcal L^{[\nu]}$, and $f_j(\mathcal B_j\cap(0,\infty))$ is its segment. By compatibility, $Q$ vanishes identically or nowhere on $\theta_j+\mathcal B_j^{[\lambda_j]}$, and the first alternative is excluded since $Q\not\equiv0$. Hence $P$ has no zeros there, and the segment of $\theta_j+\mathcal B_j^{[\lambda_j]}$ contains neither $r$ nor $R$. As it contains $f_j(\mathcal B_j\cap(0,\infty))$, the latter lies either in $(r,R)=\mathcal L\cap(0,\infty)$ or in its complement. The points $f_j(\mathcal B_j)$ in $(r,R)$, for $\mathcal B_j$ a point, form $E$, and the segments in $(r,R)$ are the $I_\alpha$. Together they cover $(r,R)$.
\end{proof}

We return to the arcs of Proposition~\ref{prop:algebraic-monotone-arcs}. The next lemma applies Lemma~\ref{lem:real-cellular-cover} to a polynomial vanishing where the curve is not locally a union of holomorphic graphs over $x$, and obtains holomorphic parametrizations of the arcs on punctured discs and annuli, the domains on which the preparation below takes place.

\begin{lemma} \label{lem:effective-scalar-lifting} Fix $\mu>1$. Let $Y$ be as in Proposition~\ref{prop:algebraic-monotone-arcs}, let $\Gamma$ be one of the arcs given by that proposition, and let $x=x_\ell/x_i$ be the affine coordinate over which $\Gamma$ is a graph. We extend $x$, by the same formula, to all points of $Y$ at which $x_i\neq0$:
\[
 Y\cap\{x_i\neq0\}\ni q \longmapsto x(q)\coloneqq \frac{q_\ell}{q_i}\in \C.
\] Then there are
\begin{enumerate}[itemsep=.4em]
\item[\emph{(a)}] a finite set $E_\Gamma\subset\Gamma$,
\item[\emph{(b)}] punctured discs or annuli $\mathcal B_\alpha$ and holomorphic maps $\gamma_\alpha:\mathcal B_\alpha^{[\mu]}\to Y\cap\{x_i\neq0\}$,
\end{enumerate}
such that
\begin{enumerate}[itemsep=.4em]
\item[\emph{(i)}] the positive real parts of the $\mathcal B_\alpha$ parametrize $\Gamma$ up to $E_\Gamma$:
\begin{equation} \label{eq:effective-scalar-lifting-cover}
\Gamma=E_\Gamma\cup\bigcup_\alpha\gamma_\alpha\left(\mathcal B_\alpha\cap(0,\infty)\right),\qquad \gamma_\alpha\left(\mathcal B_\alpha\cap(0,\infty)\right)\subset\Gamma;
\end{equation}
\item[\emph{(ii)}] in the coordinate $x$, each $\gamma_\alpha$ is a power map:
\[
x\circ\gamma_\alpha(z)=\theta_\alpha+z^{s_\alpha},\qquad \theta_\alpha\in\R,\quad 1\leq s_\alpha\leq\delta;
\]
\item[\emph{(iii)}] the affine coordinates of $\gamma_\alpha$ are bounded algebraic holomorphic functions on $\mathcal B_\alpha^{[\mu]}$ of complexity $\delta^{C_e}$, real on its positive real part;
\item[\emph{(iv)}] $\#E_\Gamma+\#\{\alpha\}\leq\delta^{C_{e,\mu}}$.
\end{enumerate}
\end{lemma}

\begin{proof} We begin by showing that, away from the zeros of a resultant, the curve is locally a union of holomorphic graphs over the coordinate $x$. Let $p$ and $f$ be given by Lemma~\ref{lem:plane-model}, so that $t=x$ on $Y\cap\{x_i\neq0\}$, and write $f=\sum_{j\leq k}a_j(t)w^j$ with $a_k\not\equiv0$, so that $1\leq k\leq\delta$. We set
\[
\Delta\coloneqq a_k\operatorname{Res}_w(f,\partial_wf)\in\R[t],
\]
where the resultant $\operatorname{Res}_w(f,\partial_wf)$ is a polynomial in $t$, computed from the coefficients of $f$ and $\partial_wf$ as polynomials in $w$, which vanishes at $t=\tau$ exactly when $f(\tau,\cdot)$ and $\partial_wf(\tau,\cdot)$ have a common root or $a_k(\tau)=0$. Since $f$ is irreducible and $\partial_wf=\sum_{j\leq k}ja_j(t)w^{j-1}$ is nonzero of degree $k-1$ in $w$, as $k\geq1$, $f$ does not divide $\partial_wf$, so the two have no common factor. Hence $\Delta\not\equiv0$, and $\deg\Delta\leq3\delta^2$. Let $U\subset\C$ be open and connected, with $\Delta\neq0$ on $U$. For $t\in U$ the polynomial $f(t,\cdot)$ has degree $k$ and $k$ simple roots, which depend holomorphically on $t$; thus
\[
\left\{(t,w)\in U\times\C: f(t,w)=0\right\}
\]
is a $k$-sheeted covering of $U$, contained in $\operatorname{Sm}(C)$. By Lemma~\ref{lem:plane-model}, $p$ maps
\[
Y_U\coloneqq\left\{q\in Y\cap\{x_i\neq0\}: x(q)\in U\right\}
\]
biholomorphically onto this covering, with $t\circ p=x$. Hence
\begin{equation} \label{eq:lifting-covering}
x:Y_U\longrightarrow U\ \text{ is a covering with at most $\delta$ sheets, and } Y_U\subset\operatorname{Sm}(Y).
\end{equation}

Next we choose cells on the line. The interval $x(\Gamma)$ is bounded, say $x(\Gamma)=(a,b)$. Apply Lemma~\ref{lem:real-cellular-cover} to $A(0,b-a)$, with $\nu=2$, the polynomial $\Delta(a+\cdot)$ and $\lambda\coloneqq \e^{C_\mu\delta}$, where $C_\mu$ is large, and translate by $a$. This gives a finite set $E_0$ and real cells $\theta_\alpha+\mathcal A_\alpha$ such that
\[
(a,b)=E_0\cup\bigcup_\alpha I_\alpha,\qquad I_\alpha\coloneqq\theta_\alpha+(\mathcal A_\alpha\cap(0,\infty)),\qquad \Delta\neq0\ \text{ on }\ U_\alpha,
\]
where $U_\alpha\coloneqq\theta_\alpha+\mathcal A_\alpha^{[\lambda]}$, or $U_\alpha\coloneqq\theta_\alpha+\mathcal A_\alpha^{[\log\lambda]}$ if $\mathcal A_\alpha$ is a disc. The preimages in $\Gamma$ of the points of $E_0$ form $E_\Gamma$. Fix $\alpha$ and put
\[
\Gamma_\alpha\coloneqq\left\{q\in\Gamma: x(q)\in I_\alpha\right\}.
\]
Then $\Gamma_\alpha$ is a connected subset of $Y_{U_\alpha}$; let $Y_\alpha$ be the component of $Y_{U_\alpha}$ containing it. By \eqref{eq:lifting-covering}, $x:Y_\alpha\to U_\alpha$ is a connected covering with $s_\alpha\leq\delta$ sheets.

We now lift these cells. If $\mathcal A_\alpha$ is a disc, then $U_\alpha$ is simply connected and $s_\alpha=1$. Then
\[
\gamma_\alpha(z)\coloneqq(x|_{Y_\alpha})^{-1}(\theta_\alpha+z)
\]
is defined on $\mathcal A_\alpha^{[\log\lambda]}$, which contains $\mathcal A_\alpha^{[\mu]}$ since $\log\lambda=C_\mu\delta>\mu$, and we let $\mathcal B_\alpha\coloneqq\mathcal A_\alpha\setminus\{0\}$. Otherwise
\[
\mathcal A_\alpha=A(\rho_1,\rho_2),\qquad U_\alpha=\theta_\alpha+A(\rho_1/\lambda,\lambda\rho_2).
\]
The connected coverings of an annulus with $s$ sheets are, up to isomorphism, the power maps $\zeta\mapsto\zeta^s$; hence there is a biholomorphism
\[
\gamma_\alpha:\mathcal B_\alpha^{[\lambda^{1/s_\alpha}]}\longrightarrow Y_\alpha,\qquad x\circ\gamma_\alpha(\zeta)=\theta_\alpha+\zeta^{s_\alpha},\qquad \mathcal B_\alpha\coloneqq A\left(\rho_1^{1/s_\alpha},\rho_2^{1/s_\alpha}\right).
\]
After a rotation $\zeta\mapsto\omega\zeta$ with $\omega^{s_\alpha}=1$, $\gamma_\alpha$ maps $\mathcal B_\alpha\cap(0,\infty)$ onto $\Gamma_\alpha$, since $\Gamma_\alpha$ is a lift of the interval $I_\alpha$. Moreover
\[
\lambda^{1/s_\alpha}\geq\lambda^{1/\delta}=\e^{C_\mu}>\mu,
\]
so $\gamma_\alpha$ is defined on $\mathcal B_\alpha^{[\mu]}$.

In both cases
\[
x\circ\gamma_\alpha(z)=\theta_\alpha+z^{s_\alpha},\qquad \gamma_\alpha\left(\mathcal B_\alpha\cap(0,\infty)\right)=\Gamma_\alpha\subset\Gamma,
\]
so \eqref{eq:effective-scalar-lifting-cover} holds. The graph of each affine coordinate of $\gamma_\alpha$ lies on the projection of the algebraic curve
\[
\left\{(z,q)\in\C\times Y: x(q)=\theta_\alpha+z^{s_\alpha}\right\},
\]
of degree $\operatorname{poly}(\delta)$, which gives the complexity bound. The coordinates are real on $\mathcal B_\alpha\cap(0,\infty)$ since $\Gamma$ is real, hence on $\mathcal B_\alpha^{[\mu]}\cap(0,\infty)$ by the identity theorem, as noted in Subsection~\ref{subsec:algebraic-cellular-preparation}. They are also bounded on $\mathcal B_\alpha^{[\mu]}$: if $\mathcal A_\alpha$ is a disc or an annulus, the closure of $\mathcal B_\alpha^{[\mu]}$ is a compact subset of the domain of $\gamma_\alpha$; if it is a punctured disc, each coordinate is algebraic and single-valued, hence meromorphic at $0$, and has no pole there, being bounded on $\mathcal B_\alpha\cap(0,\infty)$ since $\Gamma$ lies in the unit cube.

Finally, Lemma~\ref{lem:real-cellular-cover} gives
\[
\#E_\Gamma+\#\{\alpha\}\leq\operatorname{poly}(\delta,\log\lambda)=\delta^{C_{e,\mu}}.\qedhere
\]
\end{proof}

\subsection{Preparation and monomial frames} \label{subsec:annular-projective-monomial-frames} On the cells of Lemma~\ref{lem:effective-scalar-lifting}, the coordinates of an arc are bounded holomorphic functions. Proposition~\ref{prop:effective-cellular-preparation} below subdivides each cell further, so that on the new cells such functions take the form $z^n$ times a nearly constant unit. Proposition~\ref{prop:monomial-projective-frame} then finds a basis of the real span of $1$ and the affine coordinates of a chart $\gamma_\alpha$ whose elements all have this form, with pairwise different exponents $n$.

\begin{definition}\label{def:cellularly-prepared}Let $\mathcal U=A(r,R)$. A holomorphic function $f$ on $\mathcal U$ is \emph{cellularly prepared on $\mathcal U$ with exponent $n\in\Z$ and accuracy $\eta$} if either $f\equiv0$ or
\[
f(z)=z^n u(z),
\]
where $u$ is nonvanishing on $\mathcal U$, admits a single-valued logarithm there, and
\[
\diam_{\mathcal U}\log u\leq\eta.
\]
A finite family $\mathcal F$ is \emph{cellularly prepared on $\mathcal U$ with accuracy $\eta$} if each $f\in\mathcal F$ is cellularly prepared there, with some exponent $n_f\in\Z$.
\end{definition}

A function $f$ is \emph{cellularly prepared on a real cell $\theta+\mathcal U$}, with exponent $n$ and accuracy $\eta$, if $f(\theta+\cdot)$ is so on $\mathcal U$, and similarly for families.

\begin{definition}\label{def:prepared-cellular-cover} Fix $\mu,\nu>1$ and $0<\eta<1$. Let $\mathcal L$ be a punctured disc or an annulus, and let $\mathcal F$ be a finite family of algebraic holomorphic functions on $\mathcal L^{[\nu]}$. A \emph{$(\mu,\eta)$-prepared cellular cover of $\mathcal L$ for $\mathcal F$} consists of a finite set $E\subset\mathcal L\cap(0,\infty)$ and segments $I_\alpha$ of real cells $\theta_\alpha+\mathcal B_\alpha$, with $\mathcal B_\alpha$ punctured discs or annuli and $\theta_\alpha+\mathcal B_\alpha^{[\mu]}\subset\mathcal L^{[\nu]}$, such that
\[
\mathcal L\cap(0,\infty)=E\cup\bigcup_\alpha I_\alpha,\qquad I_\alpha\subset\mathcal L\cap(0,\infty),
\]
and $\mathcal F$ is cellularly prepared on each $\theta_\alpha+\mathcal B_\alpha^{[\mu]}$ with accuracy $\eta$.
\end{definition}

The following proposition is a preparation result that applies to any bounded algebraic holomorphic functions on a cell. In Proposition~\ref{prop:monomial-projective-frame} below, it will be applied to the coordinates of the charts and to functions derived from them.

\begin{proposition}\label{prop:effective-cellular-preparation} Fix $m\geq1$ and $\mu,\nu>1$. Let $\mathcal L=A(r,R)$ be a punctured disc or an annulus, and let $f_1,\ldots,f_m$ be bounded algebraic holomorphic functions on $\mathcal L^{[\nu]}$, of complexity at most $\beta$ and real on its positive real part. Then for each $0<\eta<1$ there is a $(\mu,\eta)$-prepared cellular cover of $\mathcal L$ for $\{f_1,\ldots,f_m\}$, with exceptional set $E$ and segments $I_\alpha$ of real cells $\theta_\alpha+\mathcal B_\alpha$, such that
\begin{enumerate}[itemsep=.4em]
\item[\emph{(i)}] $\#E+\#\{\alpha\}\leq\operatorname{poly}_{m,\mu,\nu}(\beta,\eta^{-1})$;
\item[\emph{(ii)}] the preparation exponents are bounded independently of $\eta$: for $f_j\not\equiv0$,
\[
f_j(\theta_\alpha+z)=z^{n_{\alpha,j}}u_{\alpha,j}(z)\quad\text{on }\mathcal B_\alpha^{[\mu]},\qquad |n_{\alpha,j}|\leq\operatorname{poly}_{m,\mu,\nu}(\beta);
\]
\item[\emph{(iii)}] all units $u_{\alpha,j}$ are real on the real axis.
\end{enumerate}
\end{proposition}

\begin{proof} For each $f_j\not\equiv0$ choose an irreducible polynomial $P_j\in\C[z,w]$ of degree $\operatorname{poly}_m(\beta)$ such that $P_j(z,f_j(z))=0$ on $\mathcal L^{[\nu]}$. Then $P_j(z,0)\not\equiv0$, since otherwise $P_j$ would be a multiple of $w$ and $f_j$ would vanish identically, and $f_j(z)=0$ implies $P_j(z,0)=0$. Thus every zero of a nonzero $f_j$ is a zero of the real polynomial
\[
P(z)\coloneqq\prod_{f_j\not\equiv0}P_j(z,0)\overline{P_j(\bar z,0)},\qquad \deg P\leq\operatorname{poly}_m(\beta).
\]
Let $\lambda\geq \e^2$ be chosen below, and let $E$ and the segments $I_\alpha$ of the real cells $\theta_\alpha+\mathcal B_\alpha$ be given by Lemma~\ref{lem:real-cellular-cover}. Each nonzero $f_j$ has no zeros on $\theta_\alpha+\mathcal B_\alpha^{[\lambda_\alpha]}$, where $\lambda_\alpha\coloneqq\lambda$, or $\lambda_\alpha\coloneqq\log\lambda$ if $\mathcal B_\alpha$ is a disc; since $\mu\leq\lambda_\alpha$, also $\theta_\alpha+\mathcal B_\alpha^{[\mu]}\subset\mathcal L^{[\nu]}$. Since $\log\lambda\geq2\mu$ for the choice of $\lambda$ below, by \eqref{eq:sigma-extension} the cell $\mathcal B_\alpha^{[\mu]}$ has a $\{\sigma_\alpha\}$-extension, in the sense of \cite{BN21}, contained in $\mathcal B_\alpha^{[\lambda_\alpha]}$, with
\[
\sigma_\alpha\lesssim_\mu\frac1{\log\lambda}.
\]
The one-dimensional monomialization lemma \cite{BN21}*{Lemma 46}, applied to each $f_j(\theta_\alpha+\cdot)$, which is bounded and has no zeros on this extension, with $\mathcal B_\alpha^{[\mu]}$ as inner cell, gives on $\mathcal B_\alpha^{[\mu]}$
\[
f_j(\theta_\alpha+z)=z^{n_{\alpha,j}}u_{\alpha,j}(z),\qquad |n_{\alpha,j}|\leq\operatorname{poly}_{m,\mu,\nu}(\beta),\qquad \diam_{\mathcal B_\alpha^{[\mu]}}\log u_{\alpha,j}\lesssim_\mu\frac{\operatorname{poly}_m(\beta)}{\log\lambda}.
\]
Choose
\[
\log\lambda\coloneqq\frac{C_{m,\mu,\nu}\operatorname{poly}_m(\beta)}{\eta},
\]
with $C_{m,\mu,\nu}$ large. Then all logarithmic diameters are at most $\eta$, and Lemma~\ref{lem:real-cellular-cover} gives
\[
\#E+\#\{\alpha\}\leq\operatorname{poly}_\nu(\deg P,\log\lambda)=\operatorname{poly}_{m,\mu,\nu}(\beta,\eta^{-1}).
\]
Discs are replaced by punctured discs at the end, which does not change the segments. Finally, the $f_j$ are real on $\mathcal L^{[\nu]}\cap\R$, by the identity theorem as noted in Subsection~\ref{subsec:algebraic-cellular-preparation}, so all units are real on the real axis.
\end{proof}

Preparation bounds the oscillation of the units, while the reduction weights of Proposition~\ref{prop:monomial-frame-PU} involve their $\mathfrak D$-derivatives. The two are related through the logarithmic coordinate
\[
z=\e^\xi,\qquad \mathfrak D=\frac{\d}{\d\xi},
\]
in which $A(r,R)$ is the strip $\{\log r<\operatorname{Re}\xi<\log R\}$ modulo $2\pi i$, and $\mathcal L^{[\mu]}$ is $\mathcal L$ widened by $\log\mu$ on each side. Cauchy estimates in the wider strip turn small oscillation into small derivatives, with constants independent of the cell.

\begin{lemma}  \label{lem:cellular-small-unit-Cauchy} Fix $\mu>1$, and let $\mathcal L$ be a punctured disc or an annulus. If $u$ has no zeros on $\mathcal L^{[\mu]}$, admits a single-valued logarithm there, and $\diam_{\mathcal L^{[\mu]}}\log u\leq\eta$, then for every $r\geq1$
\begin{equation}\label{eq:cellular-small-unit-Cauchy}
\sup_{\mathcal L}|\mathfrak D^r\log u|\lesssim_{r,\mu}\eta.
\end{equation}
\end{lemma}

\begin{proof} Fix a branch of $\log u$ and $z_0\in\mathcal L^{[\mu]}$, and put $h\coloneqq\log u-\log u(z_0)$, so that $|h|\leq\eta$ on $\mathcal L^{[\mu]}$. Each $z\in\mathcal L$ is the center of a disc of radius $c_\mu|z|$ contained in $\mathcal L^{[\mu]}$, so Cauchy's estimates give $|z^kh^{(k)}(z)|\lesssim_{k,\mu}\eta$ for $k\geq1$. Since $\mathfrak D^r$ is a fixed linear combination of the operators $z^k\partial_z^k$, $1\leq k\leq r$, this proves \eqref{eq:cellular-small-unit-Cauchy}.
\end{proof}

The next lemma integrates functions of the form $z^\rho h$ with $h$ almost constant, showing that the result has the same form. It is used in the proof of Proposition~\ref{prop:monomial-projective-frame} below, which alternates differentiation and integration. The monomials $z^k=\e^{k\xi}$ are the modes of the Laurent expansion, and $\mathfrak D=\d/\d\xi$ multiplies $z^k$ by $k$. Hence, if $\mathfrak Dg=z^\rho h$ with $h$ almost constant, then $g$ is, up to a constant, close to $z^\rho h/\rho$, since $\mathfrak D(z^\rho/\rho)=z^\rho$. The other modes $z^{\rho+k}$ of $z^\rho h$ are divided by $\rho+k$ instead of $\rho$, which costs at most a factor $|\rho|$. Moreover, $\mathfrak Dg$ has no constant Laurent coefficient, as $g$ is single-valued and a primitive of a nonzero constant in the coordinate $\xi$ is not periodic. Hence the mode with $\rho+k=0$ is absent, and $\rho\neq0$, since for $\rho=0$ this mode would be the main, nonzero coefficient of the almost constant $h$. Each integration consumes a fixed amount of room, and we fix once and for all a number $\kappa>1$ for it.

\begin{lemma} \label{lem:annular-primitive-derivation} Let $\mathcal A=A(r,R)$, with $0\leq r<R$, and let $g$ be single-valued and holomorphic on $\mathcal A^{[\kappa]}$. Suppose
\[
\mathfrak Dg=z^\rho h, \qquad\rho\in\Z,
\]
where $h$ is nonvanishing on $\mathcal A^{[\kappa]}$ and admits a single-valued logarithm there. There are $c_\kappa>0$ and $C_\kappa<\infty$, depending only on $\kappa$, such that, if
\begin{equation}\label{eq:annular-derivation-smallness}
\diam_{\mathcal A^{[\kappa]}}\log h \leq \frac{c_\kappa}{1+|\rho|},
\end{equation}
then
\[
\rho\neq0,\qquad g-c=z^\rho v,
\]
for some $c\in\C$ and some nonvanishing holomorphic function $v$ on $\mathcal A$ admitting a single-valued logarithm there. Moreover
\begin{equation} \label{eq:annular-derivation-unit-bound}
\diam_{\mathcal A}\log v \leq C_\kappa |\rho|  \diam_{\mathcal A^{[\kappa]}}\log h .
\end{equation}
If $g$ and $h$ are real-valued on the positive real axis, then $c\in\R$ and $v$ is real-valued there.
\end{lemma}

\begin{proof} Choose a branch of $\log h$ and set
\[
\omega \coloneqq \diam_{\mathcal A^{[\kappa]}}\log h .
\]
Fix $z_0\in\mathcal A^{[\kappa]}$ and set $h_*=h(z_0)$. Since $h$ is nonvanishing on $\mathcal A^{[\kappa]}$, it follows that $h_*\neq0$. Notice that
\[
\left\|\frac{h}{h_*}-1\right\|_{L^\infty(\mathcal A^{[\kappa]})} \leq \e^\omega-1 \leq 2\omega,
\]
provided $c_\kappa$ is sufficiently small. Write
\[
\frac{h(z)}{h_*}-1 = \sum_{k\in\Z}\varepsilon_k z^k.
\]
The constant Laurent coefficient of $\mathfrak Dg$ is zero by definition of $\mathfrak D$. On the other hand, by the expansion of $h/h_*$, the constant coefficient of $z^\rho h=\mathfrak Dg$ equals $h_*(\delta_{\rho,0}+\varepsilon_{-\rho})$. If $\rho=0$, this is $h_*(1+\varepsilon_0)\neq0$, since a Laurent coefficient is bounded by the supremum of the function on a circle, so that $|\varepsilon_0|\leq2\omega<1$. This contradiction shows that $\rho\neq0$. For $\rho\neq0$ it gives $\varepsilon_{-\rho}=0$, which is exactly what allows us to integrate term by term. Define
\[
v(z)\coloneqq\frac{h_*}{\rho}+h_*\sum_{k\neq-\rho}\frac{\varepsilon_k}{\rho+k}z^k.
\]
The series converges normally on $\mathcal A$, and $\mathfrak D(z^\rho v)=z^\rho h=\mathfrak Dg$. Thus $g-z^\rho v$ is constant, so $g-c=z^\rho v$ for some $c\in\C$. Since $\rho+k$ is a nonzero integer when $k\neq-\rho$,
\[
\left|\frac{\rho}{\rho+k}\right|\leq|\rho|.
\]
Cauchy's estimates on the circles of radii $\kappa R$ and $r/\kappa$ bound $|\varepsilon_kz^k|$ on $\mathcal A$ by $2\omega\kappa^{-|k|}$, and summing this geometric series gives
\[
\left\|\frac{v}{h_*/\rho}-1\right\|_{L^\infty(\mathcal A)} \leq C_\kappa|\rho|\omega.
\]
For a punctured disc, the same estimate follows after observing that the bounded function $h/h_*-1$ has a removable singularity at the puncture. The hypothesis \eqref{eq:annular-derivation-smallness}, with $c_\kappa$ sufficiently small, makes the preceding right-hand side at most $1/4$. Hence $v$ is nonvanishing and admits a single-valued logarithm on $\mathcal A$. Moreover,
\[
\diam_{\mathcal A}\log v \leq C_\kappa|\rho|\omega,
\]
which gives \eqref{eq:annular-derivation-unit-bound}. Finally, if $g$ and $h$ are real-valued on the positive real axis, choose $z_0$ there. Then the construction is invariant under complex conjugation, so $v$ is real-valued on the positive real axis. Since $g-z^\rho v$ is the constant $c$, it follows that $c\in\R$.
\end{proof}

For any $n$-dimensional space $V$ of bounded algebraic functions on a cell, Proposition~\ref{prop:monomial-projective-frame} below constructs, after subdividing the cell, a basis of $V$ cellularly prepared with pairwise distinct exponents. On punctured discs such a basis could be obtained by elimination, but not on annuli, as the following remark explains.

\begin{remark}  \label{rem:index-not-valuation} If $F=z^mu$ has no zeros on an annulus and $u$ admits a single-valued logarithm, then $m$ is the winding number of $F$ around the origin. On a punctured disc the unit of a prepared function extends holomorphically and nonvanishingly across the puncture, so $m$ is the order of $F$ at $0$, and orders satisfy $\ord_0(F+G)\geq\min\{\ord_0F,\ord_0G\}$; this is what allows bases with distinct vanishing orders. Winding numbers have no such property: on $A(2,R)$ the functions $z$ and $1-z$ have winding number $1$, while their sum has winding number $0$. Thus on annuli the exponents cannot be separated by elimination inside a space of functions, and the next proposition differentiates and integrates instead.
\end{remark}

The proof of Proposition~\ref{prop:monomial-projective-frame} is an induction on $\dim V$. Once one function $F_0\in V$ is prepared, the derivatives $(F/F_0)'$, $F\in V$, multiplied by a power of $z$ to keep them bounded, span a space of dimension $\dim V-1$, and Lemma~\ref{lem:annular-primitive-derivation} integrates a prepared basis of that space back to $V$. Since the induction changes the cells, $F_0$ is carried along in an auxiliary family $\mathcal S$ of functions prepared on the final cells.

\begin{proposition}\label{prop:monomial-projective-frame} Fix $n\geq1$ and $\mu,\nu>1$. Let $\mathcal L$ be a punctured disc or an annulus, $V$ be an $n$-dimensional real vector space of bounded algebraic holomorphic functions on $\mathcal L^{[\nu]}$, real on its positive real part, with a basis of complexity at most $\beta$, and let $\mathcal S$ be a finite family of such functions, of complexity at most $\beta$. Then for each $0<\eta<1$ there is a $(\mu,\eta)$-prepared cellular cover of $\mathcal L$ for $\mathcal S$, with exceptional set $E$ and segments $I_\alpha$ of real cells $\theta_\alpha+\mathcal B_\alpha$, such that
\begin{enumerate}[itemsep=.4em]
\item[\emph{(i)}] $\#E+\#\{\alpha\}\leq\operatorname{poly}_{n,\#\mathcal S,\mu,\nu}(\beta,\eta^{-1})$;
\item[\emph{(ii)}] for each $\alpha$, the space $V$ has a basis $G_0,\ldots,G_{n-1}$ which is cellularly prepared on $\theta_\alpha+\mathcal B_\alpha^{[\mu]}$ with accuracy $\eta$ and distinct exponents:
\begin{equation}\label{eq:monomial-frame-form}
G_i(\theta_\alpha+z)=z^{m_i}u_i(z),\qquad m_0<m_1<\cdots<m_{n-1};
\end{equation}
\item[\emph{(iii)}] the exponents of the $G_i$ and of the functions in $\mathcal S$ are bounded by $\operatorname{poly}_{n,\#\mathcal S,\mu,\nu}(\beta)$, independently of $\eta$;
\item[\emph{(iv)}] all units are real on the real axis.
\end{enumerate}
\end{proposition}

\begin{proof} Nonzero prepared functions have no zeros, and products of prepared functions add exponents, while the logarithmic diameters of their units are subadditive. We argue by induction on $n$. If $n=1$, Proposition~\ref{prop:effective-cellular-preparation} applies to a basis element of $V$ together with $\mathcal S$. Let $n\geq2$, and fix a basis $F_0,\ldots,F_{n-1}$ of $V$ of complexity at most $\beta$.

Proposition~\ref{prop:effective-cellular-preparation}, applied to $F_0$ with dilation $\kappa^2\mu$ and accuracy $\frac12$, gives segments of real cells $\theta+\mathcal A$ such that $F_0$ is prepared, hence nowhere zero, on $\theta+\mathcal A^{[\kappa^2\mu]}$. Fix one of them, translate it to $\theta=0$, and write $F_0=z^b u$ on $\mathcal A^{[\kappa^2\mu]}$. The functions
\begin{equation}\label{eq:descent-derivative}
L_j\coloneqq z^{b+1}\left(\frac{F_j}{F_0}\right)',\qquad 1\leq j\leq n-1,
\end{equation}
are holomorphic on $\mathcal A^{[\kappa^2\mu]}$, of complexity $\operatorname{poly}_{n,\#\mathcal S,\mu,\nu}(\beta)$, since quotients and derivatives of algebraic functions of complexity $\beta$ have complexity $\operatorname{poly}(\beta)$ by elimination and implicit differentiation, and linearly independent, since $(F/F_0)'=0$ forces $F\in\R F_0$, all functions being real on the real axis. They are bounded on $\mathcal A^{[\kappa\mu]}$: with $q_j\coloneqq F_j/u$, bounded on $\mathcal A^{[\kappa^2\mu]}$ since $\diam\log u\leq\frac12$, we have $L_j=-bq_j+zq_j'$, and $zq_j'$ is bounded on $\mathcal A^{[\kappa\mu]}$ by Cauchy's estimates on discs of radius $c_\kappa|z|$. The proposition in dimension $n-1$, applied to the span of the $L_j$ on $\mathcal A$, with $\kappa\mu$ in place of both $\mu$ and $\nu$, auxiliary family $\mathcal S\cup\{F_0,z^{b+1}\}$ and accuracy $\widetilde\eta$ chosen below, gives a cover of $\mathcal A$. For each of its real cells $\theta'+\mathcal B$ there are $(a_{ij})\in\operatorname{GL}_{n-1}(\R)$ and distinct integers $k_1,\ldots,k_{n-1}$ with
\begin{equation}\label{eq:monomial-frame-induction-output}
\left(\sum_{j}a_{ij}L_j\right)(\theta'+z)=z^{k_i}U_i(z),\qquad (\theta'+z)^{b+1}=z^pw(z)\qquad\text{on }\mathcal B^{[\kappa\mu]},
\end{equation}
and $F_0(\theta'+z)=z^{m_0}u_0(z)$ there, where all exponents are bounded by some $A_1=\operatorname{poly}_{n,\#\mathcal S,\mu,\nu}(\beta)$ independent of $\widetilde\eta$, and these functions, as well as those of $\mathcal S$, are prepared on $\theta'+\mathcal B^{[\kappa\mu]}$ with accuracy $\widetilde\eta$.

Put $\widetilde F_i\coloneqq\sum_ja_{ij}F_j$. By \eqref{eq:descent-derivative} and \eqref{eq:monomial-frame-induction-output}, the functions $g_i(z)\coloneqq(\widetilde F_i/F_0)(\theta'+z)$ satisfy
\[
\mathfrak Dg_i=z^{\rho_i}\frac{U_i}{w},\qquad \rho_i\coloneqq k_i-p+1,\qquad |\rho_i|\leq 2A_1+1,
\]
where $\diam\log(U_i/w)\leq2\widetilde\eta$ and the $\rho_i$ are distinct. Choose
\begin{equation}\label{eq:monomial-frame-accuracy-choice}
\widetilde\eta\coloneqq\min\left\{\frac{c_\kappa}{4(A_1+1)},\ \frac{\eta}{4\max\{1,C_\kappa(2A_1+1)\}}\right\},
\end{equation}
with $c_\kappa$ and $C_\kappa$ as in Lemma~\ref{lem:annular-primitive-derivation}. That lemma, applied on $\mathcal B^{[\mu]}$, gives $c_i\in\R$ and units $v_i$, real on the real axis, such that
\[
g_i-c_i=z^{\rho_i}v_i,\qquad \rho_i\neq0,\qquad \diam_{\mathcal B^{[\mu]}}\log v_i\leq 2C_\kappa|\rho_i|\widetilde\eta\leq\frac\eta2.
\]
Hence $G_0\coloneqq F_0$ and $G_i\coloneqq\widetilde F_i-c_iF_0$, a basis of $V$ since $(a_{ij})$ is invertible, satisfy
\[
G_0(\theta'+z)=z^{m_0}u_0(z),\qquad G_i(\theta'+z)=z^{m_0+\rho_i}u_0(z)v_i(z)\qquad\text{on }\mathcal B^{[\mu]},
\]
with distinct exponents, since the $\rho_i$ are distinct and nonzero, bounded by $3A_1+1$, and with logarithmic diameters at most $\widetilde\eta+\frac\eta2\leq\eta$.

Undoing the translation, the segments of the real cells $\theta+\theta'+\mathcal B$, for which $\theta+\theta'+\mathcal B^{[\mu]}\subset\theta+\mathcal A^{[\kappa\mu]}\subset\mathcal L^{[\nu]}$, cover the segment of $\theta+\mathcal A$ up to finitely many points. Together with the exceptional points of the first cover, they form the required cover, as $\widetilde\eta\leq\eta$, and a permutation of the $G_i$ gives \eqref{eq:monomial-frame-form}. Each level of the induction consumes a factor $\kappa^2$ of dilation, and $\widetilde\eta^{-1}\leq\operatorname{poly}_{n,\#\mathcal S,\mu,\nu}(\beta,\eta^{-1})$ by \eqref{eq:monomial-frame-accuracy-choice}, which gives the stated bounds.
\end{proof}

\subsection{Power-unit weights} \label{subsec:monomial-frames-to-PU}
Let $J$ be an interval and $Y_0=1$, $Y_1=t$, $Y_2,\ldots,Y_e$ real functions of class $\mathcal C^e$ on $J$; we call such a family a \emph{system} on $J$. Their \emph{reduction} alternates differentiation and division: $D_0h\coloneqq h$ and
\begin{equation}\label{eq:reduction}
w_q\coloneqq(D_{q-1}Y_q)',\qquad D_qh\coloneqq\frac{(D_{q-1}h)'}{w_q},\qquad 1\leq q\leq e,
\end{equation}
as long as $w_q$ has no zeros on $J$. Then
\[
w_1=1,\qquad D_qY_q=1,\qquad D_qY_p=0\quad\text{for }p<q,
\]
and $w_2,\ldots,w_e$ are the \emph{reduction weights}. For a canonical curve with continuous weights, \eqref{eq:reduction} is defined on the span of its coordinates and returns its weights and its generalized derivatives, by \eqref{eq:generalized-derivatives-coordinates} and \eqref{eq:generalized-derivative-recursion}. Replacing $Y_q$ by $-Y_q$ changes the signs of $w_q$ and $w_{q+1}$ only, so successive sign changes of $Y_2,\ldots,Y_e$ make all reduction weights positive.

The next proposition computes the reduction weights of a basis as in Proposition~\ref{prop:monomial-projective-frame}: in the parameter $t=G_1/G_0$, they are powers of $t$ times nearly constant units.

\begin{proposition}  \label{prop:monomial-frame-PU} There is $c_e>0$, depending only on $e$, such that the following holds. Let $\mathcal B$ be a punctured disc or an annulus, $L\geq2$, and
\begin{equation}\label{eq:monomial-projective-frame}
G_i=z^{m_i}u_i\quad\text{on }\mathcal B,\qquad m_0<m_1<\cdots<m_e,\qquad |m_i|\leq L,
\end{equation}
where the units $u_i$ have no zeros, admit single-valued logarithms, are real on $\mathcal B\cap(0,\infty)$, and satisfy
\begin{equation}\label{eq:monomial-unit-smallness}
\varepsilon\coloneqq\max_{0\leq i\leq e}\diam_{\mathcal B}\log u_i+\max_{\substack{0\leq i\leq e\\1\leq r\leq e}}\left\|\mathfrak D^r\log u_i\right\|_{L^\infty(\mathcal B)}\leq\frac{c_e}{L}.
\end{equation}
Then, after changing the signs of some $G_i$,
\begin{enumerate}[itemsep=.4em]
\item[\emph{(i)}] $t\coloneqq G_1/G_0$ is positive and increasing on $\mathcal B\cap(0,\infty)$, with range a bounded open interval $J\subset(0,\infty)$;
\item[\emph{(ii)}] the reduction of $(1,t,G_2/G_0,\ldots,G_e/G_0)$ is defined on $J$, with weights
\begin{equation}\label{eq:algebraic-PU-weights}
w_q=t^{a_q-1}U_q,\qquad a_q=\frac{m_q-m_{q-1}}{m_1-m_0}\in\left[\frac1{2L},2L\right],\qquad \sup_JU_q\leq2\inf_JU_q,
\end{equation}
where the units $U_q$ are positive.
\end{enumerate}
\end{proposition}

\begin{proof} Write $Y_i\coloneqq G_i/G_0=z^{m_i-m_0}v_i$, with $v_i\coloneqq u_i/u_0$, so that the $\log v_i$ satisfy \eqref{eq:monomial-unit-smallness} with $2\varepsilon$ in place of $\varepsilon$. The proof rests on the identity
\[
\mathfrak D\left(z^kV\right)=kz^kV\left(1+\frac{\mathfrak D\log V}{k}\right),\qquad k\in\Z\setminus\{0\},
\]
by which a monomial times a unit differentiates into a monomial times a unit. Changing the sign of $G_1$, we may assume $v_1>0$ on $\mathcal B\cap(0,\infty)$. There $t=Y_1$ and $\mathfrak D\log t=m_1-m_0+\mathfrak D\log v_1\geq1-2\varepsilon>0$, so $t$ is positive and increasing, $\d/\d t=\mathfrak D/\mathfrak Dt$, and \eqref{eq:reduction} becomes
\begin{equation}\label{eq:reduction-monomial-variable}
D_qh=\frac{\mathfrak D(D_{q-1}h)}{\mathfrak D(D_{q-1}Y_q)},\qquad w_q=\frac{\mathfrak D(D_{q-1}Y_q)}{\mathfrak Dt}.
\end{equation}

We claim that, for $0\leq q\leq p\leq e$,
\begin{equation}\label{eq:reduction-monomial}
D_qY_p=z^{m_p-m_q}V_{p,q},\qquad \diam_{\mathcal B}\log V_{p,q}+\max_{1\leq r\leq e-q}\left\|\mathfrak D^r\log V_{p,q}\right\|_{L^\infty(\mathcal B)}\lesssim_e\varepsilon,
\end{equation}
with units $V_{p,q}$ real on $\mathcal B\cap(0,\infty)$. For $q=0$ take $V_{p,0}\coloneqq v_p$. If \eqref{eq:reduction-monomial} holds for $q-1$, put $k_p\coloneqq m_p-m_{q-1}$, an integer which is at least $1$ for $p\geq q$. The identity and \eqref{eq:reduction-monomial-variable} give $D_qY_p=z^{m_p-m_q}V_{p,q}$, with
\[
V_{p,q}=\frac{k_pV_{p,q-1}}{k_qV_{q,q-1}}\cdot\frac{1+k_p^{-1}\mathfrak D\log V_{p,q-1}}{1+k_q^{-1}\mathfrak D\log V_{q,q-1}}.
\]
If $c_e$ is small, the last factor lies in the disc of radius $\frac12$ about $1$. With the principal branch of its logarithm, $\mathfrak D^r\log V_{p,q}$, $r\leq e-q$, is controlled by the $\mathfrak D^s\log V_{\cdot,q-1}$ with $s\leq r+1$, which proves the claim.

In the same way, with $k_1=m_1-m_0$ and $k_q=m_q-m_{q-1}$, \eqref{eq:reduction-monomial-variable} gives
\[
w_q=z^{k_q-k_1}\Theta_q,\qquad \Theta_q\coloneqq\frac{k_qV_{q,q-1}}{k_1v_1}\cdot\frac{1+k_q^{-1}\mathfrak D\log V_{q,q-1}}{1+k_1^{-1}\mathfrak D\log v_1},
\]
with $\Theta_q$ real and nonvanishing on $\mathcal B\cap(0,\infty)$, and $\diam\log|\Theta_q|\lesssim_e\varepsilon$ there. Since $z=(t/v_1)^{1/k_1}$ on $\mathcal B\cap(0,\infty)$, and $1\leq k_1,k_q\leq2L$,
\[
w_q=t^{a_q-1}v_1^{1-a_q}\Theta_q,\qquad a_q=\frac{k_q}{k_1}=\frac{m_q-m_{q-1}}{m_1-m_0}\in\left[\frac1{2L},2L\right].
\]
After sign changes of $G_2,\ldots,G_e$, every $w_q$ is positive, and $U_q\coloneqq v_1^{1-a_q}|\Theta_q|$ satisfies, by \eqref{eq:monomial-unit-smallness},
\[
\log\frac{\sup U_q}{\inf U_q}\leq2|1-a_q|\varepsilon+C_e\varepsilon\leq(4L+C_e)\varepsilon\leq\log2.
\]
Finally, $J$ is bounded since $v_1$ is.
\end{proof}

The next lemma shows that an arc whose reduction weights have the form of Proposition~\ref{prop:monomial-frame-PU} on a bounded interval $J\subset(0,\infty)$ lies on a curve of $\mathrm{PPU}(e,\Lambda,R)$, as in Definition~\ref{def:power-unit-class}.

\begin{lemma}\label{lem:fill-annular-hole}
Let $e\geq2$, let $J\subset(0,\infty)$ be a bounded open interval, and let $Y_0,\ldots,Y_e$ be a real-analytic system on $J$, with reduction weights
\[
w_q=t^{a_q-1}u_q,\qquad \Lambda^{-1}\leq a_q\leq\Lambda,\qquad u_q>0,\qquad \sup_Ju_q\leq R\inf_Ju_q,\qquad 2\leq q\leq e.
\]
Then the arc
\[
\left\{[1:t:Y_2(t):\cdots:Y_e(t)]: t\in J\right\}
\]
is contained in a curve of $\mathrm{PPU}(e,\Lambda,R)$.
\end{lemma}

\begin{proof} Let $J=(t_0,t_1)$, and extend each $u_q$ to $(0,t_0]$ by the constant $\inf_Ju_q$. This keeps $\sup u_q\leq R\inf u_q$, and $t^{a_q-1}u_q$ is integrable on $(0,t_1)$ since $a_q>0$. Definition~\ref{def:power-unit-class} asks no continuity of the units, so the canonical curve $\mathbf X$ of the weights
\[
s\longmapsto c_q(t_1s)^{a_q-1}u_q(t_1s),\qquad 0<s<1,
\]
with $c_q$ given by \eqref{e:normweights}, belongs to $\mathrm{PU}(e,\Lambda,R)$. Put $\widetilde X_p(t)\coloneqq X_p(t/t_1)$. On $J$ the weights of $\mathbf X$, read in the variable $t$, are positive multiples of the $w_q$. By \eqref{eq:generalized-derivatives-coordinates} and \eqref{eq:generalized-derivative-recursion}, each $D_q$ of \eqref{eq:reduction} is then a positive multiple of the generalized derivative $\mathcal D_q$ of $\mathbf X$, so that, for $0\leq p\leq e$, both families
\[
Y_0,\ldots,Y_p\qquad\text{and}\qquad\widetilde X_0,\ldots,\widetilde X_p
\]
lie in the kernel of $h\mapsto(D_ph)'$ on $\mathcal C^\infty(J)$. Both are linearly independent, by the triangular relations $D_rY_r=1$, $D_rY_q=0$ for $q<r$, and their counterparts for $\mathbf X$. The kernel has dimension $p+1$, since $(D_ph)'=0$ determines $h$ by $p+1$ successive integrations. Therefore $Y_p=\sum_{r\leq p}M_{pr}\widetilde X_r$ on $J$ for a constant lower-triangular $M\in\operatorname{GL}_{e+1}(\R)$, and the arc is contained in $M[\mathbf X(I)]\in\mathrm{PPU}(e,\Lambda,R)$.
\end{proof}

\subsection{Proofs of the decomposition and of the maximal estimate} We first assemble the results of this section into a proof of Theorem~\ref{thm:algebraic-PPU-decomposition}. Combined with the analytic results of Sections~\ref{sec:projective-switch-classes} and \ref{sec:power-unit}, it then gives the maximal estimate for algebraic curves, Corollary~\ref{cor:algebraic-maximal}, and Theorem~\ref{thm:intro-main}.

\begin{proof}[Proof of Theorem~\ref{thm:algebraic-PPU-decomposition}] The irreducible components of $Z^{\C}$ are points and curves, at most $D$ in number, and the degrees of the curves sum to at most $D$. The real points among the former are placed in $E$. If one of the irreducible curves, $Y$, is not fixed by complex conjugation, its real points lie in $Y\cap\overline Y$, which has at most $D^{C_d}$ points by B\'ezout \cite{BN21}*{\S2.2.3}, and they are placed in $E$. Otherwise the smallest complex projective subspace containing $Y$ is invariant under conjugation, hence the complexification of a real projective subspace $\Pi\simeq\mathbb R\mathrm P^e$, $e\leq d$, and $Y$ is nondegenerate in $\Pi^{\C}\simeq\mathbb C\mathrm P^e$, of degree $\delta\leq D$. If $e=1$, then $Y(\R)=\Pi$, and removing two points, placed in $E$, leaves two open arcs, whose closures belong to $\mathrm{PPU}(1,\Lambda_D,R_D)$ by Remark~\ref{rem:PU-dimension-one}.

Let $e\geq2$, and fix $\mu_\ast>1$ once and for all. Up to $D^{C_d}$ points, placed in $E$, Proposition~\ref{prop:algebraic-monotone-arcs} covers $Y(\R)$ by $D^{C_d}$ arcs $\Gamma$, with injective affine coordinates $x$, and Lemma~\ref{lem:effective-scalar-lifting}, with $\mu=\mu_\ast$, covers each $\Gamma$, up to $D^{C_d}$ points placed in $E$, by the positive real images of $D^{C_d}$ maps $\gamma_\alpha:\mathcal B_\alpha^{[\mu_\ast]}\to Y$. In the affine chart of $\Gamma$, with its homogeneous coordinate placed first, write $\gamma_\alpha=[1:F_{\alpha,1}:\cdots:F_{\alpha,e}]$, and let $V_\alpha$ be the real span of $1,F_{\alpha,1},\ldots,F_{\alpha,e}$. As $Y$ is nondegenerate, no nonzero linear form vanishes on the infinite set $\gamma_\alpha(\mathcal B_\alpha\cap(0,\infty))$, so $\dim V_\alpha=e+1$. Proposition~\ref{prop:monomial-projective-frame}, applied to $V_\alpha$ on $\mathcal B_\alpha$, with $\mu=\nu=\mu_\ast$, $\mathcal S=\varnothing$ and accuracy $\eta\coloneqq D^{-A_d}$, where $A_d$ is fixed below, gives $D^{C_d}$ points of $\mathcal B_\alpha\cap(0,\infty)$, whose images under $\gamma_\alpha$ are placed in $E$, and covers the rest of $\mathcal B_\alpha\cap(0,\infty)$ by $D^{C_d}$ segments $I_{\alpha,j}$ of real cells $\theta_{\alpha,j}+\mathcal B_{\alpha,j}$. On each of them, after a constant real projective change of frame, the homogeneous lift of $\gamma_\alpha(\theta_{\alpha,j}+\cdot)$ is
\begin{equation}\label{eq:algebraic-final-monomial-frame}
G_i=z^{m_i}u_i\quad\text{on }\mathcal B_{\alpha,j}^{[\mu_\ast]},\qquad m_0<\cdots<m_e,\qquad |m_i|\leq L,\qquad \diam\log u_i\leq D^{-A_d},
\end{equation}
where $L\geq2$ is the exponent bound of Proposition~\ref{prop:monomial-projective-frame}, maximized over $2\leq e\leq d$ and $\delta\leq D$. Thus $L\leq D^{C_d}$ does not depend on $A_d$, and \eqref{eq:monomial-unit-smallness} holds on $\mathcal B_{\alpha,j}$ by Lemma~\ref{lem:cellular-small-unit-Cauchy}, if $A_d$ is large. By Proposition~\ref{prop:monomial-frame-PU}, after a diagonal change of frame, $\gamma_\alpha(I_{\alpha,j})$ is the arc $\{[1:t:G_2/G_0:\cdots:G_e/G_0]: t\in J\}$, with reduction weights \eqref{eq:algebraic-PU-weights}, and Lemma~\ref{lem:fill-annular-hole} places it on a curve of $\mathrm{PPU}(e,\Lambda_D,R_D)$, where
\[
\Lambda_D\coloneqq2L\leq D^{C_d},\qquad R_D\coloneqq2.
\]

Each subdivision is polynomial in $D$, including the at most $d$ levels of the frame construction, whose counts are polynomial in $\eta^{-1}=D^{A_d}$. Hence $\#E+M\leq D^{C_d}$. The arcs $\mathcal C_\nu$ are those of the case $e=1$ and the images $\gamma_\alpha(I_{\alpha,j})$. Such an image is $\{q\in\Gamma: x(q)\in J'\}$ for an open interval $J'$, since $x$ is injective on $\Gamma$ and $x\circ\gamma_\alpha$ is continuous and injective on $I_{\alpha,j}$, hence a real algebraic arc. Its span is $\Pi$, as a real hyperplane of $\Pi$ containing it would contain infinitely many points of $Y$, hence $Y$. With $\phi_\nu$ the identification $\Pi\simeq\mathbb R\mathrm P^e$, this proves \eqref{eq:algebraic-PPU-decomposition} and \eqref{eq:algebraic-PPU-extension}.
\end{proof}

We now turn to the maximal estimate.

\begin{corollary}\label{cor:algebraic-maximal} Let $d\geq1$ and $\alpha>0$. There are $C_{d,\alpha,\varphi}<\infty$ and $C_d<\infty$ such that, for every $D\geq2$,
\begin{equation} \label{eq:algebraic-maximal}
\sup_{\substack{Z\subset\Pd\\\dim Z\leq1,\,\deg Z\leq D}}\Q_d(Z,N)\leq C_{d,\alpha,\varphi} D^{C_d} N^\alpha,\qquad N\geq1.
\end{equation}
More precisely, with $b_d$ as in Theorem~\ref{thm:projective-SECT-maximal},
\begin{equation}\label{eq:algebraic-subpower}
\sup_{\substack{Z\subset\Pd\\\dim Z\leq1,\,\deg Z\leq D}}\Q_d(Z,N)\leq C_{d,\varphi}D^{C_d}\exp\left(b_d\sqrt{\log(2+N)}\right),\qquad N\geq1.
\end{equation}
\end{corollary}

\begin{proof} Since $\exp(b_d\sqrt{\log(2+N)})\lesssim_{d,\alpha}N^\alpha$, it suffices to prove \eqref{eq:algebraic-subpower}. Apply Theorem~\ref{thm:algebraic-PPU-decomposition}. A subset of $E$ has at most $\#E$ points, so \eqref{eq:Q-trivial} gives $\Q_d(E,N)\leq C_\varphi\#E\leq C_\varphi D^{C_d}$. Fix $\nu$. By \eqref{eq:algebraic-PPU-extension} and \eqref{eq:PU-SECT-polynomial} if $e_\nu\geq2$, or by Remark~\ref{rem:PU-dimension-one} if $e_\nu=1$,
\[
\phi_\nu(\mathcal C_\nu)\subset\widetilde{\mathcal C}_\nu\in\mathrm{PSECT}(e_\nu,K),\qquad K\coloneqq\begin{cases}K_{e_\nu}(\Lambda_D,R_D)\leq D^{C_d},&e_\nu\geq2,\\1,&e_\nu=1.\end{cases}
\]
Since $\phi_\nu^{-1}$ is induced by an injective linear map $\R^{e_\nu+1}\to\R^{d+1}$ and $\Q_d(\cdot,N)$ is monotone under inclusion, Lemma~\ref{lem:Q-projective-invariance}, $d-e_\nu$ applications of Lemma~\ref{lem:hyperplane-reduction} and \eqref{eq:projective-SECT-subpower} give
\[
\begin{split}
\Q_d(\mathcal C_\nu,N)\leq\Q_{e_\nu}(\widetilde{\mathcal C}_\nu,N)&\leq C_{d,\varphi}K^{p_{e_\nu}}\exp\left(b_{e_\nu}\sqrt{\log(2+N)}\right)
\\
&\leq C_{d,\varphi}D^{C_d}\exp\left(b_d\sqrt{\log(2+N)}\right).
\end{split}
\]
Summing over the at most $D^{C_d}$ arcs with \eqref{eq:Q-trivial} completes the proof.
\end{proof}

\begin{remark} \label{rem:algebraic-maximal-polynomial} The degree enters \eqref{eq:algebraic-subpower} at three places, each polynomially. Theorem~\ref{thm:algebraic-PPU-decomposition} produces $D^{C_d}$ arcs, with $\Lambda_D\leq D^{C_d}$ and, in its proof, $R_D=2$; \eqref{eq:PU-SECT-polynomial} converts these into a switch constant $K\leq D^{C_d}$; and \eqref{eq:projective-SECT-subpower} costs $K^{p_d}$, a factor $K^k$ arising at the step of the recursion \eqref{eq:projective-SECT-recursion} in dimension $k$, for $2\leq k\leq d$. That this cost is incurred once per level, rather than once per scale, is due to the different roles of the two terms of \eqref{eq:claimed-ao}, recorded in the remark following it. In particular, the exponent of $D$ in \eqref{eq:algebraic-maximal} is independent of $\alpha$.
\end{remark}

The proof of Theorem~\ref{thm:intro-main} uses the following bound for the real points of an algebraic set.

\begin{lemma}\label{lem:real-points-degree}
Let $Y\subset\mathbb C\mathrm P^e$ be a conjugation-invariant algebraic set with $\dim Y\leq1$ and $\deg Y\leq D$. Then the Zariski closure of $Y(\R)$ has dimension at most one and degree at most $D+D^3$.
\end{lemma}

\begin{proof} Let $W$ be an irreducible component of $Y$ such that $W(\R)$ is not Zariski dense in $W$. Then $W(\R)$ is finite, and empty if $\dim W=0$. Let $\dim W=1$. If $W\neq\overline W$, then $W(\R)\subset W\cap\overline W$, which has at most $(\deg W)^2$ points by B\'ezout's theorem \cite{BN21}*{\S2.2.3}. If $W=\overline W$, the real points of $W$ are singular, since a smooth real point of a conjugation-invariant curve has a one-dimensional real neighborhood in $W(\R)$. Moreover, $W$ is not a line, since a conjugation-invariant line has infinitely many real points, and Lemma~\ref{lem:plane-model}, applied to $W$ in its span, which is defined over $\R$, gives $\#W(\R)\leq(\deg W)^3$. Hence the Zariski closure of $Y(\R)$ consists of the components $W$ in which $W(\R)$ is Zariski dense, of total degree at most $D$, and of at most $\sum_W(\deg W)^3\leq D^3$ further points.
\end{proof}

\begin{proof}[Proof of Theorem~\ref{thm:intro-main}] We may assume $D\geq2$, replacing $D$ by $\max\{D,2\}$ and absorbing the resulting dimensional constant. Let $\overline Z\subset\mathbb C\mathrm P^n$ be the projective closure of the Zariski closure of $Z$ in $\mathbb C^n$, so that $\deg\overline Z=\deg Z\leq D$, and let $Y'\subset\mathbb C\mathrm P^{n-1}$ be the closure of the image of $\overline Z\setminus\{p_0\}$ under the projection $[x_0:x]\mapsto[x]$ from $p_0=[1:0:\cdots:0]$. Then $Y'$ is defined over $\R$, has dimension at most $1$ and degree at most $D$, and $[v]\in Y'(\R)$ for every $v\in Z\setminus\{0\}$. By Lemma~\ref{lem:real-points-degree}, the Zariski closure of $Y'(\R)$ has dimension at most one and degree at most $D+D^3\leq D^{C_n}$.

Let $V\subset Z$ with $\#V\leq N$, and put $[V]\coloneqq\{[v]: v\in V\setminus\{0\}\}\subset Y'(\R)$. Since $\A_{0,r}f=\varphi(0)f$, we have $\A_Vf\leq\max\{\A_{[V]}f,|\varphi(0)||f|\}$, hence
\[
\|\A_V\|_{L^2(\R^n)\to L^2(\R^n)}\leq\Q_{n-1}\left(Y'(\R),N\right)+|\varphi(0)|.
\]
Corollary~\ref{cor:algebraic-maximal}, with $d=n-1$ and $D^{C_n}$ in place of $D$, bounds the first term. Thus \eqref{eq:algebraic-maximal} with $\alpha=\eps$ gives \eqref{eq:intro-main}. As $\log(2+N)\leq\log(\e+N)$, \eqref{eq:algebraic-subpower} gives \eqref{eq:intro-subpower} with
\[
c_n=b_{n-1}=1+\sum_{j=2}^{n-1}\sqrt{\log j},
\]
after enlarging $C_n$ and $C_{n,\eps}$.
\end{proof}

\section{Doubling canonical charts and Frenet curvatures} \label{sec:doubling-curvatures} This section illustrates the switch condition beyond algebraic curves. Remark~\ref{rem:planar-doubling-switch} gives the basic calculation for $(1,t,F(t))$: doubling of $F''(t)\,\d t$ controls the tangent error relative to the secant error. We extend that comparison to all the canonical weights, using the positive integral formulas for contact residuals and interpolation determinants. An adapted Frenet frame then gives a geometric sufficient condition for these weight measures to be doubling. The argument is independent of the proof of Theorem~\ref{thm:intro-main}. Recall that a canonical curve $\mathbf X=(X_0,\ldots,X_d)$ is generated by a normalized weight vector $\mathbf w=(w_0,\ldots,w_d)$ as in \eqref{e:normweights}, with $w_0=w_1=1$ and
\[
X_j=\mathcal I_{w_1}\cdots\mathcal I_{w_j}\mathbf 1_I,\qquad 2\leq j\leq d.
\]
For such a chart we define the \emph{transition measures} $\mu_q$, $2\leq q\leq d$, by
\[
\d\mu_q(t)\coloneqq w_q(t)\,\d t.
\]
We first impose the doubling condition at the level of a canonical chart. This condition is not itself intrinsic to the projective arc; the geometric class considered below consists of those projective arcs which admit such a chart after a projective change of coordinates. In particular, membership in the resulting projective class depends only on the unparametrized projective arc.

\begin{definition} \label{def:doubling-measure}
Let $K\geq2$. A finite measure $\nu$ on a compact interval $J_0$ is \emph{$K$-doubling} if
\[
\nu(J)\leq K\min\left\{\nu(J^-),\nu(J^+)\right\}
\]
for every interval $J\subset J_0$, where $J^-$ and $J^+$ are the two halves of $J$.
\end{definition}

\begin{lemma}  \label{lem:iterated-doubling}
Let $\nu$ be $K$-doubling on $J_0$. If $J'\subset J\subset J_0$ and $|J'|\geq2^{1-m}|J|$ for an integer $m\geq1$, then
\[
\nu(J)\leq K^m\nu(J').
\]
Moreover, if $\nu\neq0$, then $\nu$ charges every nondegenerate subinterval of $J_0$.
\end{lemma}

\begin{proof} Bisecting $J$ successively $m$ times, one of the resulting intervals of length $2^{-m}|J|$ is contained in $J'$, and iteration of the doubling inequality gives the first assertion. For the second, given a nondegenerate interval $J'\subset J_0$, choose $m$ with $|J'|\geq2^{1-m}|J_0|$; then $0<\nu(J_0)\leq K^m\nu(J')$.
\end{proof}

\begin{definition} \label{def:doubling-canonical}
Let $d\geq2$ and $K\geq 2$. We write $\mathbf X\in\mathrm{DSECT}(d,K)$ if $\mathbf X$ is a canonical curve with normalized weights as in \eqref{e:normweights} and each of its transition measures $\mu_q$, $2\leq q\leq d$, is $K$-doubling. We write $\mathcal C\in\mathrm{PDSECT}(d,K)$ if
\[
B\mathcal C=[\mathbf X(I)]
\]
for some $B\in\operatorname{PGL}_{d+1}(\R)$ and some $\mathbf X\in\mathrm{DSECT}(d,K)$. For $d=1$ there are no transition measures, and we set $\mathrm{DSECT}(1,K)\coloneqq\mathrm{SECT}(1,K)$, the single curve $t\mapsto(1,t)$, $t\in I$, and $\mathrm{PDSECT}(1,K)\coloneqq\mathrm{PSECT}(1,K)$, its projective images.
\end{definition}

The next theorem places these classes inside $\mathrm{SECT}$ and $\mathrm{PSECT}$, respectively.

\begin{theorem}  \label{thm:doubling-canonical-switch}
Let $d\geq2$ and $\mathbf X\in\mathrm{DSECT}(d,K)$. Then
\[
\max_{2\leq j\leq d}\mathsf S_j\left(\mathbf X^{[j]}\right)\leq C_dK^{C_d}.
\]
Consequently,
\[
\mathrm{DSECT}(d,K)\subset\mathrm{SECT}\left(d,C_dK^{C_d}\right),
\qquad
\mathrm{PDSECT}(d,K)\subset\mathrm{PSECT}\left(d,C_dK^{C_d}\right).
\]
\end{theorem}

\begin{proof} Fix $2\leq j\leq d$ and $t_0<\cdots<t_j$. Put
\[
\Pi(\mathbf t)\coloneqq\prod_{m=1}^j(t_m-t_{m-1})\prod_{q=2}^j\prod_{m=1}^{j-q+1}\mu_q\left(\left[t_{m-1},t_{m+q-1}\right]\right).
\]
For an increasing tuple of $e+1$ nodes, $1\leq e\leq j$, $\Pi$ is defined by the same formula with $j$ replaced by $e$; in particular $\Pi(\partial_k\mathbf t)$ is the quantity of order $j-1$. The staircase formula \eqref{eq:staircase-vandermonde} and doubling give, at every order,
\begin{equation} \label{eq:doubling-staircase-short}
c_jK^{-C_j}\Pi(\mathbf t)\leq\Delta_j(\mathbf t)\leq\Pi(\mathbf t).
\end{equation}
Indeed, the upper bound follows by enlarging the staircase domain to the corresponding product of intervals. For the lower bound, put $G_m^{(0)}\coloneqq[t_{m-1},t_m]$. Recursively, let $G_m^{(q)}$ be the interval between the right endpoint of the middle half of $G_m^{(q-1)}$ and the left endpoint of the middle half of $G_{m+1}^{(q-1)}$; by induction, $G_m^{(q-1)}$ lies to the left of $G_{m+1}^{(q-1)}$, so $G_m^{(q)}$ is a nondegenerate interval. An induction gives
\[
|G_m^{(q)}|\geq4^{-q}(t_{m+q}-t_{m-1}).
\]
If $H_m^{(q)}$ denotes the middle half of $G_m^{(q-1)}$, then the product of the intervals $H_m^{(q)}$ is contained in the staircase domain and
\[
|H_m^{(q)}|\geq2^{1-2q}(t_{m+q-1}-t_{m-1}).
\]
Lemma~\ref{lem:iterated-doubling} therefore gives
\[
\mu_q(H_m^{(q)})\geq K^{-2q}\mu_q([t_{m-1},t_{m+q-1}]).
\]
Multiplying over the staircase variables yields \eqref{eq:doubling-staircase-short}.
Likewise, \eqref{eq:canonical-contact-integral} gives, for $J=[\sigma\wedge t,\sigma\vee t]$,
\begin{equation} \label{eq:doubling-residual-short}
|\mathcal R_j(t;\sigma)|\leq|t-\sigma|\prod_{q=2}^j\mu_q(J).
\end{equation}
By \eqref{eq:barycentric-raw-determinant} and \eqref{eq:doubling-staircase-short},
\[
\beta_k^{(j)}(\mathbf t)\leq C_jK^{C_j}\frac{\Pi(\partial_k\mathbf t)}{\Pi(\mathbf t)}.
\]
Hence, by Lemma~\ref{lem:switch-node-reduction}, it suffices to bound
\[
|\mathcal R_j(t_k;t_i)|\frac{\Pi(\partial_k\mathbf t)}{\Pi(\mathbf t)}
\]
for one of the neighboring nodes $|i-k|=1$. Write
\[
a_m\coloneqq t_m-t_{m-1},\qquad x_i^{(q)}\coloneqq\mu_q([t_{i-1},t_i]),\qquad S_\ell^{(q)}\coloneqq\sum_{i=\ell+1}^{\ell+q}x_i^{(q)}.
\]
For $k=0$, \eqref{eq:doubling-residual-short} gives directly
\[
|\mathcal R_j(t_0;t_1)|\frac{\Pi(\partial_0\mathbf t)}{\Pi(\mathbf t)}
\leq\frac{a_1\prod_{q=2}^jx_1^{(q)}}{a_1S_0^{(j)}\prod_{q=2}^{j-1}S_0^{(q)}}\leq1,
\]
and $k=j$ follows by reflection: this bound, like the one below for $1\leq k\leq j-1$, involves only $\Pi$, the right-hand side of \eqref{eq:doubling-residual-short} and the doubling of the $\mu_q$, which keep their form under $t\mapsto1-t$, with the order of the nodes reversed and each $\mu_q$ replaced by its reflection, although the chart itself is not reflection invariant.
Let $1\leq k\leq j-1$, and suppose $a_k\leq a_{k+1}$; the other case follows by the same reflection. Put
\[
r_q\coloneqq\max(0,k-q),\qquad u_q\coloneqq\min(k,j-q).
\]
Deleting $t_k$ merges the two adjacent gaps, and cancellation in the quotient of the two products gives
\[
\frac{\Pi(\partial_k\mathbf t)}{\Pi(\mathbf t)}
=\frac{a_k+a_{k+1}}{a_ka_{k+1}}\frac{1}{S_0^{(j)}}
\prod_{q=2}^{j-1}\left[\frac{1}{S_{u_q}^{(q)}}\prod_{\ell=r_q}^{u_q-1}\frac{S_\ell^{(q)}+x_{\ell+q+1}^{(q)}}{S_\ell^{(q)}}\right].
\]
For $2\leq q\leq j-1$, $1\leq k\leq j-1$ and positive $x_1,\ldots,x_j$, if $S_\ell=x_{\ell+1}+\cdots+x_{\ell+q}$, $r=\max(0,k-q)$ and $u=\min(k,j-q)$, then
\begin{equation} \label{eq:doubling-sliding-short}
\frac{x_k}{S_u}\prod_{\ell=r}^{u-1}\frac{S_\ell+x_{\ell+q+1}}{S_\ell}\leq1+\frac{x_k}{x_{k+1}}.
\end{equation}
To prove \eqref{eq:doubling-sliding-short}, write $P$ for the product. The indices occurring in its numerators are $\ell+q+1\geq r+q+1>k\geq u$, so $P$ is nonincreasing in each of $x_{r+1},\ldots,x_u$. If $u<k$, let these variables tend to zero. The resulting product telescopes to at most $S_u/x_k$, and hence $x_kP/S_u\leq1$. If $u=k$, let only $x_{r+1},\ldots,x_{u-1}$ tend to zero. The resulting product telescopes to at most $(S_u+x_k)/x_k$, and hence
\[
\frac{x_k}{S_u}P\leq1+\frac{x_k}{S_u}.
\]
Finally $S_u\geq x_{k+1}$, which proves \eqref{eq:doubling-sliding-short}.
Using \eqref{eq:doubling-residual-short} with the left neighbor $i=k-1$, the length factor is at most $2$, the level $j$ factor is at most $1$, and \eqref{eq:doubling-sliding-short} bounds the contribution of each level $2\leq q\leq j-1$ by
\[
1+\frac{x_k^{(q)}}{x_{k+1}^{(q)}}.
\]
Since $a_k\leq a_{k+1}$, the intervals $[t_{k-1},t_k]$ and $[t_k,t_k+a_k]$ are the two halves of one interval, and hence
\[
x_k^{(q)}\leq Kx_{k+1}^{(q)}.
\]
Thus every summand is bounded by $C_jK^{C_j}$. Summing over the $j+1$ nodes proves the theorem.
\end{proof}

\begin{corollary} \label{cor:projective-doubling-maximal}
Let $d\geq1$, $K\geq2$ and $\alpha>0$. Then
\[
\sup_{\mathcal C\in\mathrm{PDSECT}(d,K)}\Q_d(\mathcal C,N)\leq C_{d,\alpha,\varphi}K^{C_d}N^\alpha,\qquad N\geq1.
\]
\end{corollary}

\begin{proof} For $d=1$ the conclusion follows directly from the convention in Definition~\ref{def:doubling-canonical} and Theorem~\ref{thm:projective-SECT-maximal}. For $d\geq2$, apply Theorem~\ref{thm:doubling-canonical-switch} followed by Theorem~\ref{thm:projective-SECT-maximal}.
\end{proof}

\begin{example}The doubling condition includes many nonalgebraic canonical curves. For instance, any family of normalized weights satisfying
\[
0<c_q\leq w_q(t)\leq C_q<\infty
\]
gives a doubling canonical chart, with constant controlled by $\max_qC_q/c_q$. In dimension two,
\[
w_2(t)=\frac{\e^t}{\e-1}
\]
gives
\[
\mathbf X(t)=\left(1,t,\frac{\e^t-1-t}{\e-1}\right),
\]
whose image is not contained in an algebraic curve: a nonzero polynomial relation between the two affine coordinates would make $\e^t$ algebraic over $\R(t)$.

There are also doubling charts which do not belong to any power-unit class. Consider
\[
w_2(t)=\frac{1-\log t}{2},\qquad 0<t\leq1.
\]
This is normalized since $\int_0^1(1-\log t)\,\d t=2$. To see that $\mu_2=w_2\,\d t$ is doubling, write $J=[a,a+2h]$. Since $w_2$ is decreasing it suffices to compare the two halves. If $a\geq h$, then
\[
\frac{\mu_2([a,a+h])}{\mu_2([a+h,a+2h])}\leq\frac{1-\log a}{1-\log(a+2h)}\leq1+\log3.
\]
If $a<h$, then
\[
\mu_2([a,a+h])\leq h(2-\log(2h)),
\]
while
\[
\mu_2([a+h,a+2h])\geq
\begin{cases}
\frac h2(1-\log(3h)),&3h\leq1,
\\
\frac h2,&3h>1,
\end{cases}
\]
and the quotient is less than $5$. Thus $\mu_2$ is $6$-doubling. The corresponding curve is
\[
\mathbf X(t)=\left(1,t,\frac58t^2-\frac14t^2\log t\right)\in\mathrm{DSECT}(2,6).
\]
It belongs to no $\mathrm{PU}(2,\Lambda,R)$: for every $a>0$, the quotient $(1-\log t)t^{1-a}$ tends to $0$ if $a<1$ and to $\infty$ if $a\geq1$ as $t\to0^+$, and hence cannot be a unit with bounded multiplicative oscillation. The same holds for the projective arc: a projective change of frame, together with the corresponding canonical reparametrization, alters $w_2$ near $t=0$ only by factors bounded above and below, whereas the weights of a curve in $\mathrm{PU}(2,\Lambda,R)$ are of power-unit form near $t=0$ and, away from $t=0$, bounded above and below; so $[\mathbf X(I)]$ lies on no curve of $\mathrm{PPU}(2,\Lambda,R)$, for any $\Lambda$ and $R$. Finally, closed curves are covered by subadditivity. For instance, if $0<\varepsilon<(10\e)^{-1}$, the loop 
\[
\Gamma_\varepsilon=\left\{\left[1:\cos\theta:\sin\theta+\varepsilon\e^{\cos\theta}\right]:\, \theta\in\R/2\pi\Z\right\} 
\] 
is the union of four arcs which, in suitable affine charts, are graphs with second derivative bounded above and below by absolute constants. It is not algebraic: on its upper arc $y-\sqrt{1-x^2}=\varepsilon\e^x$, so an algebraic relation between $x$ and $y$ would make $\e^x$ algebraic over $\R(x)$. Each arc therefore lies in $\mathrm{PDSECT}(2,K_0)$ for an absolute constant $K_0$, and Corollary~\ref{cor:projective-doubling-maximal} gives $\Q_2(\Gamma_\varepsilon,N)\lesssim_{\alpha,\varphi}N^\alpha$.
\end{example}

\begin{remark} \label{rem:turning-vs-algebraic} The doubling condition provides a criterion of a different nature from the power-unit structure used in the algebraic argument. Even for algebraic examples which are already canonical curves in the sense of Section~\ref{sec:canonical-ect-overlap}, a direct use of the doubling constant can give substantially worse quantitative dependence. For the normalized weight measure
\[
\d\mu_2=(m+1)t^m\,\d t,\qquad m\in\{0,1,2,\ldots\},
\]
one has
\[
\frac{\mu_2([0,r])}{\mu_2([0,r/2])}=2^{m+1},
\]
whereas the corresponding generalized moment curve has $\mathsf S_2\lesssim(1+m)^8$ by Proposition~\ref{prop:muntz-switch}. Thus, on this family, a direct application of Corollary~\ref{cor:projective-doubling-maximal} gives constants exponential in $m$, whereas the power-unit analysis retains polynomial dependence on the order of vanishing. A finer decomposition of the parameter interval can reduce this loss, but we do not pursue this.
\end{remark}

Theorem~\ref{thm:doubling-canonical-switch} applies to curves already in canonical form, and its hypothesis is a condition on their weights. The next theorem gives a geometric condition on a curve in arbitrary position: doubling of its Frenet curvature measures $\kappa_i\,\d s$. This allows the curvatures to vanish or blow up at a point, like powers or logarithms, but not to vanish to infinite order, like $\e^{-1/s}$. The idea of the proof is to cut the curve into short arcs which turn little, in the sense that their total curvatures $\int\kappa_i\,\d s$ are small, and to choose on each arc the projective frame given by the Frenet frame at one of its points. The osculating subspaces along the arc then stay close to the coordinate subspaces of this frame, which are the osculating subspaces at the chosen point. In the resulting canonical chart, the transition measures, written in arclength, are comparable to the curvature measures, up to constant factors, so that doubling passes from the latter to the former.

\begin{theorem} \label{thm:doubling-curvatures} Let $d\geq2$, $K\geq2$ and $\alpha>0$. Let $\gamma:[0,\ell]\to {{\mathbb S}^d}$ be continuous and ${{\mathcal C}^d}$ on $(0,\ell)$, parametrized there by arclength, and assume that
\[
\gamma,\gamma',\ldots,\gamma^{(d)}
\]
are linearly independent on $(0,\ell)$, so that the Frenet curvatures $\kappa_1,\ldots,\kappa_{d-1}$ of Subsection~\ref{subsec:intro-doubling} are defined and positive there. Suppose that each measure $\kappa_i\,\d s$, $1\leq i\leq d-1$, is $K$-doubling. Writing
\[
M\coloneqq1+\ell+\sum_{i=1}^{d-1}\int_0^\ell \kappa_i (s)\,\d s,
\]
the image $\Gamma$ of $\gamma$ in $\Pd$ satisfies
\[
\Q_d(\Gamma,N)\leq C_{d,\alpha,\varphi}MK^{C_d}N^\alpha,\qquad N\geq1.
\]
\end{theorem}

\begin{proof} The case $\ell=0$ is the single-direction estimate, so let $\ell>0$. We first cut the curve into short arcs which turn little. Let
\[
\d\nu(s)\coloneqq \left(1+\sum_{i=1}^{d-1}\kappa_i(s)\right)\,\d s, \qquad \nu([0,\ell])=M-1.
\]
This measure is finite and nonatomic, and its cumulative mass is continuous and strictly increasing. Cutting at consecutive mass levels spaced by $10^{-2}$ partitions $[0,\ell]$ into at most $1+100(M-1)\leq100M$ nondegenerate closed intervals. On each interval $[a,b]$ both its length and every curvature mass are at most $10^{-2}$, and the restricted curvature measures are $K$-doubling on $[a,b]$.

We now show that each of these arcs belongs to $\mathrm{PDSECT}(d,CK^3)$. Fix one of the intervals $[a,b]$. Let $E_0,\ldots,E_d$ be the Frenet frame on $(0,\ell)$, let $c$ be the midpoint of $[a,b]$, put $e_i\coloneqq E_i(c)$ and
\[
\omega_j\coloneqq E_0\wedge\cdots\wedge E_j,\qquad a_j(s)\coloneqq\langle\omega_j(s),\omega_j(c)\rangle.
\]
By the Gram--Schmidt construction in Subsection~\ref{subsec:intro-doubling}, the curvatures $\kappa_1,\ldots,\kappa_{d-1}$ are continuous and positive on $(0,\ell)$. For $0\leq j<d$, the exterior frame $\omega_j$ is continuously differentiable in the interior, with $|\omega_j'|=\kappa_j$, where $\kappa_0\coloneqq1$; $\omega_d$ is constant. The stated ${{\mathcal C}^d}$ regularity suffices here, since differentiating $\omega_j$ uses derivatives of $\gamma$ only through order $j+1\leq d$. The finite curvature masses therefore give limits of these frames at the endpoints, and, for $s\in[a,b]$,
\[
|\omega_j(s)-\omega_j(c)|\leq10^{-2},\qquad a_j(s)=1-\tfrac12|\omega_j(s)-\omega_j(c)|^2\in[1-10^{-4},1].
\]
The differential identities below are first used in the interior; the extension of the canonical coordinates to the endpoints is justified below. Writing
\[
\Omega_j=\prod_{i=1}^{j-1}\kappa_i^{j-i}, \qquad V_j=\left\langle\gamma\wedge\cdots\wedge\gamma^{(j)},e_0\wedge\cdots\wedge e_j\right\rangle=\Omega_ja_j,
\]
let
\[
\tau(s)\coloneqq\frac{\langle\gamma(s),e_1\rangle}{\langle\gamma(s),e_0\rangle}.
\]
Then
\[
\frac{\d\tau}{\d s}=V_1V_0^{-2}=\frac{a_1}{a_0^2}>0.
\]

The functions $\langle\gamma,e_0\rangle,\ldots,\langle\gamma,e_d\rangle$ have Wronskians $V_0,\ldots,V_d$, which are positive. By \cite{KaSt}*{Chapter XI, Theorem 1.1, Remark 1.2 and Theorem 1.2}, after a lower-triangular change of coordinates which leaves $\langle\gamma,e_0\rangle$ and the Wronskians unchanged, they are iterated integrals of positive weights, as in \cite{KaSt}*{Chapter XI, (1.5)}, and the weights are given by \cite{KaSt}*{Chapter XI, (1.13)}:
\[
w_0=\langle\gamma,e_0\rangle,\qquad w_1=\frac{V_1}{V_0^2}=\frac{\d\tau}{\d s},\qquad w_q=\frac{V_qV_{q-2}}{V_{q-1}^2},\quad 2\leq q\leq d.
\]
Dividing by $w_0$ and using $\tau$ as parameter, we obtain a canonical curve whose transition measures satisfy
\[
\tau^\ast\mu_q=c_q\frac{V_qV_{q-2}}{V_{q-1}^2}\,\d s,\qquad 2\leq q\leq d,
\]
for constants $c_q>0$ coming from the normalization \eqref{e:normweights}. Since
\[
\frac{\Omega_q\Omega_{q-2}}{\Omega_{q-1}^2}=\kappa_{q-1},
\]
we obtain
\[
\tau^\ast\mu_q=c_q\kappa_{q-1}\frac{a_qa_{q-2}}{a_{q-1}^2}\,\d s.
\]
Both $\d\tau/\d s$ and $a_qa_{q-2}a_{q-1}^{-2}$ lie in $[1-10^{-3},1+10^{-3}]$. If $J'$ is either half of a $\tau$-interval $J\subset\tau([a,b])$, then
\[
\frac{|\tau^{-1}(J')|}{|\tau^{-1}(J)|}>\frac14.
\]
Lemma~\ref{lem:iterated-doubling}, with three iterations, therefore gives
\[
\mu_q(J)\leq CK^3\mu_q(J').
\]
Thus every transition measure is $CK^3$-doubling.

The representation of \cite{KaSt}*{Chapter XI, Theorem 1.2} is applied on compact subintervals of $(a,b)$, on which the weights are given by the same formulas. To include the endpoints, first reconstruct the canonical basis with the interior base point $\tau(c)$. The preceding formula shows that each reconstructed weight is positive and integrable on the full closed $\tau$-interval $\tau([a,b])$, even if a curvature vanishes or is unbounded at an endpoint. Every iterated Volterra integral based at $\tau(c)$ therefore extends continuously to that interval. The identities expressing the affine coordinates of $\gamma$ in this basis hold in the interior and extend to both endpoints by continuity. Also $a_0$ stays positive and $\d\tau/\d s$ stays bounded above and below, so the affine chart and the parameter remain valid on the closed arc.

Changing the base point of these iterated integrals to $\tau(a)$ is a constant invertible lower-triangular change of coordinates preserving each prefix space. Rescale the parameter interval to $I=[0,1]$ and normalize the weight masses. The resulting chart belongs to $\mathrm{DSECT}(d,CK^3)$, as these operations preserve doubling constants. Hence the closed projective arc $[\gamma([a,b])]$ belongs to $\mathrm{PDSECT}(d,CK^3)$.

Finally, by Corollary~\ref{cor:projective-doubling-maximal}, each of the at most $100M$ arcs satisfies
\[
\Q_d\left([\gamma([a,b])],N\right)\leq C_{d,\alpha,\varphi}K^{C_d}N^\alpha,
\]
and the subadditivity of $\Q_d$ gives the result.
\end{proof}

The Frenet curvatures, and hence the hypothesis of Theorem~\ref{thm:doubling-curvatures}, depend on the choice of projective frame, through the spherical normalization of the curve. The next proposition shows that the doubling condition does not: a projective change of coordinates $L\in\mathrm{GL}_{d+1}(\R)$ changes the doubling constants only by factors depending on $\|L\|\|L^{-1}\|$.

\begin{proposition}  \label{prop:projective-curvature-doubling} Under the hypotheses of Theorem~\ref{thm:doubling-curvatures}, let $L\in\mathrm{GL}_{d+1}(\R)$ and set \[\widetilde\gamma(s)\coloneqq \frac{L\gamma(s)}{|L\gamma(s)|}.\] After reparametrization by its arclength $\sigma$, the curve $\widetilde\gamma$ again has finite doubling Frenet curvature measures. More precisely, if $\chi=\|L\|\|L^{-1}\|$, a common doubling constant is
\[
 \widetilde K=\chi^{4d}K^{\lceil\log_2(4\chi^4)\rceil}.
\]
Consequently, the existence of a finite doubling constant does not depend on the projective frame: it is a property of the projective arc.
\end{proposition}

\begin{proof} All derivatives in this proof are initially taken with respect to $s$. Put
\[
W_j=\gamma\wedge\gamma'\wedge\cdots\wedge\gamma^{(j)}, \qquad \Delta_j=|W_j|,\qquad 0\leq j\leq d,
\]
and define $\widetilde W_j,\widetilde\Delta_j$ in the same way. Since $\widetilde\gamma=|L\gamma|^{-1}L\gamma$, the Leibniz rule shows that $\widetilde\gamma^{(k)}$ is $|L\gamma|^{-1}L\gamma^{(k)}$ plus a linear combination of $L\gamma,\ldots,L\gamma^{(k-1)}$; these lower-order terms do not contribute to the wedge product, and hence
\[
\widetilde W_j=|L\gamma|^{-(j+1)}\left(\mathop{\bigwedge}\nolimits^{j+1}L\right)W_j.
\]
In particular, nondegeneracy is preserved. The Frenet identities give, for $2\leq q\leq d$,
\[
\kappa_{q-1}\,\d s =\frac{\Delta_q\Delta_{q-2}}{\Delta_{q-1}^2}\,\d s,\qquad\sigma^\ast\left(\widetilde\kappa_{q-1}\,\d\sigma\right) =\frac{\widetilde\Delta_q\widetilde\Delta_{q-2}} {\widetilde\Delta_{q-1}^2}\,\d s.
\]
The second identity includes the change of arclength: writing $\widetilde\gamma(s)$ as its arclength parametrization composed with $\sigma(s)$ contributes a factor $(\d\sigma/\d s)^{j(j+1)/2}$ to a wedge of order $j$, so the displayed ratio acquires exactly one factor $\d\sigma/\d s$. The powers of $|L\gamma|$ cancel in the ratio. Bounding the exterior powers of $L$ by its largest and smallest singular values therefore gives
\begin{equation}\label{eq:projective-curvature-comparison}
\chi^{-2q}\kappa_{q-1}\,\d s \leq \sigma^\ast\left(\widetilde\kappa_{q-1}\,\d\sigma\right) \leq \chi^{2q}\kappa_{q-1}\,\d s.
\end{equation}
Also, since $\Delta_1=1$,
\[
\chi^{-2}\leq\frac{\d\sigma}{\d s} =\frac{|(\mathop{\bigwedge}\nolimits^2 L)W_1|}{|L\gamma|^2} \leq\chi^2.
\]
Thus the new length and all curvature masses are finite. If $J'$ is either half of a $\sigma$-interval $J$, its preimage has length at least $|\sigma^{-1}(J)|/(2\chi^4)$. Apply Lemma~\ref{lem:iterated-doubling} with $m=\lceil\log_2(4\chi^4)\rceil$ and then \eqref{eq:projective-curvature-comparison}. Each new curvature measure is $\chi^{4q}K^m$-doubling, hence $\widetilde K$-doubling. Applying the same argument to $L^{-1}$ proves the converse.
\end{proof}

\begin{example}\label{ex:rotation-orbit} Let $\lambda\in\R\setminus\{0,\pm1\}$ and $\theta_0>0$, and consider the curve in $\mathbb S^3$
\[
\gamma(\theta)=\frac1{\sqrt2}\left(\cos\theta,\ \sin\theta,\ \cos\lambda\theta,\ \sin\lambda\theta\right),\qquad 0\leq\theta\leq\theta_0.
\]
It is the orbit $\gamma(\theta)=\exp(\theta A)\gamma(0)$ of a one-parameter group of rotations, with $A$ skew-symmetric. Hence $\gamma^{(k)}(\theta)=\exp(\theta A)A^k\gamma(0)$, so $|\gamma'|=\sqrt{(1+\lambda^2)/2}$ is constant, and the Gram--Schmidt frame of $\gamma,\gamma',\gamma'',\gamma'''$ at $\theta$ is the image under the rotation $\exp(\theta A)$ of the frame at $0$. After reparametrization by arclength, the Frenet curvatures $\kappa_1,\kappa_2$ are therefore constant. They are positive: $\gamma(0),A\gamma(0),A^2\gamma(0),A^3\gamma(0)$ are linearly independent, since their determinant is a nonzero multiple of $\lambda(1-\lambda^2)^2$. The curvature measures are multiples of arclength, hence $2$-doubling, and Theorem~\ref{thm:doubling-curvatures} gives
\[
\Q_3\left([\gamma([0,\theta_0])],N\right)\leq C_{\alpha,\varphi,\lambda}(1+\theta_0)N^\alpha,\qquad N\geq1.
\]
If $\lambda$ is irrational, no algebraic curve contains $\gamma([0,\theta_0])$, nor its projectivization. Indeed, a polynomial vanishing on this arc vanishes on the whole orbit, by analyticity in $\theta$, hence on its closure, which is the two-dimensional surface $\{x_1^2+x_2^2=x_3^2+x_4^2=\tfrac12\}$. The curve is given in natural coordinates, far from canonical form, and its canonical weights are not explicit, while the curvature hypothesis is verified by symmetry alone.
\end{example}

\section{Applications}\label{sec:applications} This section derives Corollaries~\ref{cor:intro-arbitrary-3d}--\ref{cor:intro-CZ-3d} from Theorem~\ref{thm:intro-main}, using in addition \cite{DP-PAdv}*{Theorem E} for Corollary~\ref{cor:intro-arbitrary-3d}, and the exponential good-$\lambda$ comparison for Corollaries~\ref{cor:intro-CZ-algebraic} and~\ref{cor:intro-CZ-3d}.

\subsection{Arbitrary directions in three dimensions} Let
\[
F(N)\coloneqq \sup_{\substack{V\subset {\mathbb S^2}\\\#V\leq N}}\|\A_V\|_{2\to2}, \qquad N\geq1.
\]
By \eqref{eq:homogeneity}, normalizing the vectors of $V$ does not change $\A_V$, so $F(N)$ is the left-hand side of \eqref{eq:intro-arbitrary-3d}. The proof is a recursion in the number of directions, based on the algebraic almost orthogonality principle of \cite{DP-PAdv}*{Theorem E}. We first show that there are constants $A,c_0,C_0<\infty$, independent of $\eta>0$, such that
\begin{equation}\label{eq:ambient-recursion}
F(N)\leq C_\eta D^A N^\eta+C_0D^{1/2}F \left(\left\lceil c_0N/D^2\right\rceil\right)
\end{equation}
for integers $2\leq D\leq N^{1/2}$.

We use an exact polynomial partition on ${\mathbb S^2}$, whose construction we recall for completeness. The restrictions to ${\mathbb S^2}$ of polynomials of degree at most $k$ form a space of dimension $(k+1)^2$. Successively applying the polynomial ham-sandwich theorem in these spaces, with degrees of order $2^{j/2}$ at the $j$th bisection, produces a polynomial $p$ of degree $O(D)$, not identically zero on ${\mathbb S^2}$, such that each component of ${\mathbb S^2}\setminus Z(p)$ contains at most $c_0N/D^2$ of the selected directions. Indeed, after $J$ bisections each sign region contains at most $N/2^J$ points, whereas the degree of the product is $O(2^{J/2})$; choose $2^{J/2}$ comparable to $D$. Points lying on $Z(p)$ are assigned to the wall. This standard argument allows all scales to be treated without approximating directions near a wall.

Apply \cite{DP-PAdv}*{Theorem E} with the set of scales $S\coloneqq(0,\infty)$, that is, to averages of all lengths. That theorem is stated for the smooth averages of \cite{DP-PAdv}, defined by a fixed real even multiplier with small support, without positivity of its inverse Fourier transform; its proof uses only this support and the decay of the kernel, and therefore applies to the multiplier $c\varphi(2^8\cdot)$, $c>0$, which defines the maximal operator $c\A_V$. The wall ${\mathbb S^2}\cap Z(p)$ is the set of real points of the conjugation-invariant curve $\{x\in\C^3: x_1^2+x_2^2+x_3^2=1,\ p(x)=0\}$, whose projective closure has degree $O(D)$ by B\'ezout's theorem and no real points at infinity. By Lemma~\ref{lem:real-points-degree}, the wall has dimension at most one and degree $O(D^3)$, isolated real points included. Hence Theorem~\ref{thm:intro-main} bounds its contribution by $C_\eta D^A N^\eta$. The cellular contribution is the second term of \eqref{eq:ambient-recursion}. This proves that inequality.

\begin{proof}[Proof of Corollary~\ref{cor:intro-arbitrary-3d}]
Fix $\eta=1/16$ and $\delta=1/(16(A+1))$. For sufficiently large $N$ take $D=\lfloor N^\delta\rfloor$ and $m=\lceil c_0N/D^2\rceil$. Then
\[
D^A N^\eta\leq N^{1/8},\qquad m\leq C N^{1-2\delta},\qquad D^{1/2}m^{1/4}\leq C N^{1/4}.
\]
Increasing a fixed threshold $N_0$ gives $m<N$ and
\[
\log(2+m)\leq(1-\delta)\log(2+N),\qquad N\geq N_0.
\]
Choose $B$ so large that the constant multiplying the cellular term in \eqref{eq:ambient-recursion}, including the constants in the preceding display, times $(1-\delta)^B$ is at most $1/2$. Induction on $N$, applied to $F(N)\leq C_*N^{1/4}(\log(2+N))^B$, now closes: the cellular term is at most half the proposed bound, and the wall term is at most the other half after choosing $C_*$ large enough. The finitely many cases $N<N_0$ follow from the elementary bound $F(N)\lesssim N^{1/2}$.
\end{proof}

\begin{remark}\label{rmrk:higherd} The argument in the proof of Corollary~\ref{cor:intro-arbitrary-3d} applies in every ambient dimension $n$ up to its last step: polynomial partitioning and \cite{DP-PAdv}*{Theorem E} reduce matters to the directions on the wall $Z$, the zero set in $\mathbb S^{n-1}$ of the partitioning polynomial, which is an algebraic subvariety of dimension at most $n-2$. For $n=3$ the variety $Z$ is a curve and Theorem~\ref{thm:intro-main} completes the proof. For $n\geq4$, a natural way to continue would be to perform a polynomial partition of $Z$ itself, but the overlap argument does not carry over. In the proof of \cite{DP-PAdv}*{Theorem E}, the cells are open in $\mathbb S^{n-1}$: each cell contains a surface ball around each of its directions whose radius is a small multiple of the distance of that direction from the wall. A hyperplane passing close to a direction, relative to this distance, therefore meets the cell of that direction, so the overlap of the frequency bands of the cells is bounded by the number of cells that a hyperplane meets. The cells of a partition of $Z$ are not open in $\mathbb S^{n-1}$ and contain no such balls; such a cell may lie in a thin band around a hyperplane that does not meet it, and counting crossings no longer controls the overlap. The same phenomenon occurs for curves, and dealing with it is one of the main difficulties of this paper: the finite-crossing property alone does not suffice, and we rely instead on a quantitative form of it, encoded in the uniform boundedness of the switch moduli of Section~\ref{sec:canonical-ect-overlap}. Extending this approach to varieties of dimension two or higher would require both a suitable quantitative crossing property for such varieties and an effective algebraic preparation of the variety into canonical-like charts, analogous to the one carried out for curves in Section~\ref{sec:algebraic}, each of which is a challenging problem in its own right.
\end{remark}

\subsection{Local maximal averages and Sobolev regularity} We use a finite-frequency discretization which does not require any regularity of the vector field. Let $\lambda\geq2$ and let $g\in\mathcal S(\R^n)$ satisfy $\supp\widehat g\subset B(0,\lambda)$. Choose a $\lambda^{-1}$-net $E_\lambda\subset\Gamma$ with $\#E_\lambda\leq C_n\lambda^{n-1}$; the elementary packing bound on the ambient sphere suffices. Let
\[
\zeta_\lambda(x)\coloneqq\lambda^n(1+\lambda|x|)^{-n-1}, \qquad h_\lambda\coloneqq|g|*\zeta_\lambda.
\]
Then
\begin{equation}\label{eq:local-net-domination}
\M_\Gamma^{\mathrm{loc}}g \leq C_n\M_{E_\lambda}h_\lambda, \qquad \|h_\lambda\|_2\leq C_n\|g\|_2.
\end{equation}
To see this, write $g=g*\rho_\lambda$, where $\rho_\lambda\coloneqq\lambda^n\rho(\lambda\cdot)$ and $\rho$ is a Schwartz kernel whose Fourier transform equals one on the unit ball. For $|v-w|\leq\lambda^{-1}$ and $|t|\leq1$, rapid decay gives
\[
|g(x-tv)|\leq C_nh_\lambda(x-tw),
\]
because $|t(v-w)|\leq\lambda^{-1}$. Average over $[-r,r]$, with $r\leq1$, and take suprema. Young's inequality proves the norm bound.

Theorem~\ref{thm:intro-main}, applied to $E_\lambda$, now gives, for every $\sigma>0$,
\begin{equation}\label{eq:local-bandlimited-bound}
\|\M_\Gamma^{\mathrm{loc}}g\|_2 \leq C_{n,\sigma}D^{C_n}\lambda^\sigma\|g\|_2.
\end{equation}
Here choose the exponent in the theorem to be $\sigma/(n-1)$.

\begin{proof}[Proof of Corollary~\ref{cor:intro-sobolev}] Let $f=\sum_{j\geq0}P_jf$ be a smooth inhomogeneous Littlewood--Paley decomposition, with $P_jf$ supported in frequencies $|\xi|\lesssim2^j$. Apply \eqref{eq:local-bandlimited-bound} with $\sigma=s/2$, then use sublinearity and Cauchy--Schwarz:
\[
\begin{split}
\|\M_\Gamma^{\mathrm{loc}}f\|_2 &\leq C_{n,s}D^{C_n}\sum_{j\geq0}2^{js/2}\|P_jf\|_2
\\
&\leq C_{n,s}D^{C_n} \left(\sum_{j\geq0}2^{2js}\|P_jf\|_2^2\right)^{1/2} \leq C_{n,s}D^{C_n}\|f\|_{H^s}.
\end{split}
\]
For a measurable field $v$ taking values in $\Gamma$, $\M_{v(\cdot)}^{\mathrm{loc}}f\leq\M_\Gamma^{\mathrm{loc}}f$ on Schwartz functions; for such $f$ the average $\frac1{2r}\int_{-r}^r|f(x-tu)|\,\d t$ is continuous in $(x,u,r)$, so the supremum may be taken over rational $r$, and $\M_{v(\cdot)}^{\mathrm{loc}}f$, a countable supremum of continuous functions composed with the measurable map $x\mapsto(x,v(x))$, is measurable. The inequality $|\M_{v(\cdot)}^{\mathrm{loc}}f-\M_{v(\cdot)}^{\mathrm{loc}}g| \leq\M_{v(\cdot)}^{\mathrm{loc}}(f-g)$ then gives a unique continuous extension to $H^s$. The same approximation controls the directional averages uniformly in $0<r\leq1$. For smooth functions those averages
converge to the function as $r\downarrow0$; the maximal estimate and Chebyshev's inequality pass this conclusion to their density extensions for $H^s$ functions, with $r$ restricted to a countable dense subset of $(0,1]$. This proves the differentiation assertion following the corollary.
\end{proof}

\subsection{Directional Calder\'on--Zygmund operators} \label{subsec:directional-CZ} We record the comparison underlying the two singular-integral corollaries. For every nonempty finite $V\subset\R^n\setminus\{0\}$, with $\#V\leq N$,
\begin{equation}\label{eq:CZ-maximal-comparison}
\|T_V\|_{L^2(\R^n)\to L^2(\R^n)} \leq C_T\log(2+N)\|\M_V\|_{L^2(\R^n)\to L^2(\R^n)}.
\end{equation}
Here $\M_V$ has the same definition for nonzero vectors of arbitrary length, since its supremum over averaging lengths absorbs their norms.

Following \cite{DemeterDiPlinio2014}*{Proposition 2.2 and \S2.2}, we use a fiberwise argument. Let $M$ be the one-dimensional Hardy--Littlewood maximal operator and $T_*$ the maximal truncation of the singular integral $T$. The operator $\mathcal R_Tg=M(Tg)+T_*g$ is bounded on $L^2(\R)$, and
\begin{equation}\label{eq:CZ-good-lambda}
\left|\{|Tg|>a\lambda,\ Mg\leq c_T\lambda\}\right| \leq C_T \e^{-b_Ta} \left|\{\mathcal R_Tg>\lambda\}\right|, \qquad a\geq a_T,\quad \lambda>0,
\end{equation}
with constants depending only on the stated data of $T$. The proof of \cite{DemeterDiPlinio2014}*{Proposition 2.2} is written for convolution kernels; combined with the local exponential estimate \cite{OrtizPerezRela2013}*{Theorem 3}, it gives \eqref{eq:CZ-good-lambda} for general Calder\'on--Zygmund kernels. Apply \eqref{eq:CZ-good-lambda} on the lines $y+\R v$, use Fubini, and take the union over $v\in V$. Denoting the corresponding fiber operator by $\mathcal R_{T,v}$ gives
\[
|\{T_Vf>a\lambda\}| \leq |\{\M_Vf>c_T\lambda\}|+C_T\e^{-b_Ta}\sum_{v\in V}|\{\mathcal R_{T,v}f>\lambda\}|.
\]
The coordinate change $x=y+tv$ shows that $\|\mathcal R_{T,v}f\|_2\leq C_T\|f\|_2$, uniformly in $v$. Integrating against $2\lambda\,\d\lambda$ therefore yields
\[
a^{-2}\|T_Vf\|_2^2\leq C_T\|\M_Vf\|_2^2+C_TN \e^{-b_Ta}\|f\|_2^2.
\]
Take $a=C_T'\log(2+N)$ sufficiently large. Since $V\neq\varnothing$ implies $|f|\leq\M_Vf$ almost everywhere, this proves \eqref{eq:CZ-maximal-comparison}, first for smooth compactly supported functions and then by density.

\begin{proof}[Proof of Corollaries~\ref{cor:intro-CZ-algebraic} and~\ref{cor:intro-CZ-3d}] Use \eqref{eq:CZ-maximal-comparison} and the equivalence of the positive and smooth maximal models. Theorem~\ref{thm:intro-main}, with exponent $\eps/2$, gives \eqref{eq:intro-CZ-algebraic} after absorbing the logarithm; \eqref{eq:intro-subpower} gives its stronger subpower form. For $\R^3$, Corollary~\ref{cor:intro-arbitrary-3d} gives the bound with $(\log(2+N))^{B+1}$, hence \eqref{eq:intro-CZ-3d}. Here radial normalization of $V$ leaves the positive comparison operator unchanged; the definition of $T_v$ in \eqref{eq:intro-directional-CZ} stays fixed.
\end{proof}

\section*{Acknowledgment}\sloppy Part of this project was carried out during F. Di Plinio's visit to the Basque Center for Applied Mathematics (BCAM) in Bilbao in March--April 2025 and I. Parissis's stay at Università degli Studi di Napoli Federico II in February--March 2026. The former visit was partially supported by grant CEX2021-001142-S funded by {MICIU/AEI/10.13039/501100011033}. The latter stay was supported by the Program for Visiting Professors of Istituto Nazionale di Alta Matematica ``F. Severi'' (INdAM). The authors gratefully acknowledge the support and hospitality of all institutions involved.

\section*{Statements and declarations}

\subsection*{Use of AI tools}
The main argument of this paper was laid out by the authors between March 2025 and March 2026, in a form based on Frenet frames, and presented in talks, including at the Analysis and Applications Seminar at UAM on 21 November 2025 \cite{ParICMAT25} and the BIRS--IMAG workshop in Granada on 10 June 2026 \cite{DPBIRS26}. The present version was developed to explicitly rely on projective invariance, to obtain bounds polynomial in the degree, and to correct inaccuracies in the earlier decomposition of algebraic curves into basic pieces. In preparing it, we used the AI assistant Claude (Anthropic) to streamline, check and improve definitions, notation and proofs throughout. This applies in particular to the decomposition of algebraic curves in Section~\ref{sec:algebraic}, whose mechanism, the reduction to power-unit charts through the cellular preparation of Binyamini and Novikov, is fully due to the authors. We also used it to verify references and computations, and to edit the exposition. The authors have independently checked all mathematical content and take full responsibility for it.

\subsection*{Competing interests}
The authors have no competing interests to declare.

\bibliographystyle{amsplain}
\bibliography{conical}

\end{document}